 \pdfoutput=1
\documentclass[final,executivepaper,leqno,onefignum,onetabnum]{siamltex}
\usepackage{graphicx}
\usepackage[notref,notcite]{showkeys}
\usepackage{amssymb,amsmath,amsbsy,mathrsfs,eucal,amsfonts}
\usepackage{lscape}
\usepackage{fancybox}
\usepackage{algorithmic}
\usepackage{epstopdf}
\usepackage{color}
\usepackage{float}
\numberwithin{equation}{section}
\usepackage{graphicx}
\usepackage[notref,notcite]{showkeys}
\usepackage{amssymb,amsmath,amsbsy,mathrsfs,eucal,amsfonts}
\usepackage{lscape}
\usepackage{fancybox}
\usepackage{algorithmic}
\usepackage{epstopdf}
\usepackage{color}
\numberwithin{equation}{section}
\newtheorem{example}[theorem]{Example}
\newtheorem{remark}{Remark}[section]
\newcommand{\bld}[1]{\boldsymbol{#1}}
\newcommand{\bx}{\bld{x}}
\newcommand{\bn}{\bld{n}}
\newcommand{\qn}{\boldsymbol{q}}

\newcommand{\vn}{\boldsymbol{v}}
\newcommand{\normal}{\boldsymbol{n}}
\newcommand{\Vn}{\boldsymbol{V}}
\newcommand{\mn}{\boldsymbol{m}}

\newcommand{\Kn}{\textbf{\textrm{K}}}
\newcommand{\In}{\textbf{\textrm{I}}}

\newcommand{\divv}{\nabla \cdot}
\newcommand{\trian}{\textsf{D}_h}

\newcommand{\ltwoint}{{\mbox{\tiny int}}}
\newcommand{\ltwotilde}{{L^2(\widetilde{\textsf{D}}_h)}}

\newcommand{\tildeE}{{L^2(\widetilde{\mathcal{E}}_h)}}

\definecolor{blue}{rgb}{0,0,1}

\usepackage{float}

\newcommand{\xx}{\bld{x}}

\newcommand{\divh}{\nabla_{h} \cdot}

\title{A high order HDG method for curved-interface problems via approximations from straight triangulations}

\author{ Weifeng Qiu 
        \thanks{Department of Mathematics, City University of Hong Kong, Hong Kong, email: {\tt weifeqiu@cityu.edu.hk}.}
\and
       Manuel Solano
        \thanks{Corresponding author. Departamento de Ingenier\'ia Matem\'atica 
and CI$^2$MA, Universidad de Concepci\'on, Chile, email: {\tt msolano@ing-mat.udec.cl}.   }
\and 
Patrick Vega
 \thanks{Departamento de Ingenier\'ia Matem\'atica 
and CI$^2$MA, Universidad de Concepci\'on, Chile, email: {\tt pvega@ing-mat.udec.cl}.}
}

\begin{document}

\maketitle

\begin{abstract}
We generalize the technique of [{\em Solving Dirichlet boundary-value problems on curved domains 
by extensions from subdomains},  SIAM J. Sci. Comput. 34, pp. A497--A519 (2012)] to elliptic problems with mixed boundary conditions and elliptic interface problems involving a non-polygonal interface. 
We study first the treatment of the Neumann boundary data since it is crucial to understand the applicability of the technique to curved interfaces. We provide numerical results showing that, 
in order to obtain optimal high order convergence, it is desirable 
to construct the computational domain by interpolating the boundary/interface  using piecewise linear segments. In this case the distance of the computational domain to the exact boundary is only $O(h^2)$.
\end{abstract}



\section{Introduction}
\label{Introduction}
In this paper we present a technique to numerically solve second order elliptic problems in domains $\Omega$ which are 
not necessarily polygonal. In addition, we deal with domains divided in two regions by a curved interface $\Sigma$. In particular we use a high order hybridizable discontinuous Galerkin method (HDG) \cite{CGL,CGS} 
where the computational domain do not exactly fit the curved boundary or interface. The main motivation of this technique is being able 
to use high order polynomial approximations and keep high order accuracy using triangular meshes having only straight elements.

One of the first ideas in this direction  
was introduced by \cite{CGR} for the one-dimensional case and then extended to higher space dimensions for pure diffusion  
\cite{CQS,CS} and convection-diffusion \cite{CS} equations. In their work, the mesh does not fit the domain and the distance between the computational domain and the boundary $\Gamma := \partial \Omega$ is of only  order $O(h)$, making this method attractive from a computational point of view. In addition, \cite{CSS} applied this method to couple boundary element 
and HDG methods to solve exterior diffusion problems. However, only  Dirichlet boundary value problems have been considered since  Neumann data can not be handled in the same way as we will explain below. We will see that for the Neumann boundary case the proposed technique works properly if the computational domain is order $O(h^2)$ away from the actual boundary. 

The work presented here focuses first on the  treatment of part of the boundary where a Neumann data is prescribed. It is important to understand this situation in order to extend the ideas to problems having a curved interface. In fact, the transmission conditions at the interface involve jumps of the scalar variable and jumps of the normal component of the flux. The computational jump of the scalar variable can be treated considering the {\it transferring} technique of \cite{CS} and the computational jump of the normal component of the flux can be handled using the {\it extrapolation} method for the Neumann data that we will describe in the following sections.\\

One of the first methods that approximate Neumann boundary conditions on curved domains considering non-fitted meshes  was introduced by \cite{BE3}. Here, a piecewise linear finite element method was considered and optimal convergence in the $H^1$-norm was shown. In addition, the same authors solved a semi-definite Neumann problem on curved domains using a similar technique (\cite{BE1}). They showed optimal behavior of the errors in $H^1$ and $L^2$-norms using again piecewise linear elements. On the other hand, higher order approximation finite element methods require to properly fit the boundary in order to keep  high order accuracy. For instance, isoparametric element can be considered (\cite{BE1},\cite{Lenoir}). In the case of elliptic interface problems, usually the curve describing the interface is interpolated by a piecewise linear computational interface. Hence, super-parametric elements near the interface must be considered in order to achieve high order accuracy (\cite{MIT1}).

This article aims to develop a high order method based on a triangulation of the domain involving only straight elements. As we will discuss, the boundary/interface must be interpolated by piecewise linear function in order to obtain the expected rates of convergence. Since most of the methods based on linear fitting are only  second order accurate, we believe our method constitutes a competitive alternative.

The rest of the manuscript is organized  as follows. We will begin by setting notation. Then, we will describe the technique for a boundary-value problem where Neumann data is prescribed in part of the boundary. In particular, we will discuss the proper choice of the {\it paths} that will {\it transfer} the Dirichlet and impose the Neumann data. We will provide numerical simulations showing the performance of the method.  Then, we will adapt these ideas in order to solve a elliptic interface problem and show numerical experiments validating the technique.

%
%

\section{Mesh construction and notation}
\label{sec:mesh_notation}
Let $\textsf{D}_h$ be a a triangulation constructed by the union of disjoint straight triangles that approximates a bounded domain $\Omega\subset \mathbb{R}^2$ and does not necessarily fit its boundary. The Dirichlet and Neumann part of the boundary $\Gamma$ are denoted by $\Gamma_D$ and $\Gamma_N$ ($\Gamma_D \cap \Gamma_N =\emptyset$, $\Gamma_D \cup \Gamma_N = \Gamma$). We also assume that the computational boundary, $\Gamma^h$, satisfies $\Gamma^h = \Gamma_D^h \cup \Gamma_N^h$ and $\Gamma_D^h \cap \Gamma_N^h =\emptyset$ where $\Gamma_D^h$ and $\Gamma_N^h$ are part of $\Gamma^h$ with Dirichlet ($\widetilde{g}_D$) and  Neumman ($\widetilde{g}_N$) data, respectively. Let $d(\Gamma,\Gamma^h)$ be the distance between $\Gamma$ and $\Gamma^h$. We denote by $h_K$ the diameter of the element $K\in\textsf{D}_h$ and
by $\normal$ its outward unit normal.  The meshsize $h$ is defined as $\max_{K\in \textsf{D}_h} h_K$. Let $\mathcal{E}^0_h$
be the set of interior edges of $\trian$ and
$\mathcal{E}^{\partial}_h$ the edges at the boundary. We say that
an edge $e\in \mathcal{E}^{0}_h$ if there are two elements $K^+$
and $K^-$ in $\trian$ such that $e=\partial K^+ \cap \partial
K^-$. Also, we say that $e\in \mathcal{E}^{\partial}_h$ if there
is an element $K \in \trian$ such that $e=\partial K\cap
\Gamma^h$. We set $\mathcal{E}_h=\mathcal{E}^0_h
\cup\mathcal{E}^{\partial}_h$. For each element $K$ in the triangulation $\trian$, we
denote by $\mathcal{P}^{k}(K)$ the space of polynomials of
degree at most $k$ defined on the element $K$. For each edge $e$ in $\mathcal{E}_h$
$\mathcal{P}^{k}(e)$ is the space of polynomials of degree at most
$k$ defined on the edge $e$. Given an element $K$,
$(\cdot,\cdot)_K$ and $\langle\cdot, \cdot\rangle_{\partial K}$
denote the $L^2(K)=\{v: \int_K v^2 < \infty\}$ and $L^2(\partial
K)=\{\xi: \int_{\partial K} \xi^2 < \infty\}$ products,
respectively. Thus, for each $\xi$ and $\psi$ we def\/ine
\begin{eqnarray*}
(\xi,\psi)_{\trian}=\sum_{K\in\trian} (\xi,\psi)_K  \quad
\textrm{and} \quad
\langle\xi,\psi\rangle_{\partial\trian}=\sum_{K\in\trian}
\langle\xi,\psi\rangle_{\partial K}.
\end{eqnarray*}

\section{Boundary value problem with mixed boundary conditions}
\label{sec:mbc}

We consider the following model problem:
\begin{subequations}\label{model_problem}
\begin{eqnarray}
- \divv\: \qn &=&f \:\: \textrm{in} \:\: \Omega,\label{problem1} \\
\qn+\Kn \nabla u&=& 0\:\: \: \textrm{in} \:\: \Omega, \label{problem2}\\
u&=&g_D \:\: \: \textrm{on} \:\: \Gamma_D,\label{problem3}\\
\qn\cdot \normal&=&g_N \:\: \: \textrm{on} \:\: \Gamma_N.\label{problem4}
\end{eqnarray}
\end{subequations}
Here $g_D\in H^{1/2}(\Gamma_D)$ and  $g_N\in H^{-1/2}(\Gamma_N)$ are given data at the border, $f\in L^2(\Omega)$ is a
source term and $\Kn\in [L^\infty(\Omega)]^{2\times 2}$ is a symmetric and positive definite tensor.   

In the computational domain  $\textsf{D}_h$, problem   (\ref{model_problem}) can be written as follows:
\begin{subequations}\label{model_problem_h}
\begin{eqnarray}
- \divv\: \qn &=&f \:\: \textrm{in} \:\: \textsf{D}_h,\label{problem1_h} \\
\qn+\Kn \nabla u&=& 0\:\: \: \textrm{in} \:\: \textsf{D}_h, \label{problem2_h}\\
u&=&\widetilde{g}_D \:\: \: \textrm{on} \:\: \Gamma_D^h,\label{problem3_h}\\
\qn\cdot \normal&=&\widetilde{g}_N \:\: \: \textrm{on} \:\: \Gamma_N^h.\label{problem4_h}
\end{eqnarray}
\end{subequations}

Here $\widetilde{g}_D$ and $\widetilde{g}_N$ are unknowns. As we mentioned before, $\widetilde{g}_D$ can be calculated following \cite{CGR,CSS,CS}, i.e.,
\begin{eqnarray}
\widetilde{g}_D(\bx): = g_D(\bar{\bx}) + \int_{\sigma(\bx)} \Kn^{-1} \qn
\cdot \mn \:ds, \label{new_bc}
\end{eqnarray}
where $\sigma(\bx)$, is a path starting at $\bx \in \Gamma_D^h$ and ending at $\bar{\bx}\in \Gamma_D$; and $\mn$ is the
tangent vector to $\sigma(\bx)$. This expression comes from integrating (\ref{problem2}) along the path $\sigma(\bx)$ (see \cite{CS} for details). 

In principle, any kind of numerical method using polygonal domains can be used to solve the equations in $\textsf{D}_h$. However, it is desirable to consider
those methods where an accurate approximation of $\qn$ is obtained, since the boundary condition (\ref{new_bc}) depends on that flux. We also notice from (\ref{new_bc}) that the same idea will not work for $\widetilde{g}_N$ since a similar expression will involve derivatives of $\qn$ which are not well approximated by the numerical method.

\subsection{The HDG method} \label{sec:hdg}

The method seeks an approximation
$(\qn_h,u_h, \widehat{u}_h)$ of the exact solution $(\qn, u,
u|_{\mathcal{E}_h})$ in the space  $\Vn_h \times W_h \times M_h$
given by
\begin{subequations}
\begin{alignat}{3}
\Vn_h&=\{\vn \in [L^2(\textsf{D}_h)]^2: &&\quad\vn|_K\in  [\mathcal{P}^{k}(K)]^2&&\quad\forall K\in \textsf{D}_h\},\\
W_h&=\{w\in L^2(\textsf{D}_h): &&\quad w|_{K}\in \mathcal{P}^{k}(K)
&&\quad\forall K\in \textsf{D}_h\},\\
M_h&=\{\mu \in L^2(\mathcal{E}_h): &&\quad \mu|_e \in \mathcal{P}^{k}(e)
&&\quad \forall e
\in \mathcal{E}_h\}.
\end{alignat}
\end{subequations}
It is defined by requiring that it satisfies the equations
\begin{subequations}\label{hdg_eq}
\begin{eqnarray}
-(\Kn^{-1}\qn_h,\nabla w)_{\trian} + \langle \widehat{\qn}_h \cdot
\normal,w\rangle_{\partial\textsf{D}_h}&=&(f,w)_{\trian}\label{HDG1}\\
(\qn_h,\vn)_{\trian} -(u_h, \divv \vn)_{\trian} + \langle
\widehat{u}_h,\vn \cdot \normal \rangle_{\partial\textsf{D}_h}&=&0
,\label{HDG2}\\
 \langle \mu, \widehat{\qn}_h\cdot\boldsymbol{\nu}
\rangle_{\partial\textsf{D}_h\setminus\Gamma^h}&=&0,
\label{conserv}
\\
 \langle \mu, \widehat{u}_h \rangle_{\Gamma_D^h}&=& \langle \mu,
g_D^h \rangle_{\Gamma_D^h}, \label{DBC}\\
 \langle \mu, \widehat{\qn}_h \rangle_{\Gamma_N^h}&=& \langle \mu,
g_N^h \rangle_{\Gamma_N^h}, \label{NBC}
\end{eqnarray}
for all $(\vn,w, \mu)\in \Vn_h \times W_h\times M_h$. Here $g_D^h$ is the approximation of $\widetilde{g}_D$ proposed by \cite{CS}. More precisely, let $K\in \textsf{D}_h$. We define the operator $E^K: [\mathcal{P}^k(K)]^2 \rightarrow [\mathcal{P}^k(\mathbb{R}^2)]^2$ such that $E^K(\vn) = \vn$ for all $\vn \in [\mathcal{P}^k(K)]^2$. Then, for $\bx\in e \subset \Gamma_D^h$,
\begin{eqnarray}\label{gh_D}
\widetilde{g}_D(\bx)\approx g_D^h(\bx) : = g_D(\bar{\bx}) + \int_{\sigma(\bx)} \Kn^{-1}E^{K_e}(\qn_h)
\cdot \mn \:ds,
\end{eqnarray}
where $K_e$ is the triangle where $e$ belongs. In other words, $E^{K_e}$ is the standard extension of a polynomial to the whole $\mathbb{R}^{2}$ space.
 On the other hand, $g_N^h$ is an approximation of $\widetilde{g}_N$ which is still unknown.  In Subsection \ref{sec:gNh} we propose to replace (\ref{NBC}) by an equation involving  known quantities at the right hand side.

Finally, to complete
the definition of the HDG method we must specify the definition of
numerical trace $\widehat{\qn}_h$ on $\partial\textsf{D}_h$, which we takes of the form
\begin{eqnarray}
\widehat{\qn}_h&=& \qn_h+\tau (u_h-\widehat{u}_h)\normal,\label{uhat2}
\end{eqnarray}\label{hdg}
\end{subequations}
where $\tau: \partial \trian \rightarrow (0,\infty)$ is a stabilization parameter that guaranties solvability of (\ref{hdg_eq})  and can be set as $\tau|_K = \|\Kn\|_{L^\infty(K)}$ on each element $K$ (\cite{CGL,MIT2}).

\subsection{Definition of the family of paths}
The representation of $g_D^h$ in (\ref{gh_D}) is independent on the integration path. Let $\bx$ be a point on a boundary edge $e$. Previous work have proposed two ways to determine a point $\bar{\bx}$ in $\Gamma$ and hence construct $\sigma(\bx)$:
\begin{enumerate}
\item[{\bf (P1)}] If $\bx$ is a vertex, an algorithm developed by \cite{CS} uniquely determines $\bar{\bx}$ as the closest point to $\bx$ such that $\sigma(\bx)$ does not intersect another path before terminating at $\Gamma$ and does not intersect the interior of the domain $\Omega$. In addition, if $\bx$ is not a vertex, its corresponding path is defined as convex combination of those paths associated to the vertices of $e$.  For the Dirichlet boundary value problem, the authors in \cite{CS} numerically showed optimal rates of convergence with this choice of $\sigma(\bx)$ when $d(\Gamma,\Gamma^h)$ is of order $h$, that is, order $k+1$ for $u_h$ and $\qn_h$ and order $k+2$ for the numerical trace $\widehat{u}_h$.
\item[{\bf (P2)}]   On the other hand, \cite{CQS} proposed to determine $\bar{\bx}$ such that  $\mn$ is normal to the edge $e$. In this case these authors theoretically proved that if $d(\Gamma,\Gamma^h)$ is of order $h$, the order of convergence for $u_h$ and $\qn_h$ is indeed $k+1$, but the order for  $\widehat{u}_h$ is only $k+3/2$. However, if $d(\Gamma,\Gamma^h)$ is of order $h^{5/4}$ the numerical trace also superconverges with order $k+2$. Moreover, they also showed numerical evidence indicating that the numerical trace optimally superconverges even though $d(\Gamma,\Gamma^h)$ is of order $h$.\\
\end{enumerate}

\begin{figure*}[h!]
\begin{center}
\includegraphics[scale=0.8]{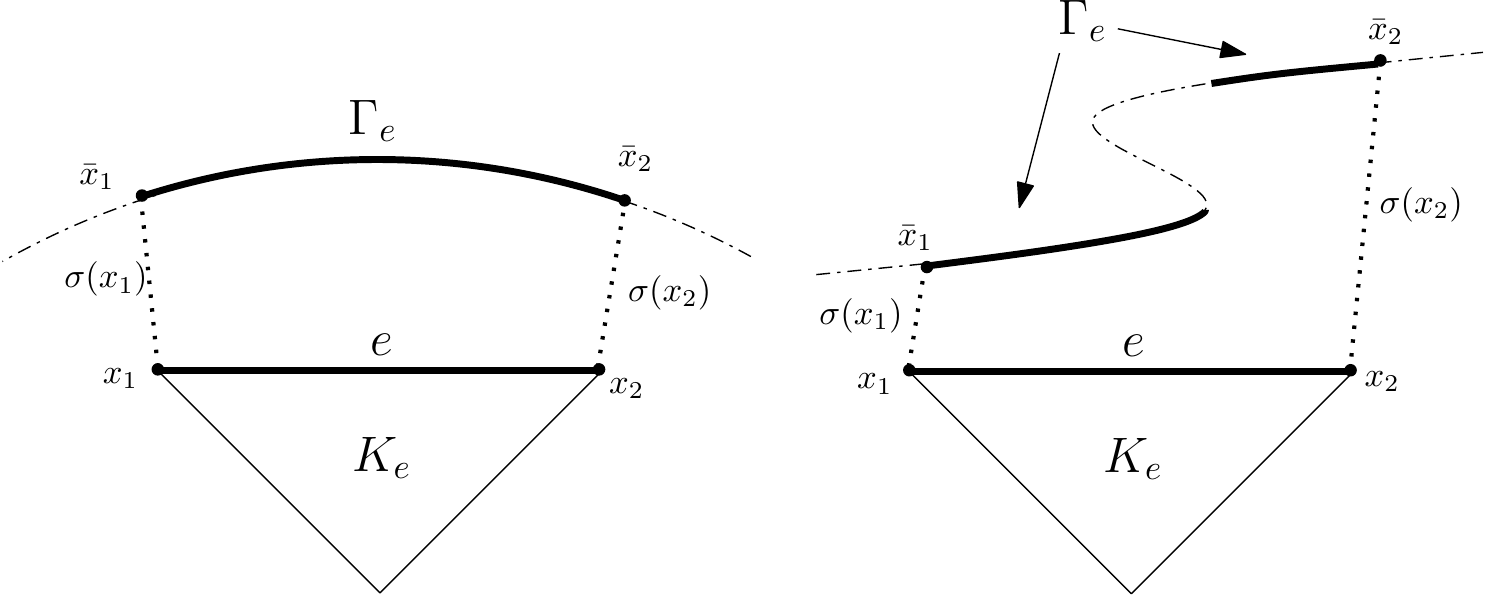}
\caption{Examples of a boundary edge $e$ with vertices $\bx_1$ and $\bx_2$. $\Gamma_e$ is the segment of $\Gamma_N$ determined by $\bar{\bx}_1$ and $\bar{\bx}_2$.}
\label{fig1}
\end{center}
\end{figure*}

Let now be $e$ a boundary edge with vertices $\bx_1$ and $\bx_2$. We denote by $\Gamma_e$ the part of $\Gamma$ determined by $\bar{\bx}_1$ and $\bar{\bx}_2$ as it is shown in Fig. \ref{fig1}. 
In this paper we assume that if $e\subset \Gamma_D^h$ (or $e\subset \Gamma_N^h)$ then $\Gamma_e\subset \Gamma_D$ $(\Gamma_e\subset \Gamma_N)$. The algorithm in {\bf (P1)} can be easily modified to satisfy this assumption. On the other hand, the paths defined in {\bf (P2)} will not always satisfy this condition.

\subsection{Approximation of the Neumann boundary condition}\label{sec:gNh}
Let $e \subset \Gamma_N^h$ a Neumann boundary face and $\Gamma_e\subset \Gamma_N$ the part of $\Gamma_N$ associated to $e$. We denote by $K_e$ the element of the triangulation where $e$ belongs. 

The main idea is to characterize $\Gamma_e$ using the parameterization {\it induced} by the family of paths. More precisely, Let $e=\{\bx: \bx(\theta) = (\bx_2-\bx_1)\theta +\bx_1, \theta \in [0,1]\}$. Then
\begin{eqnarray}
\Gamma_e = \{\bar{\bx}=\boldsymbol{\phi}(\theta): \boldsymbol{\phi}(\theta)=\bx(\theta) + |\sigma(\bx(\theta))|\mn(\theta), \theta \in [0,1]\},
\end{eqnarray}
where we recall that $|\sigma(\bx(\theta))|$ and $\mn(\theta)$ are the length and tangent vector of the segment joining $\bx(\theta)$ and $\bar{\bx}(\theta)$. 
 We define the space

\begin{eqnarray}
M_{\boldsymbol{\phi}}(\Gamma_e):=\Bigg\{\mu \in L^2(\Gamma_e): \mu = \displaystyle\frac{\tilde{\mu}\circ \boldsymbol{\phi}^{-1}}{\|\boldsymbol{\phi}'\circ \boldsymbol{\phi}^{-1}\|_2} \:\textrm{with}\:\tilde{\mu}\in \mathcal{P}_k([0,1])\Bigg\}.
\label{spaceM}
\end{eqnarray}

Equation (\ref{NBC}) is then replaced by imposing the following condition over $\qn_h$:
\begin{eqnarray}
\langle E^{K_e}(\qn_h)\cdot \bn, \mu  \rangle_{\Gamma_e} = \langle g_N , \mu \rangle_{\Gamma_e} \quad \forall \mu \in M_{\boldsymbol{\phi}}(\Gamma_e).
\label{Neumann}
\end{eqnarray}

Notice that (\ref{Neumann}) becomes
\begin{eqnarray}
\int_0^1\big(E^{K_e}(\boldsymbol{q}_h) \cdot\boldsymbol{n}) \circ \boldsymbol{\phi}\big)(\theta)\: \tilde{\mu}(\theta) d\theta= \int_0^1 (g_N \circ \boldsymbol{\phi}) \: \tilde{\mu}(\theta) d\theta 
\label{Neumann_param}
\end{eqnarray}
for all $\tilde{\mu}((\theta)) \in \mathcal{P}_k([0,1])$; hence, there is no need of computing the derivative of $\boldsymbol{\phi}$.\\

On the other hand, we observe that if $\mn$ and $\sigma$ were independent of $\theta$ (for example, if $\Gamma_e$ were polygonal and $\mn$ perpendicular to $e$), then $\|\boldsymbol{\phi}'\circ \boldsymbol{\phi}^{-1}\|_2$ would be constant and hence $M_{\boldsymbol{\phi}}(\Gamma_e)$ becomes a standard space of polynomials through pulling back polynomials from the interval $[0,1]$. As we will see in the numerical experiments provided in next section,  this technique performs optimally if  $\mn$ and $\normal$ have the same direction.

\subsection{Numerical results: boundary-value problem} \label{sec:numerical_results}

In this section we present numerical experiments showing the performance the extrapolation technique and the influence of the choice of paths.  Since the size of the computational domain changes with $h$, we measure the
errors $e_u:=u-u_h$,   $e_{\qn}:= \qn-\qn_h$ and $e_{\widehat{u}} :=u-\widehat{u}_h$ by using the following norms:
\begin{eqnarray*}
\|e_u\|_{\ltwoint}: &=& \frac{\|e_u\|_{L^2(\textsf{D}_h)}}{{|\textsf{D}_h|}^{1/2}}, \:\|e_{\qn}\|_{\ltwoint}: = \frac{\|e_{\qn}\|_{[L^2(\textsf{D}_h)]^2}}{|\textsf{D}_h|^{1/2}},\\
\|e_{\hat{u}}\|_{\mathcal{E}_h}:&=&\left(\frac{ \sum_{K\in
\textsf{D}_h} h_K \|\textsf{P}_{\partial} u - \hat{u}_h
\|^2_{L^2(\partial K)}}{ \sum_{K\in\textsf{D}_h} h_k |\partial
K|}\right)^{1/2}.
\end{eqnarray*}
Here $\textsf{P}_{\partial} $ is the $L^2-$projection over
$\mathcal{P}^k(e)$ with $e\subset \partial K$.

In addition we compute an element-by-element
postprocessing, denoted by $u^*_h$, of the approximate solution $u_h$, which
provides a better approximation for the scalar variable when $k\ge 1$ (\cite{CGS,CockburnGuzmanWang09}). Given an element $K$ we construct 
$u^*_h=\bar{u}_h + \tilde{u}_h$
as
the only function in $\mathcal{P}^{k+1}(K)$ such that
\begin{equation*}
\bar{u}_h = \left\{%
\begin{array}{cc}
  \frac{1}{3} \sum_{e\in \partial K} {\widehat{u}_h|_e} & \textrm{if} \: k=0 ,\\
  &\\
  \frac{1}{|K|} \int_K u_h dx &  \textrm{if} \: k>0 ,\\
\end{array}%
\right.
\end{equation*}
and $\tilde{u}_h$ is the polynomial in $\mathcal{P}^{k+1}_0(K)$ (set
of functions in $\mathcal{P}^{k+1}(K)$ with mean zero) satisfying
$$
(\nabla \tilde{u}_h, \nabla w )_K = -(\qn_h,\nabla w)_K  \quad \forall w \in
\mathcal{P}^{k+1}(K).
$$
In the purely diffusive case, this new approximation of $u$ has been proven to converge with order $k+2$ for $k\geq1$ when the
domain is polygonal (\cite{CGS,CockburnGuzmanWang09}), and also when it has curved Dirichlet boundary (\cite{CQS,CS}).\\

We set $\Kn = \mathbf{I}$  in all the experiments of this section. 
In Subsection~\ref{subsec:h} we show that deteriorate convergence can happen if $d(\Gamma,\Gamma^h)=O(h)$. However, we will see
in Subsection~\ref{subsec:h2} that optimal convergence is obtained when $d(\Gamma,\Gamma^h)=O(h^{2})$. 

\subsubsection{Computational domain at a distance $d(\Gamma,\Gamma^h)=O(h)$}\label{subsec:h}
In the following examples the computational domain is constructed in such a way that the distance $d(\Gamma,\Gamma^h)$ is of order $h$. Moreover,  $f$, $g_D$ and $g_N$ are chosen in order that $u(x,y)=\sin(x)\sin(y)$ is solution the exact of (\ref{model_problem}).\\


\begin{example}\label{ex1}
{\rm Our first example consist of approximating a squared domain $\Omega = (0,1)$ by a squared subdomain satisfying 
$d(\Gamma,\Gamma^h)=O(h)$ as Fig. \ref{fig_ex1} shows. Let $\Gamma_N = \{x: x = 0\}$, 
$\Gamma_D=\partial \Omega \setminus \Gamma_N$ and  the family of paths is computed according to {\bf (P2)}.

\begin{figure}[ht!]
\begin{center}
\includegraphics[scale=0.25]{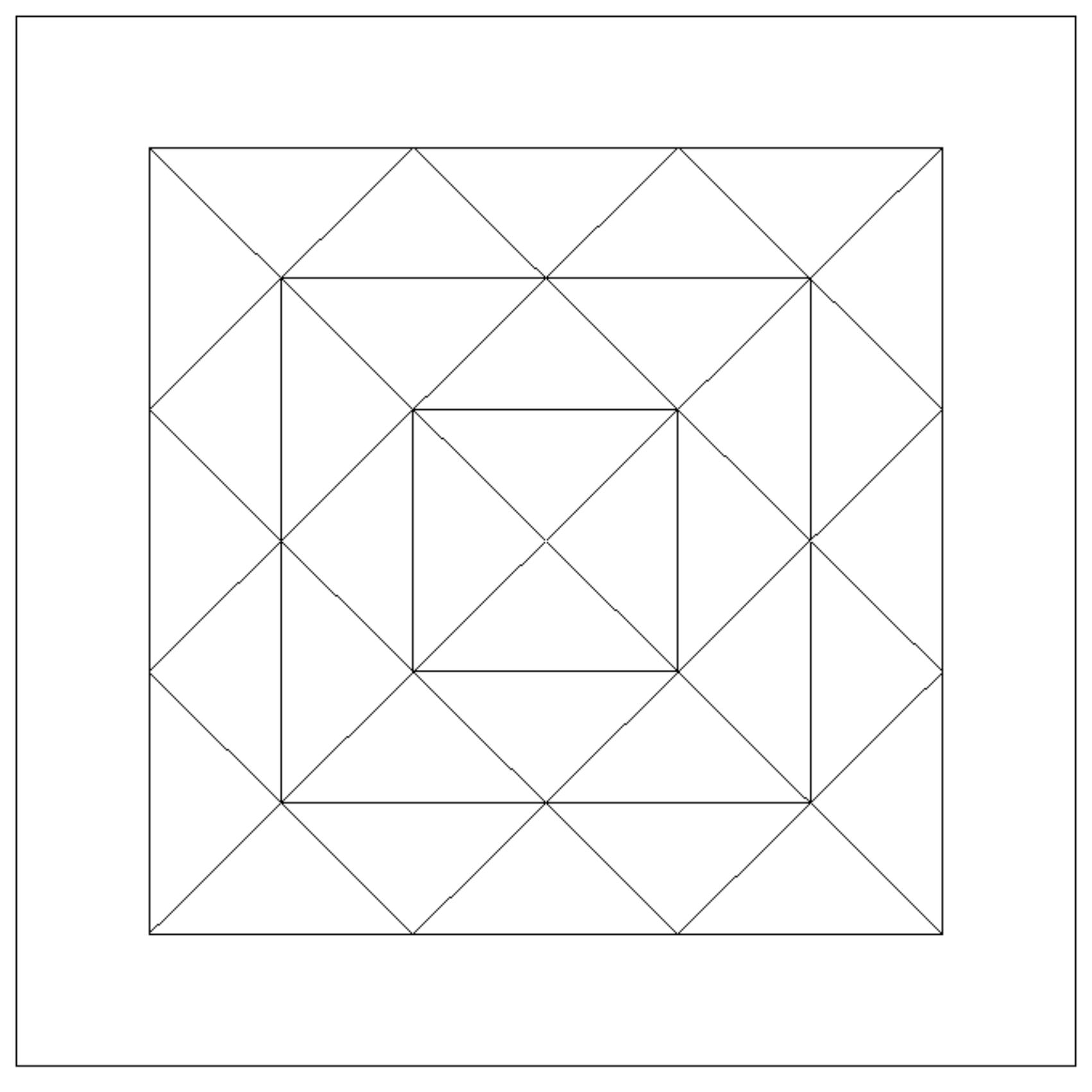}
\includegraphics[scale=0.25]{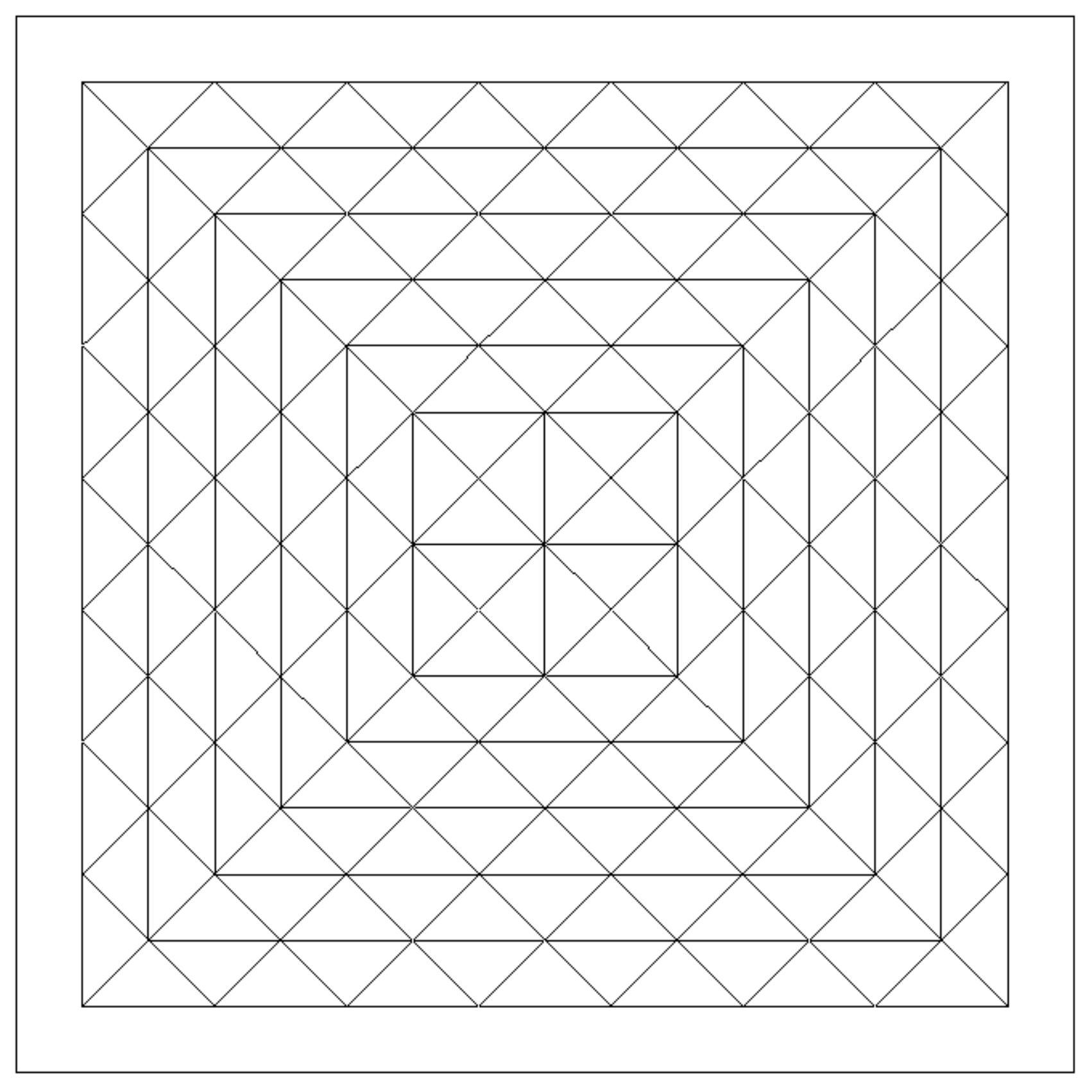}
\caption{Two consecutive meshes ($h=1/4$ and $h=1/8$) approximating the domain of Example \ref{ex1}. (Figure obtained from \cite{CQS})}
\label{fig_ex1}
\end{center}
\end{figure}

In Table~\ref{table:R1_square} we display the history of convergence for different polynomial degree ($k=0,1,2$ and $3$) and meshsizes ($h = 1/2, 1/4, 1/8, 1/16$ and $1/32$). We observe that the error of $u$ and $\qn$ behaves optimally with convergence rate of order $k+1$. Moreover the error of numerical trace and postprocessed solution also converge with order $k+1$, which is not optimal for the standard HDG method on polygonal domains. 
Even though, the errors $e_{u^*}$ are always small than $e_u$. We attribute this lack of superconvergence to the fact that the Neumann condition (\ref{Neumann}) is being imposed on $\qn_h$ and not on $\widehat{\qn}_h$ as in the standard HDG method.  \\

\begin{table}[htbpa!]\renewcommand{\arraystretch}{1.3}\addtolength{\tabcolsep}{-5pt}
\begin{center}
{\scriptsize\begin{tabular}{c|c|cc|cc|cc|cc}
  \hline \hline
  \multicolumn{2}{c}{}     &
  \multicolumn{2}{c}{$\|e_u\|_{\ltwoint}$}     &
  \multicolumn{2}{c}{$\|e_{\qn}\|_{\ltwoint}$} &
  \multicolumn{2}{c}{$\|e_{\widehat{u}}\|_{\mathcal{E}_h}$}&
  \multicolumn{2}{c}{$ \|e_{u^*}\|_{\ltwoint}$}\\
  $k$& $h$& error & order &error&order &error &order &error &order \\
  \hline \hline
                    		  & 1/2  &4.58E-03  & -         & 6.59E-02&  -    & 2.13E-02& -    & 7.50E-03&  -    \\
							  & 1/4  &6.09E-03  & -0.41  & 4.77E-02&0.46& 5.75E-03&1.89 &6.60E-03&0.18 \\
	0						  & 1/8  &4.62E-03  & 0.40    & 2.74E-02&0.80 & 1.75E-03&1.71 &4.71E-03&0.49  \\
                     		  & 1/16&2.78E-03  & 0.73    & 1.46E-02&0.91 	& 6.18E-04&1.51 &2.80E-03&0.75  \\
                     		  & 1/32&1.52E-03  & 0.87    & 7.52E-03&0.96 	& 2.51E-04&1.30 &1.53E-03&0.88  \\
  \hline
                    		  & 1/2  &1.54E-03  & -       &9.89E-03 &  -  		  &3.70E-03 & -    &1.67E-03&  -    \\
							  & 1/4  &5.67E-04   &1.44  &2.55E-03 &1.96	  &6.31E-04 &2.55&4.68E-04&1.84\\
1							  & 1/8  &1.69E-04  &1.75   &7.09E-04 &1.85 	  &1.50E-04 &2.07&1.31E-04&1.83  \\
                     		  & 1/16&4.62E-05  &1.86   &1.94E-04 &1.87 	  &3.84E-05 &1.97&3.60E-05&1.87\\
                     		  & 1/32&1.21E-05  &1.93   &5.13E-05 &1.92 	  &9.83E-06 &1.97&9.52E-06&1.92\\
   \hline
                    		  & 1/2  &2.29E-04  & -       &1.20E-03 &  -  		  &5.23E-04 & -        &2.17E-04&   -    \\
							  & 1/4  &2.82E-05  & 3.02  &1.24E-04 & 3.28	  &3.36E-05 &3.96 	  &2.44E-05 &3.16\\
2							  & 1/8  &3.43E-06  & 3.03  &1.36E-05 & 3.19	  &3.22E-06 &3.38 	  &2.81E-06 &3.12  \\
                     		  & 1/16&4.25E-07  & 3.01  &1.63E-06 & 3.06	  &3.61E-07 &3.16	  &3.38E-07 &3.05\\   
                       		  & 1/32&5.28E-08  & 3.01  &2.02E-07 & 3.01	  &4.26E-08 &3.08	  &4.13E-08 &3.03\\                     		                   		  
  \hline
                    		  & 1/2 &  3.37E-05& -		 &1.51E-04 &  -        &7.55E-05 & -         &3.39E-05&  -    \\
							  & 1/4 & 2.30E-06&3.87   &9.32E-06 &4.02	      &3.12E-06 & 4.59	  &2.30E-06&3.88\\
3							  & 1/8  & 1.55E-07&3.89   &6.74E-07 &3.79 	  &1.78E-07 & 	4.14   &1.55E-07&3.89  \\
                     		  & 1/16& 1.05E-08 &3.89   &4.76E-08 &3.82 	  &1.12E-08 &	3.99	  & 1.05E-08&3.89  \\
							  & 1/32& 6.90E-10 &3.92   &3.22E-09 &3.89 	  &7.13E-10 &	3.97	  & 6.90E-010&3.92  \\                     		  
 \hline \hline
 \end{tabular}  }
\end{center}
\caption{History of convergence of the approximation in Example \ref{ex1}.}\label{table:R1_square}
\end{table}
}
\end{example}


\begin{example}\label{ex2}
{\rm
We now consider an annular domain $\Omega =\{(x,y) \in \mathbb{R}^2: 14^2<x^2+y^2<20^2\}$ that is being approximated by a polygonal subdomain satisfying $d(\Gamma,\Gamma^h) = O(h)$ as shown in Fig. \ref{fig_ex2}. We consider Neumman data in the outer boundary  $\Gamma_N = \{(x,y) \in \mathbb{R}^2 : x^2+y^2 = 20^2\}$ and Dirichlet data in the inner boundary $\Gamma_D = \{(x,y): x^2+y^2 = 14^2\}$. Here the paths are computed according to {\bf (P2)}. 

\begin{figure}[ht!]
\begin{center}
\includegraphics[scale=0.25]{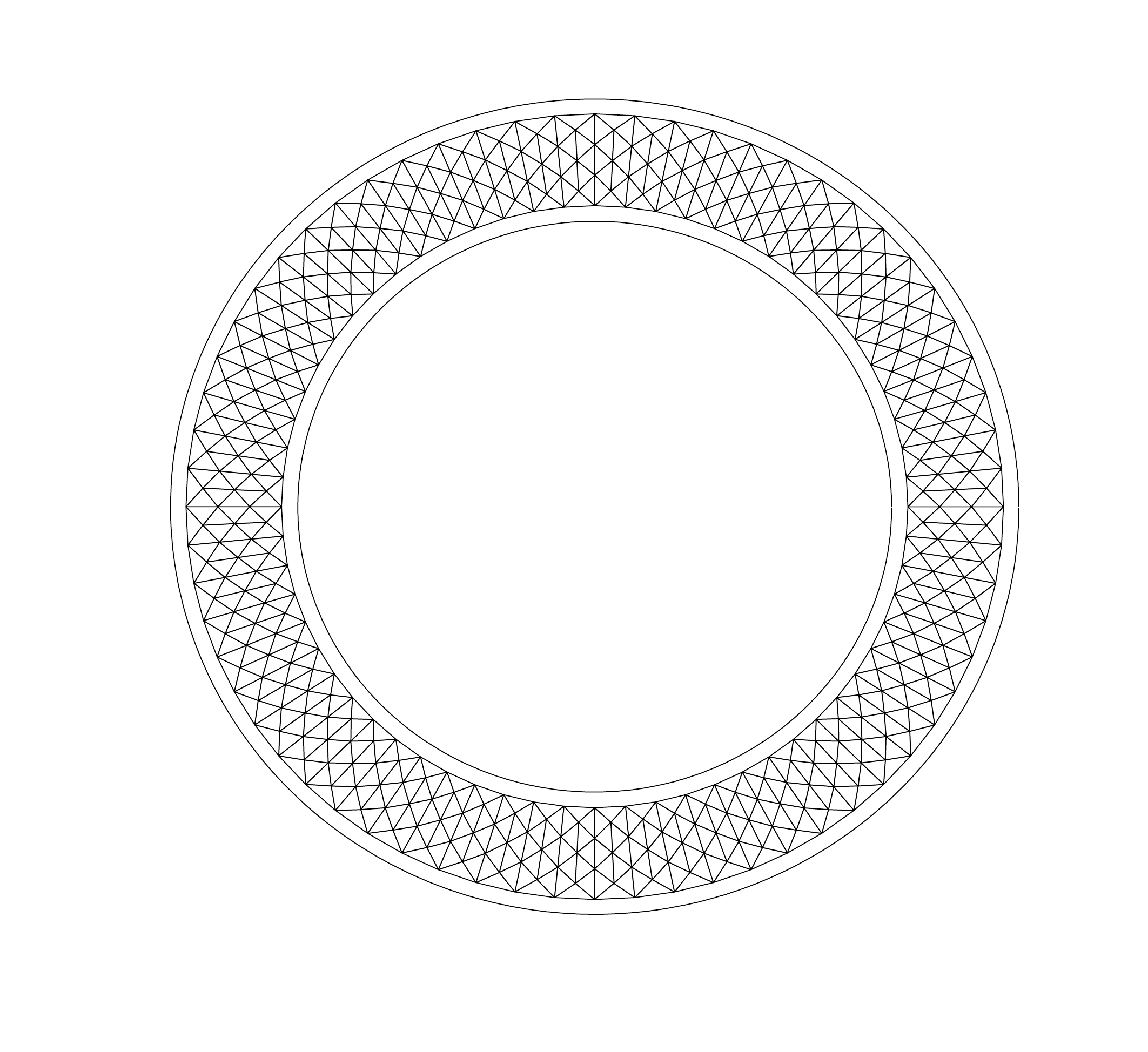}
\caption{Annular domain and mesh in Example \ref{ex2}.}
\label{fig_ex2}
\end{center}
\end{figure}

The behavior of the $L^2$-norm of the error displayed in Table~\ref{table:R1_ring} is similar to the one obtained in the previous example, i.e., the rate of convergence of the error in all the variables is of order $k+1$. Thus, this example suggests that our technique performs properly when the boundary is actually non-polygonal.\\

\begin{table}[htbpa!]\renewcommand{\arraystretch}{1.3}\addtolength{\tabcolsep}{-5pt}
\begin{center}
{\scriptsize\begin{tabular}{c|c|cc|cc|cc|cc}
  \hline \hline
  \multicolumn{2}{c}{}     &
  \multicolumn{2}{c}{$\|e_u\|_{\ltwoint}$}     &
  \multicolumn{2}{c}{$\|e_{\qn}\|_{\ltwoint}$} &
  \multicolumn{2}{c}{$\|e_{\widehat{u}}\|_{\mathcal{E}_h}$}&
  \multicolumn{2}{c}{$ \|e_{u^*}\|_{\ltwoint}$}\\
  $k$& $h$& error & order &error&order &error &order &error &order \\
  \hline \hline
                		    &1.89&9.56E+00& -    &8.79E+00&  -   &4.66E-01& -    &9.80E+00&  -    \\
							&0.96&8.47E+00&0.18&5.82E+00&0.61&3.72E-01&0.33&8.50E+00&0.21\\
0							&0.49&5.72E+00&0.57&3.38E+00&0.79&2.42E-01&0.63&5.72E+00&0.56   \\
                  			&0.24&3.29E+00&0.81&1.82E+00&0.90&1.37E-01&0.83&3.29E+00&0.81  \\
                 			&0.12&1.76E+00&0.91&9.42E-01&0.91 &7.26E-02&0.92&1.76E+00&0.91  \\
  \hline
                  			 & 1.89&2.03E+01& -    &7.85E+00&  -   &9.56E-01& -    &2.04E+01&  -    \\
				 			 & 0.96&5.94E+00&1.82&2.12E+00&1.94&2.58E-01&1.94&5.96E+00&1.82\\
1							 & 0.49&1.43E+00&2.08&5.03E-01&2.10&6.00E-02&2.13&1.43E+00&2.08   \\
                  			 & 0.24&3.40E-01&2.09 &1.20E-01&2.08&1.40E-02&2.11&3.40E-01&2.09  \\
                  			 & 0.12&8.19E-02&2.06 &2.92E-02&2.06&3.35E-03&2.11&8.20E-02&2.06  \\
  \hline
                  			 & 1.89&4.04E+00& -  &1.82E+00&  - &1.90E-01& -  &4.04E+00&  -    \\
		          			 & 0.96&6.80E-01&2.64&3.42E-01&2.46&2.95E-02&2.76&6.81E-01&2.64\\
2							 & 0.49&1.41E-01&2.30&5.86E-02&2.58&5.89E-03&2.36&1.41E-01&2.30   \\
                 			 & 0.24&2.12E-02&2.75&8.33E-03&2.83&8.75E-04&2.77&2.12E-02&2.75  \\
							 & 0.12&2.88E-03&2.89&1.10E-03&2.93&1.16E-04&2.90&2.88E-03&2.93  \\
  \hline
                 			& 1.89&4.12E+00& -    &1.52E+00&  -  &1.93E-01& -  &4.12E+00&  -    \\
							& 0.96&3.17E-01&3.80&1.07E-01&3.93&1.37E-03&3.92&3.17E-01&3.80\\
3							& 0.49&1.89E-02&4.13&6.29E-03&4.15&7.89E-04&4.18&1.89E-02&4.13   \\
                 	        & 0.24&1.10E-03&4.13&3.70E-04&4.12&4.53E-05&4.15&1.10E-03&4.13  \\   
                 	        & 0.12&6.56E-05&4.08&2.23E-05&4.07&2.68E-06&4.09&6.56E-05&4.08  \\  
 \hline \hline
 \end{tabular}  }
\end{center}
\caption{History of convergence of the approximation in Example \ref{ex2}.}\label{table:R1_ring}
\end{table}
}
\end{example}

\begin{figure*}[h!]
\begin{center}
\includegraphics[scale=0.9]{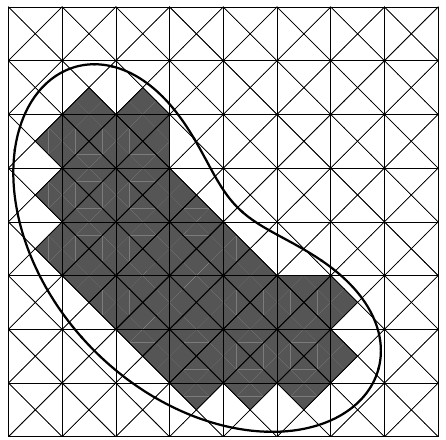}
\caption{Left: Domain $\Omega$, its boundary $\Gamma$ (solid line), a
background mesh $\mathcal{B}_h$ and the polygonal subdomain $\textsf{D}_h$ (gray). Right: Dirichlet data $g$ on $\Gamma$ transferred to
$\varphi$ on $\Gamma_h$. (Figure taken from \cite{CS})}
\label{fig:mesh_Dir}
\end{center}
\end{figure*}

\begin{remark}
The construction of the family of paths according to {\bf (P1)} in Examples \ref{ex1} and \ref{ex2} deliver similar results since the difference between  {\bf (P1)} and  {\bf (P2)} is not significant for these domains. That is why we do not display the convergence tables for this case.  This numerical evidence indicates that the technique proposed provides optimal 
rate of convergence when $d(\Gamma,\Gamma^h)=O(h)$ and the family of paths is constructed according to  {\bf(P1)} or {\bf(P2)}. 
However, in practice, this condition over the distance can not be satisfied in general, unless the mesh is constructed properly 
to do so.

A practical construction of the computational domain $\trian$ was described in \cite{CS}. It consists of ``immersing'' the domain 
in a Cartesian background mesh and set $\trian$ as the union of all the elements that are completely inside of $\Omega$ as it is 
shown in Fig. \ref{fig:mesh_Dir}. Here $d(\Gamma,\Gamma^h)=O(h)$. In this case it is not convenient to construct the paths according 
to {\bf(P2)}. In fact, given a point $\bx \in \mathcal{E}^{\partial}_h$ it might happen that $\bar{\bx}$ is extremely far from 
$\bx$, specially in parts of $\Gamma$ where the domain is non-convex. Since both procedures deliver similar results in previous examples, we will consider from now on {\bf(P1)}.
 \end{remark}\\
 
\begin{figure}[ht!]
\begin{center}
\includegraphics[scale=0.26]{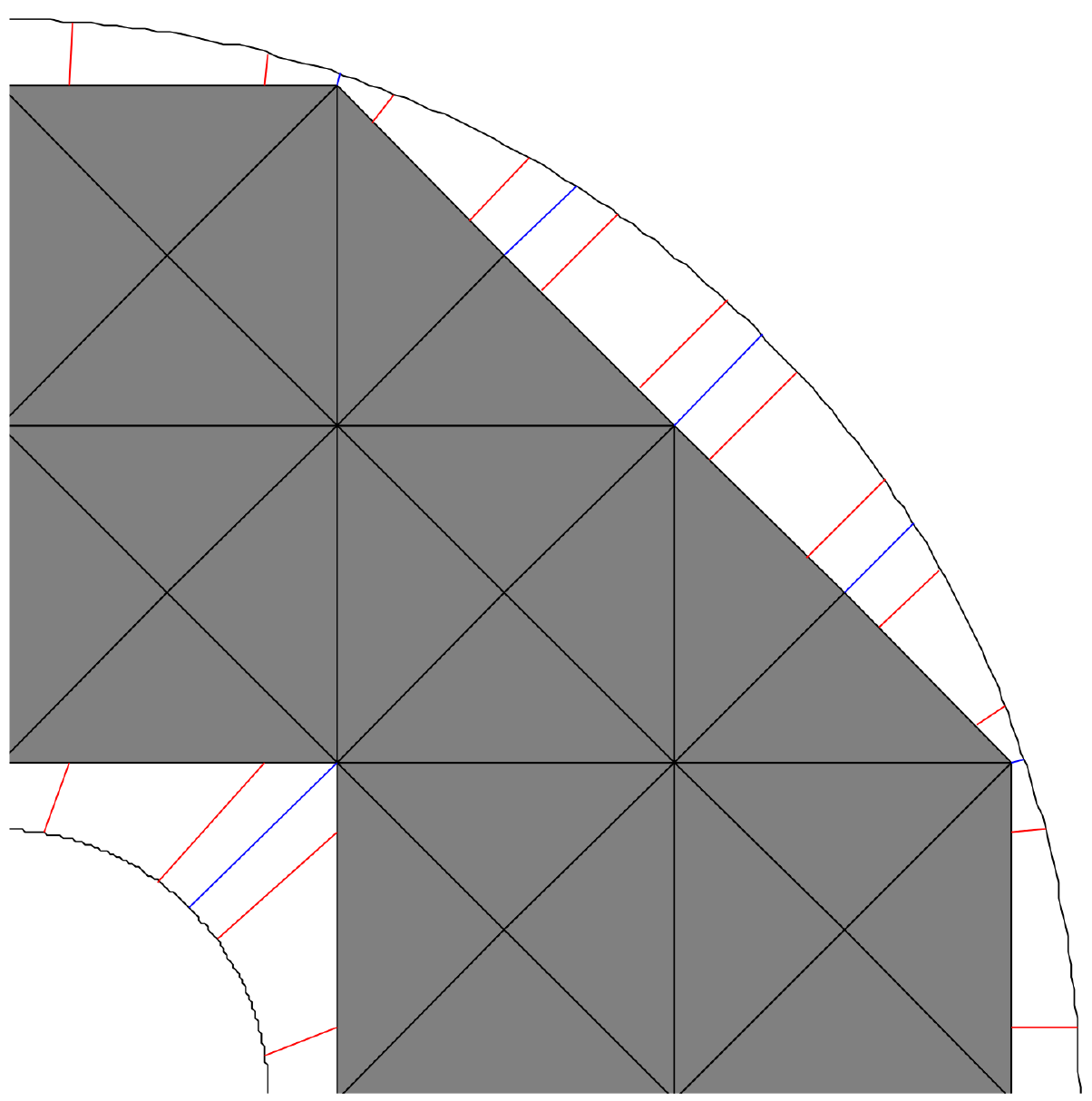}
\includegraphics[scale=0.26]{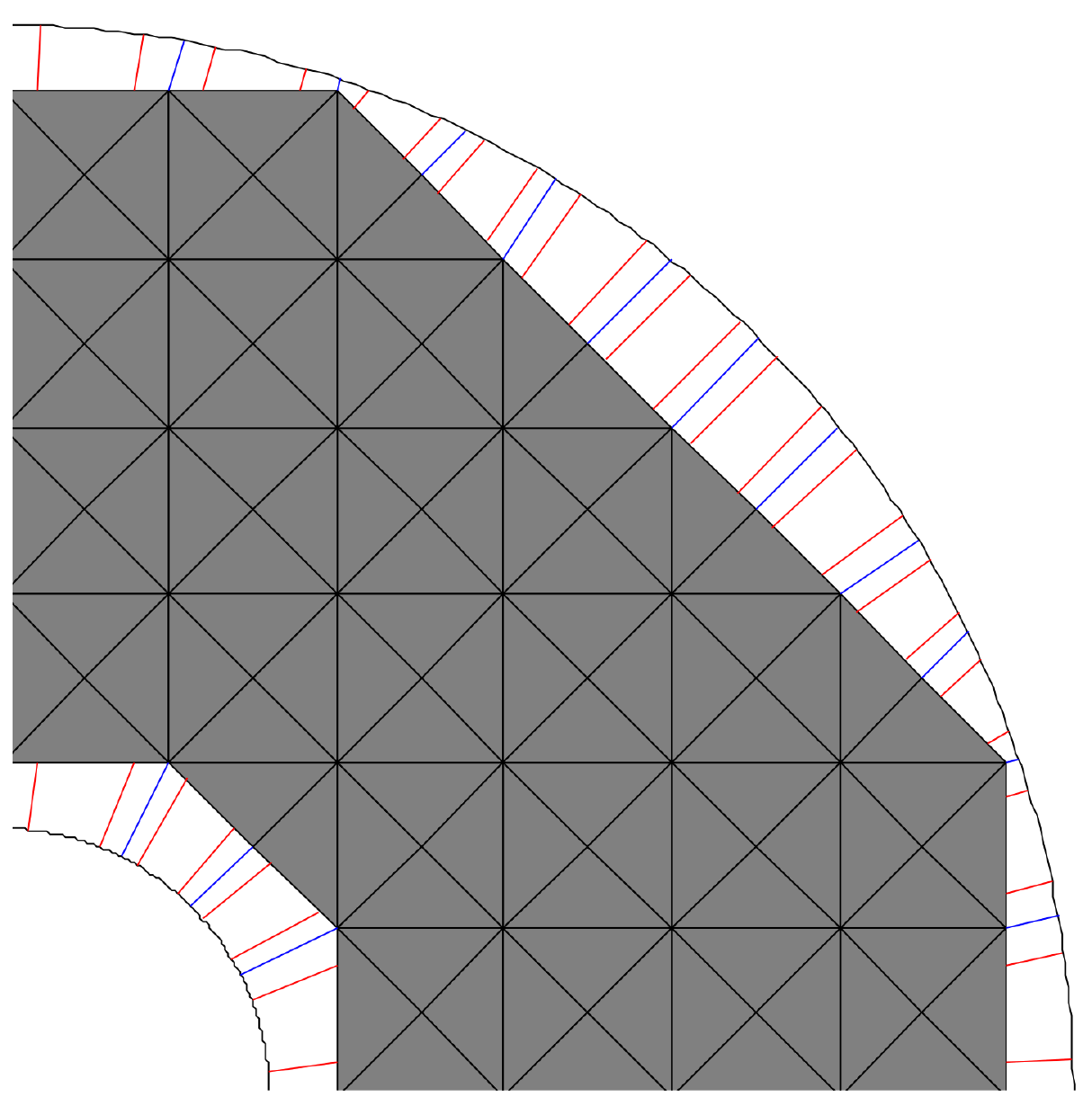}\\
\includegraphics[scale=0.26]{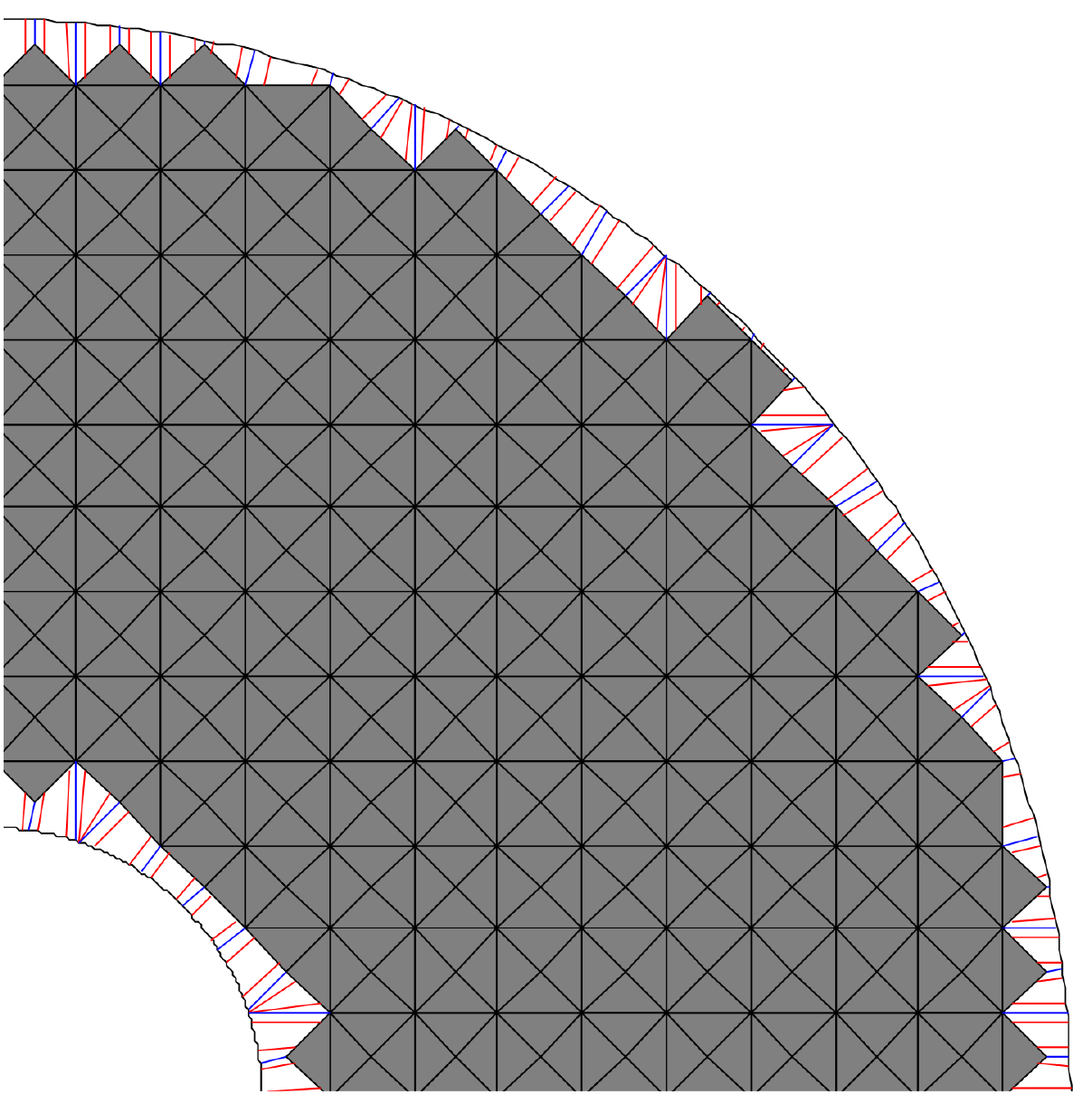}
\caption{Zoom at the upper-right corner of three consecutive meshes of Example \ref{ex3}. Mesh (grey region) constructed considering the procedure in \cite{CS} and family of paths determined according to {\bf (P1)}. Blue lines: paths from the vertices. Red lines: paths from quadrature points of the boundary edges ($k=1$). }
\label{fig_ex3}
\end{center}
\end{figure}

\begin{example}\label{ex3} 
{\rm
In order to observe the performance of the method where the mesh satisfies  $d(\Gamma,\Gamma^h)=O(h)$ and  the paths are given by {\bf(P1)}, we consider the ring $\Omega =\{(x,y) \in \mathbb{R}^2: 0.25^2<(x-0.5)^2+(y-0.5)^2<1\}$  with $\Gamma_N = \{(x,y) \in \mathbb{R}^2 : x^2+y^2 = 1\}$ and $\Gamma_D = \{(x,y): x^2+y^2 = 0.25^2\}$. In Fig. \ref{fig_ex3} we show a zoom at the upper-right corner of three consecutive meshes. We also plot the family paths from vertices and quadrature points on the boundary edges.  In Table~\ref{table:ex_3} we display the history of convergence. Even though the method is still convergent for $k=0$, $1$ and $2$, the rates deteriorate. Moreover, there is no convergence when $k=3$ . For the Dirichlet boundary value problem this non-optimal behavior does not occur as \cite{CS} showed. This example suggests that in a practical situation (meshes satisfying $d(\Gamma,\Gamma^h)=O(h)$ and paths constructed using  {\bf (P1)}, the method does not perform properly. So, it seems that for Neumann boundary data, the family of paths needs to be build according to {\bf(P2)}. Even though we have no theoretical support that explains this behavior, we believe it might be related to the oscillatory nature of high degree polynomials. In fact, for the Dirichlet problem, \cite{CQS} showed error estimates where some of the constants depend on the polynomial degree. In addition, \cite{CS2} numerically studied the robustness of this method applied to a convection-diffusion problem with Dirichlet boundary data. The concluded that, even though   $d(\Gamma,\Gamma^h)=O(h)$, $\Gamma$ and $\Gamma^h$ must be ``close enough'  when $k\geq 1$. 

One way of always being able to construct the paths  using  {\bf(P2)}  is to interpolate the boundary by a piecewise linear function. In this case $d(\Gamma,\Gamma^h)=O(h^2)$. }\\
\end{example}

\begin{remark}
In Example \ref{ex3} it is not possible to construct the family of path by {\bf(P1)}. In fact, a path perpendicular to an inner boundary edge might not intersect the inner ring . Moreover, a path perpendicular  to an outer boundary edge might intersect the outer boundary extremely ``far'' as would happen in the third mesh of Fig. \ref{fig_ex3}.
\end{remark}

\begin{table}[htbpa!]\renewcommand{\arraystretch}{1.3}\addtolength{\tabcolsep}{-5pt}
\begin{center}
{\scriptsize\begin{tabular}{c|c|cc|cc|cc|cc}
  \hline \hline
  \multicolumn{2}{c}{}     &
  \multicolumn{2}{c}{$\|e_u\|_{\ltwoint}$}     &
  \multicolumn{2}{c}{$\|e_{\qn}\|_{\ltwoint}$} &
  \multicolumn{2}{c}{$\|e_{\widehat{u}}\|_{\mathcal{E}_h}$}&
  \multicolumn{2}{c}{$ \|e_{u^*}\|_{\ltwoint}$}\\
  $k$& $h$& error & order &error&order &error &order &error &order \\
  \hline \hline
                		    &0.312&4.12E-02& -    &1.83E-01 &  -     &4.40E-02& -    &  4.15E-02 &- \\
							&0.156&3.70E-02&0.16&1.27E-01 &0.53  &3.26E-02&0.43&  3.69E-02& 0.17   \\
0							&0.078&1.69E-02&1.13&1.37E-01 &-0.11&1.50E-02&1.12&  1.69E-02& 1.13   \\
                  			&0.039&9.11E-03&0.89&7.00E-02 &0.96  &7.61E-03&0.97&  9.11E-03& 0.89   \\
                 			&0.019&8.50E-03&0.10&4.92E-02 &0.51  &5.66E-03&0.43&  8.50E-03& 0.10   \\
  \hline
                		    &0.312&6.13E-03&-       &1.82E-02& -     &3.75E-03&  -     &   5.71E-03&-   \\
							&0.156&3.44E-03&0.84  &1.06E-02&0.77 &2.18E-03&0.78  &  3.37E-03& 0.76  \\
1							&0.078&3.86E-03&-0.17&9.41E-03&0.18  &2.36E-03&-0.11& 3.86E-03&-0.20     \\
                  			&0.039&1.16E-03&1.74  &2.68E-03&1.81  &6.88E-04&1.78 &  1.16E-03& 1.73   \\
                 			&0.019&5.17E-04&1.16  &1.16E-03&1.20  &3.04E-04&1.18 &  5.16E-04& 1.16   \\
  \hline
                		    &0.312&4.68E-04& -    &1.25E-03&  -   &3.03E-04& -       & 4.60E-04&-\\
							&0.156&2.25E-04&1.06&5.89E-04&1.08&1.45E-04&   1.06& 2.24E-04& 1.04  \\
2							&0.078&1.21E-04&0.89&3.24E-04&0.86&7.39E-05&   0.97& 1.21E-04& 0.89 \\
                  			&0.039&1.31E-05&3.20&3.60E-05&3.17&7.79E-06&  3.25 & 1.31E-05& 3.21  \\
                 			&0.019&2.63E-06&2.32&7.03E-06&2.35&1.54E-06&  2.33 & 2.63E-06& 2.32   \\
  \hline
                		    &0.312& 3.02E-05&-       &8.78E-05&  -    &1.98E-05& -      & 3.00E-05&- \\
							&0.156& 1.11E-05&1.44  &3.45E-05&1.35 &7.19E-06&1.45 & 1.10E-05& 1.44   \\
3							&0.078& 1.65E-06&2.75  &5.37E-06&2.67 &1.01E-06&2.83  & 1.65E-06& 2.75   \\
                  			&0.039& 6.69E-06&-       &1.53E-05&-       &3.98E-06&-       &6.70E-06&-   \\
                 			&0.019& 8.03E-03&-       &2.26E-02&-	     &4.73E-03&-        & 8.04E-03&-\\
 \hline \hline
 \end{tabular}  }
\end{center}
\caption{History of convergence of the approximation in Example \ref{ex3}.}\label{table:ex_3}
\end{table}


\subsection{Computational domain at a distance $d(\Gamma,\Gamma^h)=O(h^2)$}\label{subsec:h2}

Another practical construction of $\trian$ is defining first $\Gamma^h$ by interpolating $\Gamma$ using piecewise linear segments. Then, $\trian$ is the domain enclosed by $\Gamma_h$ as Fig. \ref{fig_ex4} shows. In this case $d(\Gamma,\Gamma_h) = O(h^2)$ and the family of paths can be easily defined according to {\bf (P2)}.\\


\begin{example} \label{ex4}
{\rm We consider the domain $\Omega =\{(x,y) \in \mathbb{R}^2: 1<(x-0.5)^2+(y-0.5)^2<4\}$  with $\Gamma_N = \{(x,y) \in \mathbb{R}^2 : x^2+y^2 = 1\}$ and $\Gamma_D = \{(x,y): x^2+y^2 = 4\}$.  In Table~\ref{table:ex_4} we observe again that the order of convergence in all the variables in $k+1$. We point out that part of the computational domain is outside of $\Omega$ as it can be observed in the inner circle in Fig.  \ref{fig_ex4}. This was never the case in the examples provided by \cite{CS} and \cite{CQS}. Thus,  these results indicates that their technique also works when $\Omega^c\cap \trian \neq\emptyset$. In Fig. \ref{fig2_ex4} we show the approximated solution $p_h$ considering  $h=1.10$ (left)  and $0.55$ (right) and using polynomials of degree $k=0,1$ and $2$. We clearly see an improvement either when the mesh is refined or the polynomial degree increases.

\begin{figure}[ht!]
\begin{center}
\includegraphics[scale=0.25]{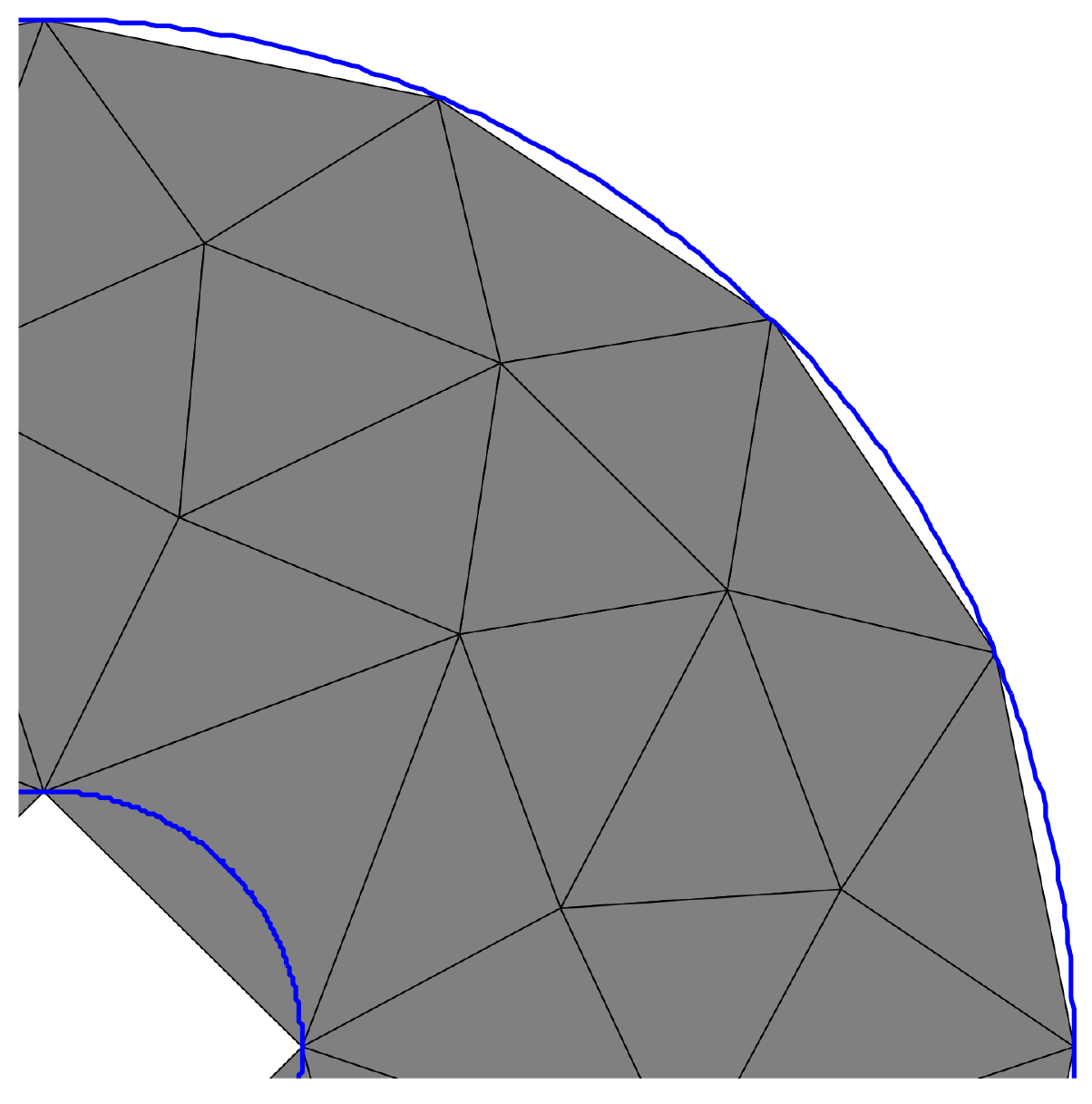}
\caption{Zoom at the upper-right corner of Example \ref{ex4}. Blue line: boundary $\Gamma$. Grey region: mesh.}
\label{fig_ex4}
\end{center}
\end{figure}

\begin{table}[htbpa!]\renewcommand{\arraystretch}{1.3}\addtolength{\tabcolsep}{-5pt}
\begin{center}
{\scriptsize\begin{tabular}{c|c|cc|cc|cc|cc}
  \hline \hline
  \multicolumn{2}{c}{}     &
  \multicolumn{2}{c}{$\|e_u\|_{\ltwoint}$}     &
  \multicolumn{2}{c}{$\|e_{\qn}\|_{\ltwoint}$} &
  \multicolumn{2}{c}{$\|e_{\widehat{u}}\|_{\mathcal{E}_h}$}&
  \multicolumn{2}{c}{$ \|e_{u^*}\|_{\ltwoint}$}\\
  $k$& $h$& error & order &error&order &error &order &error &order \\
  \hline \hline
                		    &1.72 &5.31E-01 & -    & 2.14E+00&  -     &2.22E-01& -     &  6.63E-01 &- \\
							&1.10 &2.87E-01& 1.37& 1.19E+00& 1.3   &1.14E-01& 1.48 &   3.00E-01& 1.77 \\ 
0						    &0.55 &1.45E-01& 0.99& 6.13E-01& 0.95  &5.76E-02& 1.00&  1.46E-01& 1.04\\ 
                  			&0.29 &8.10E-02& 0.89& 3.31E-01& 0.95  &3.10E-02& 0.95 &   8.05E-02& 0.91 \\                    			 
                  			&0.15 &4.36E-02& 0.98& 1.69E-01& 1.07  &1.60E-02& 1.05&  4.34E-02& 0.98 \\  
                  			&0.08 &2.24E-02& 0.99& 8.48E-02& 1.02  &8.12E-03& 1.01 &   2.23E-02& 1.00 \\  
  \hline
                		    &1.72 &2.59E-01 & -    & 9.51E-03&  -     &9.51E-03& -      & 1.22E-01 & -     \\
							&1.10 &7.11E-02& 2.89& 1.61E-03& 3.97 &1.61E-03& 3.97&1.80E-02& 4.27 \\ 
1						    &0.55 &1.77E-02& 2.01& 2.50E-04& 2.68 &2.50E-04& 2.68&2.54E-03& 2.82 \\ 
                  			&0.29 &4.45E-03& 2.12& 5.92E-05& 2.22 &5.92E-05& 2.22&4.23E-04& 2.76\\                    			 
                  			&0.15 &1.08E-03& 2.26& 1.43E-05& 2.25 &1.43E-05& 2.25&9.03E-05& 2.45 \\  
                  			&0.08 &2.66E-04& 2.08& 4.24E-06& 1.81 &4.24E-06& 1.81&2.69E-05& 1.80  \\  
  \hline
                		    &1.72 &4.59E-02 & -    & 6.22E-02&  -     &1.43E-03& -     &   1.04E-02&- \\
							&1.10 &6.55E-03& 4.35& 9.09E-03& 4.29 &1.95E-04& 4.44&   1.35E-03& 4.56 \\ 
2						    &0.55 &8.37E-04& 2.97& 1.26E-03& 2.85 &1.10E-05& 4.15&   8.25E-05& 4.03 \\ 
                  			&0.29 &1.12E-04& 3.09 &1.71E-04& 3.07 &2.14E-06& 2.52&   1.44E-05& 2.67 \\                    			 
                  			&0.15 &1.42E-05& 3.29 &2.11E-05& 3.32 &2.01E-07& 3.75&   1.34E-06& 3.77 \\  
                  			&0.08 &1.77E-06& 3.10 &2.63E-06& 3.10 &3.37E-08& 2.66&   2.22E-07 &2.68 \\  
  \hline
                		    &1.72 &5.61E-03& -     & 8.48E-03&   -    &.57E-04& -     &1.28E-03&   - \\
							&1.10 &4.47E-04& 5.65& 6.59E-04& 5.71&6.52E-06& 7.11&4.82E-05& 7.32 \\ 
3						    &0.55 &3.31E-05& 3.75& 4.77E-05& 3.78&1.77E-07& 5.20&1.42E-06& 5.08 \\ 
                  			&0.29 &2.26E-06& 4.12& 3.30E-06& 4.11&1.51E-08& 3.78&1.04E-07& 4.01 \\                    			 
                  			&0.15 &1.37E-07& 4.46& 2.12E-07& 4.36&9.59E-10& 4.39&6.42E-09& 4.43\\  
                  			&0.08 &8.47E-09 &4.14& 1.32E-08& 4.13&9.52E-11& 3.43&6.28E-10  &3.46 \\  
\hline \hline                  			
 \end{tabular}  }
\end{center}
\caption{History of convergence of the approximation in Example \ref{ex4}.}\label{table:ex_4}
\end{table}
    
\newpage   
\begin{figure}[ht!]
\begin{center}
\includegraphics[scale=0.15]{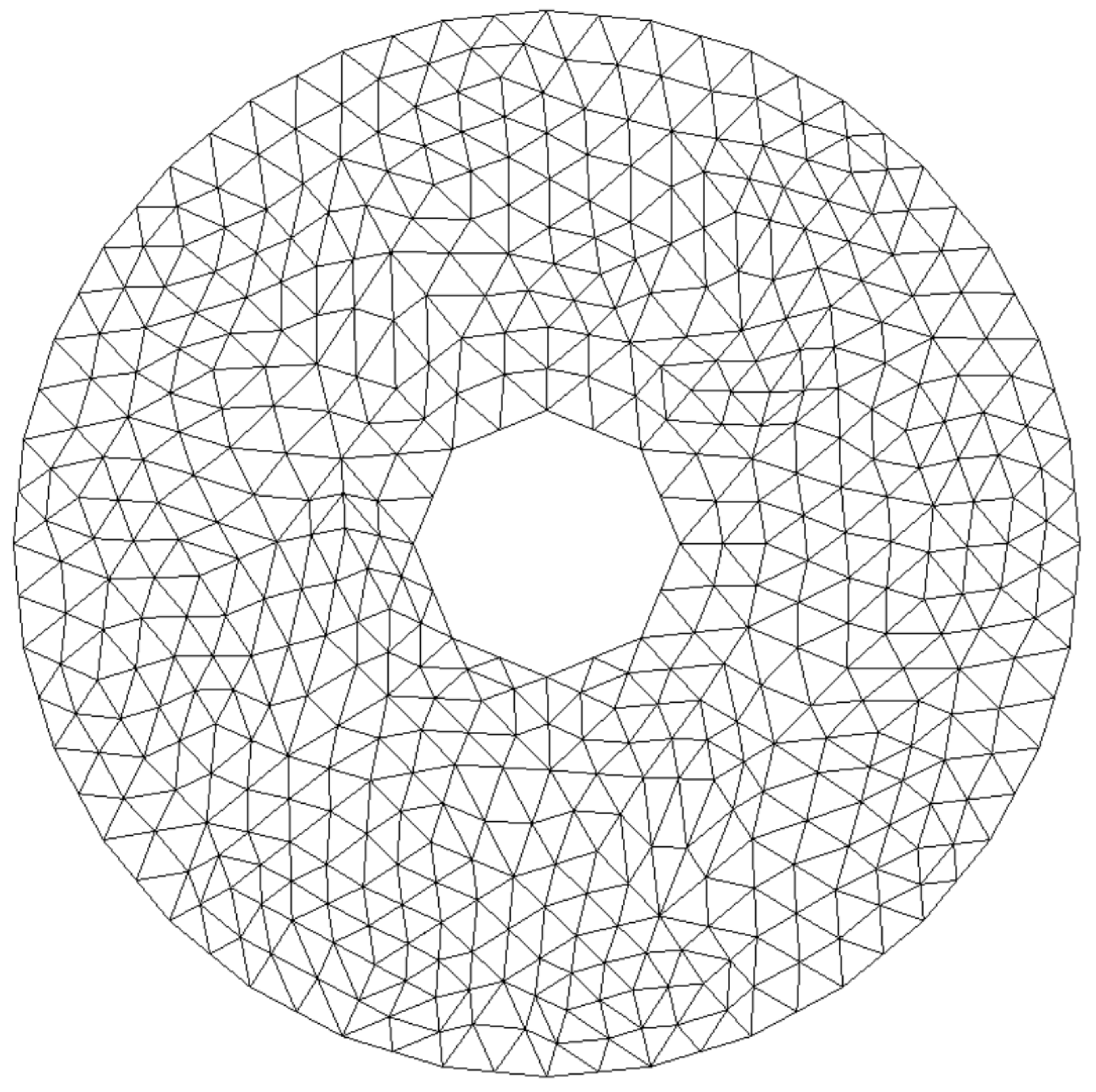}
\includegraphics[scale=0.15]{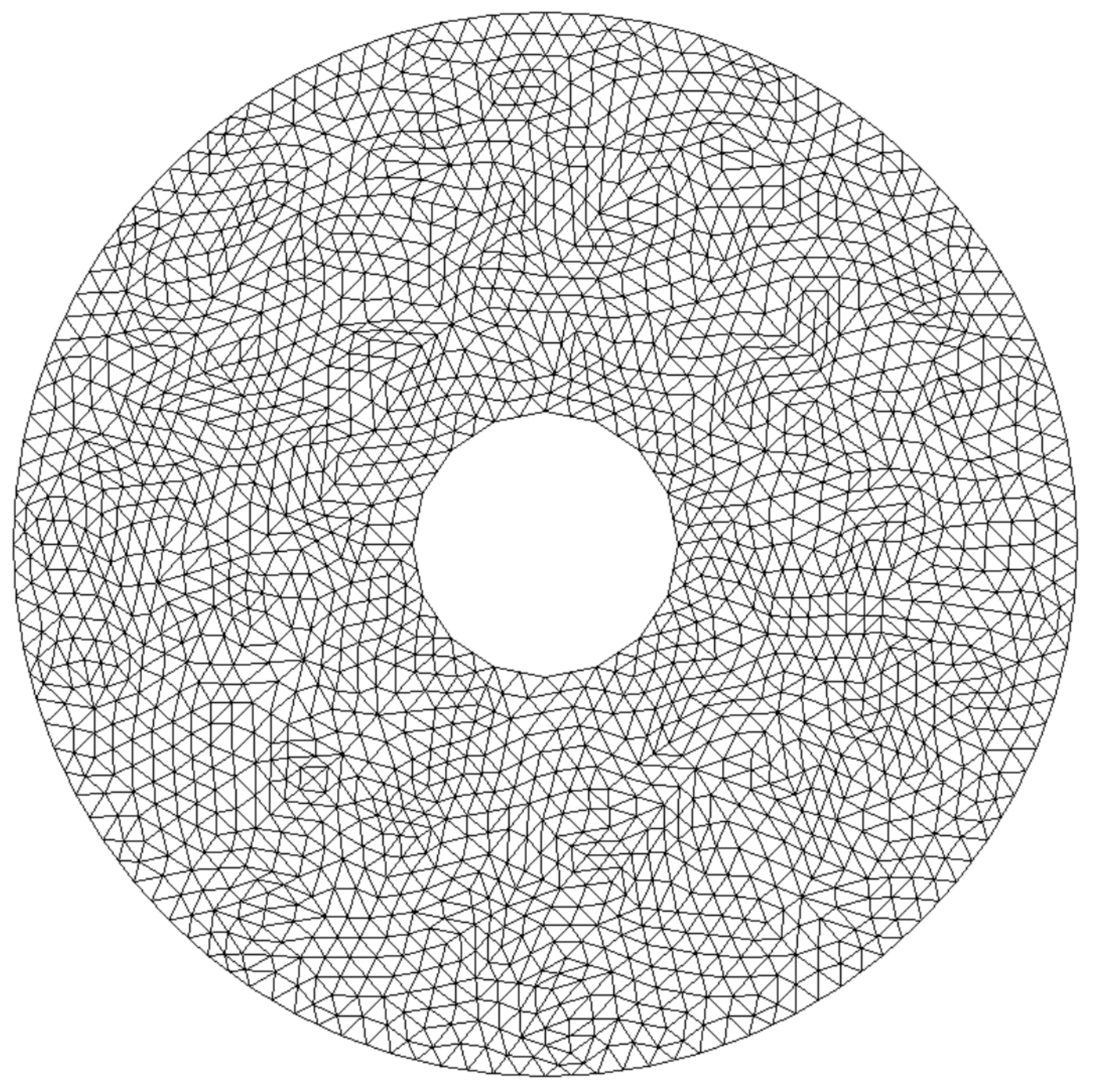}\\
\includegraphics[scale=0.15]{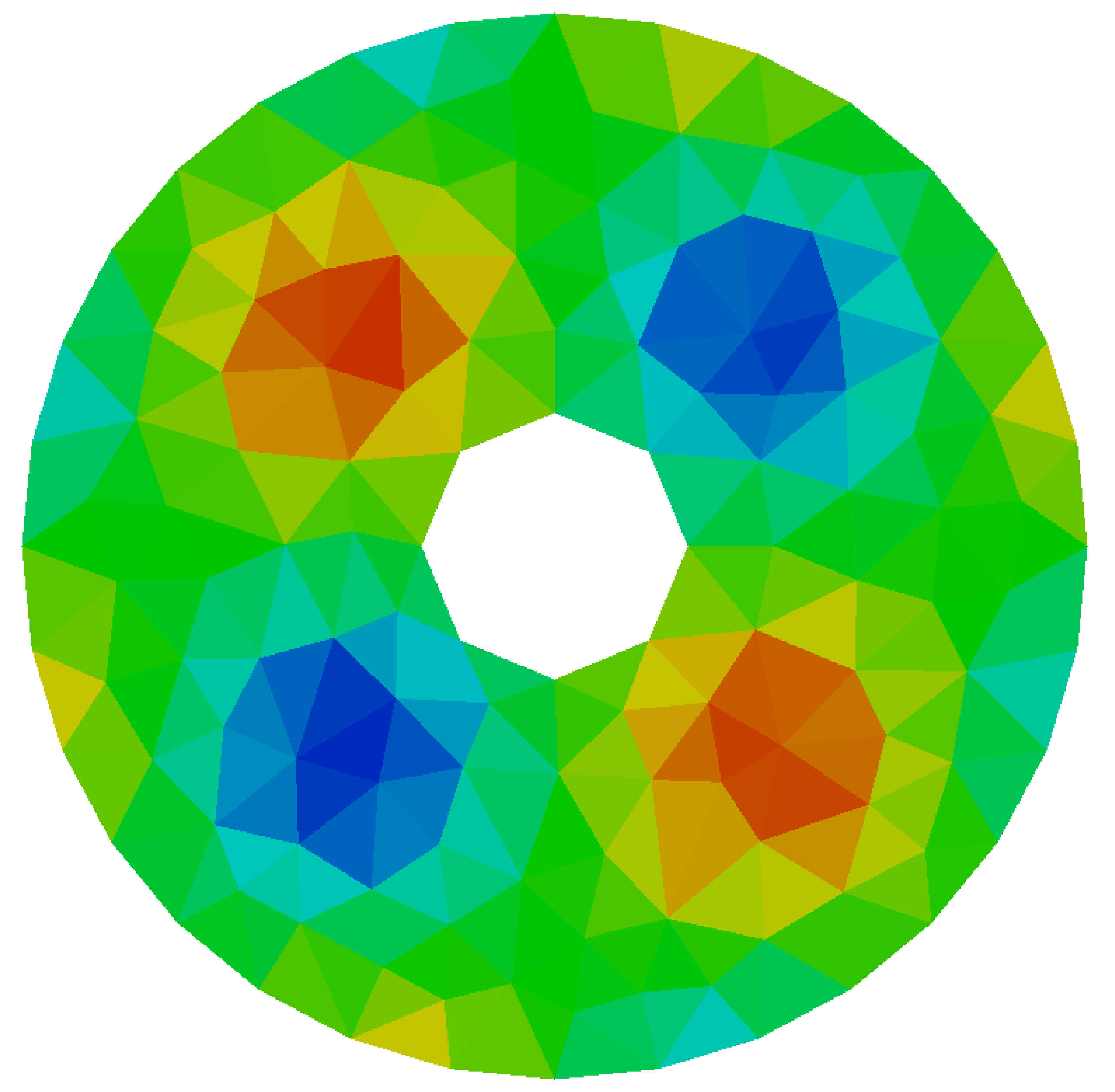}
\includegraphics[scale=0.15]{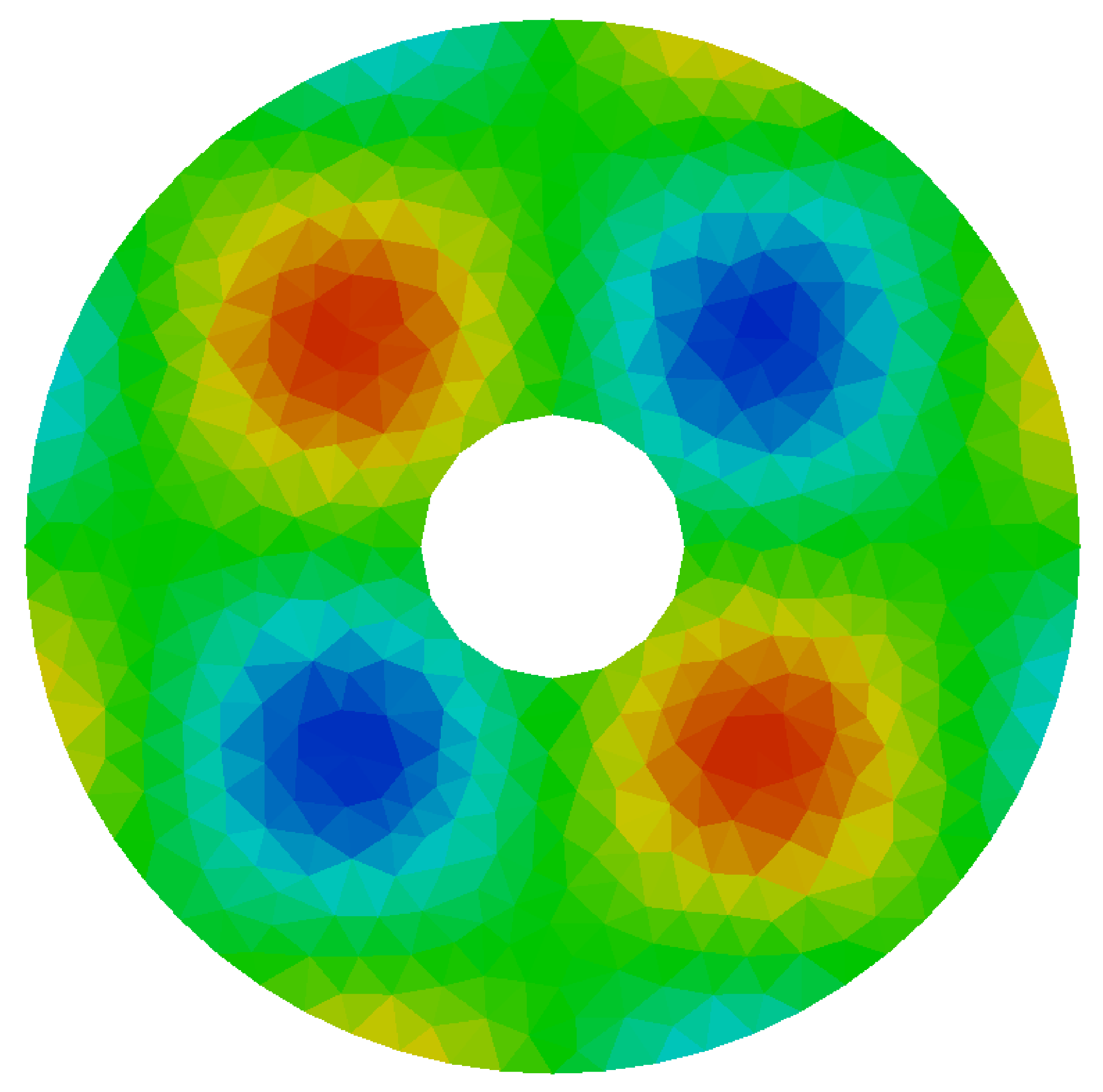}\\
\includegraphics[scale=0.15]{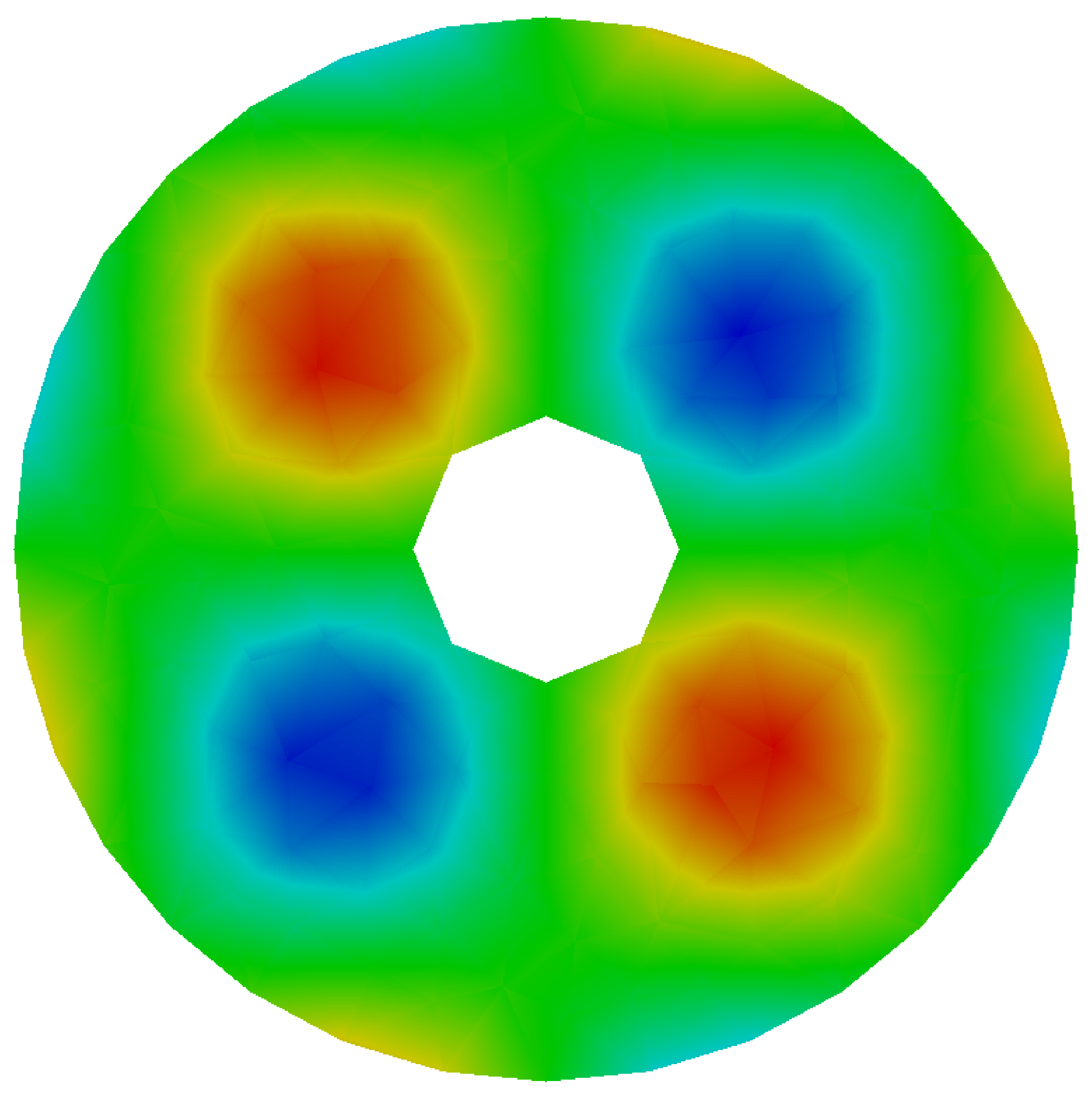}
\includegraphics[scale=0.15]{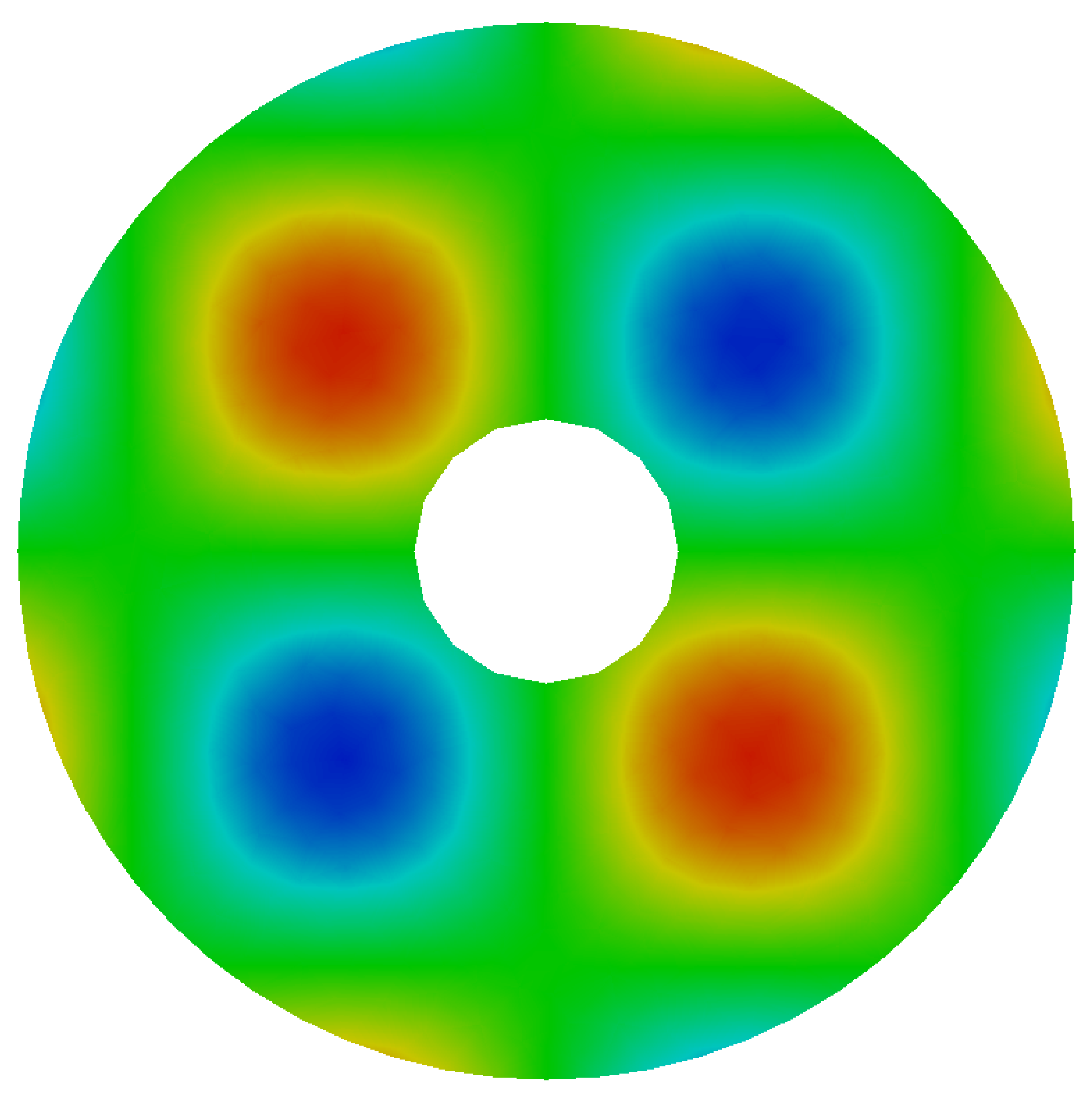}\\
\includegraphics[scale=0.15]{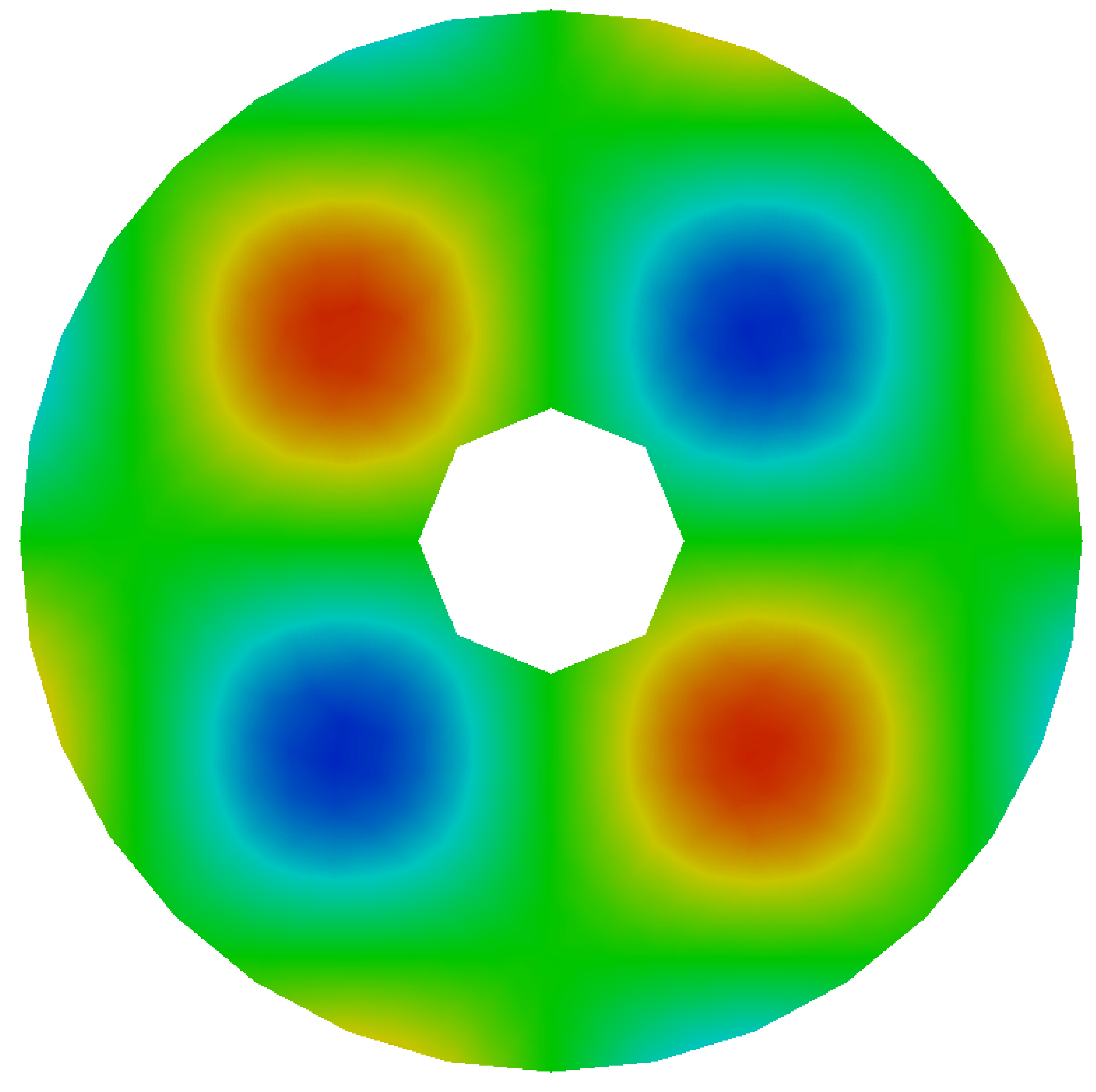}
\includegraphics[scale=0.15]{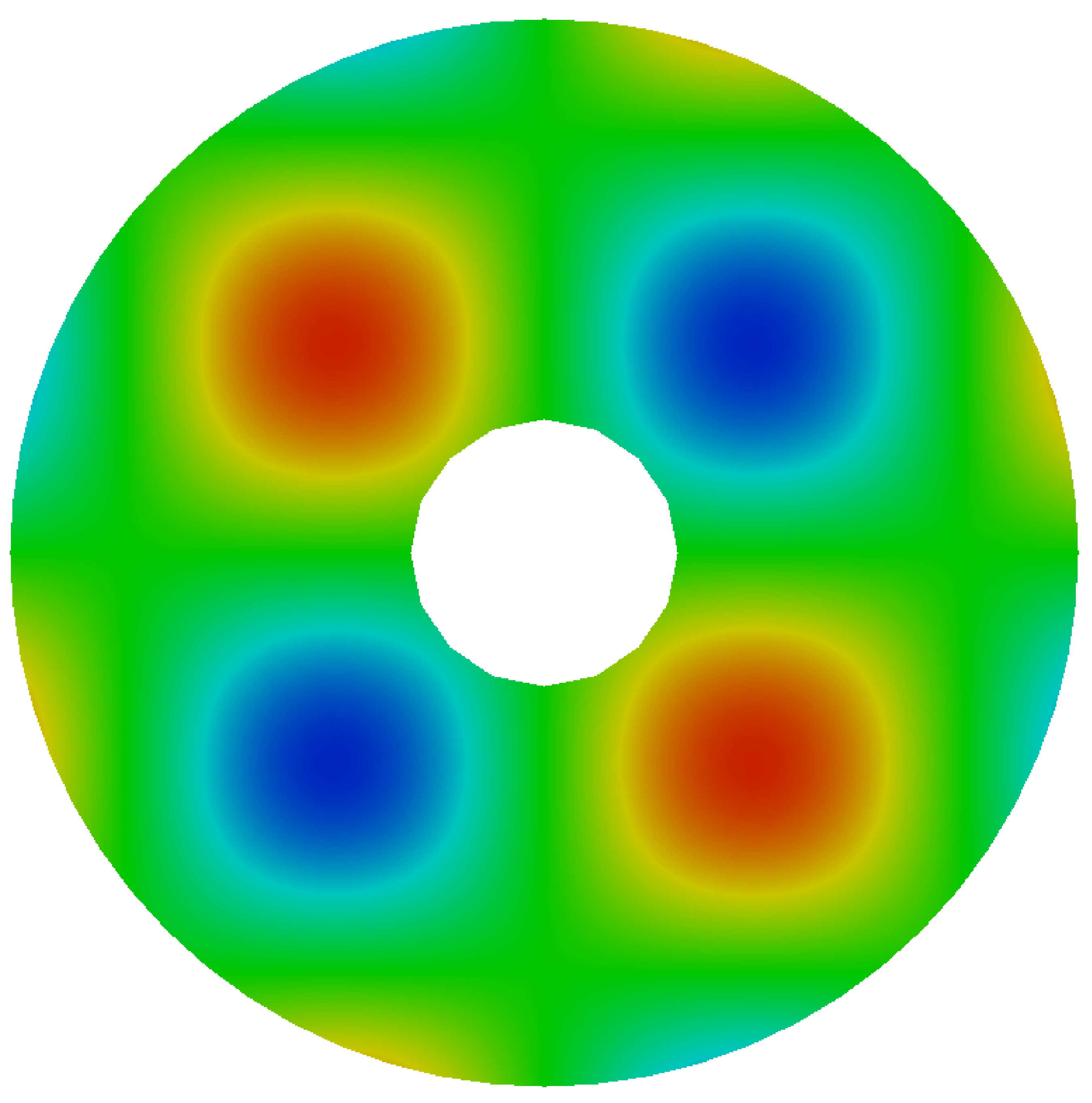}
\caption{Approximation of the scalar variable in Example \ref{ex4}. Columns: meshsize $h=1.10$ and $0.55$. Rows: Polynomial of degree $k=0,1$ and $2$.}
\label{fig2_ex4}
\end{center}
\end{figure}
}
\end{example}
     
\begin{figure}[ht!]
\begin{center}
\includegraphics[scale=0.3]{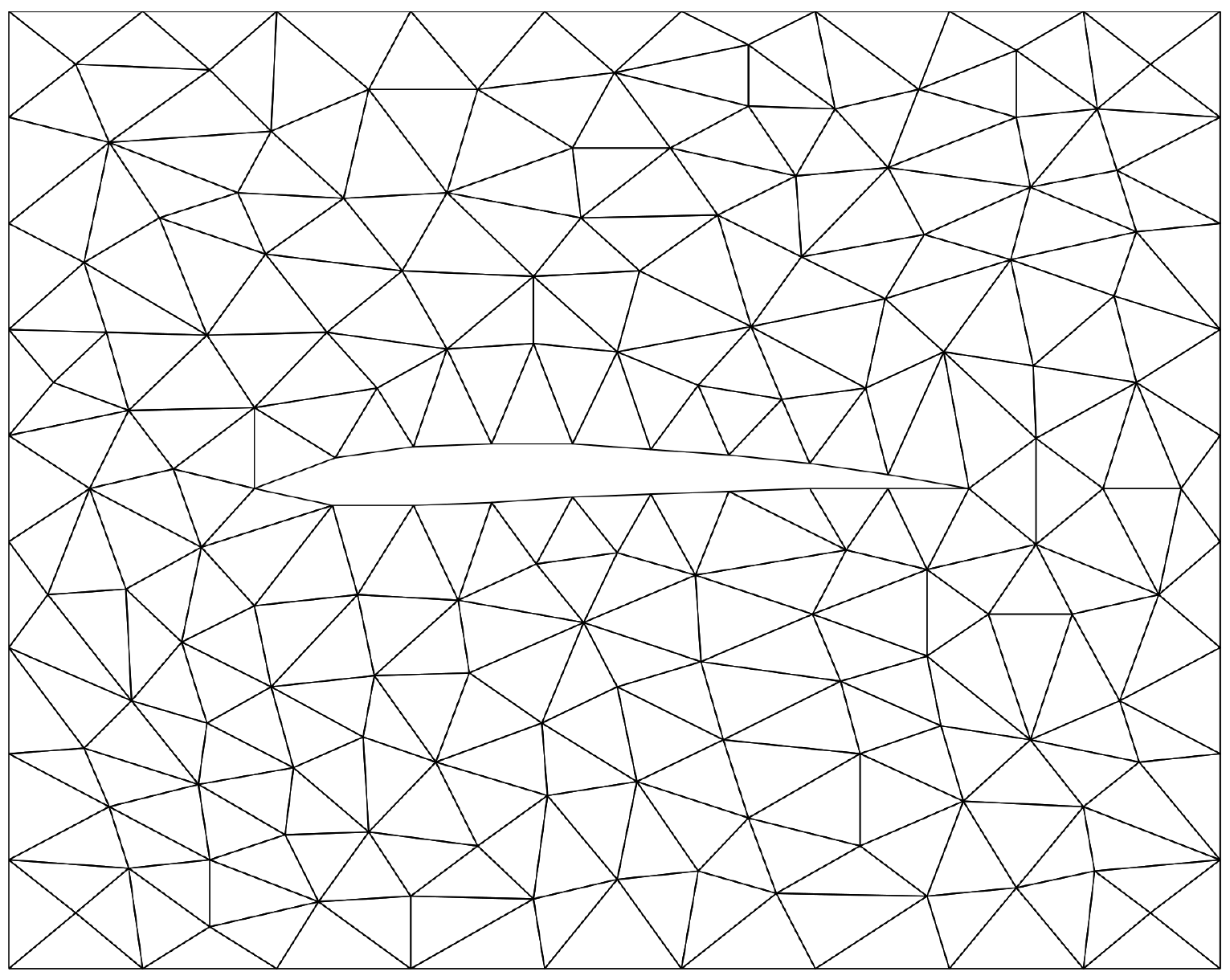}
\includegraphics[scale=0.3]{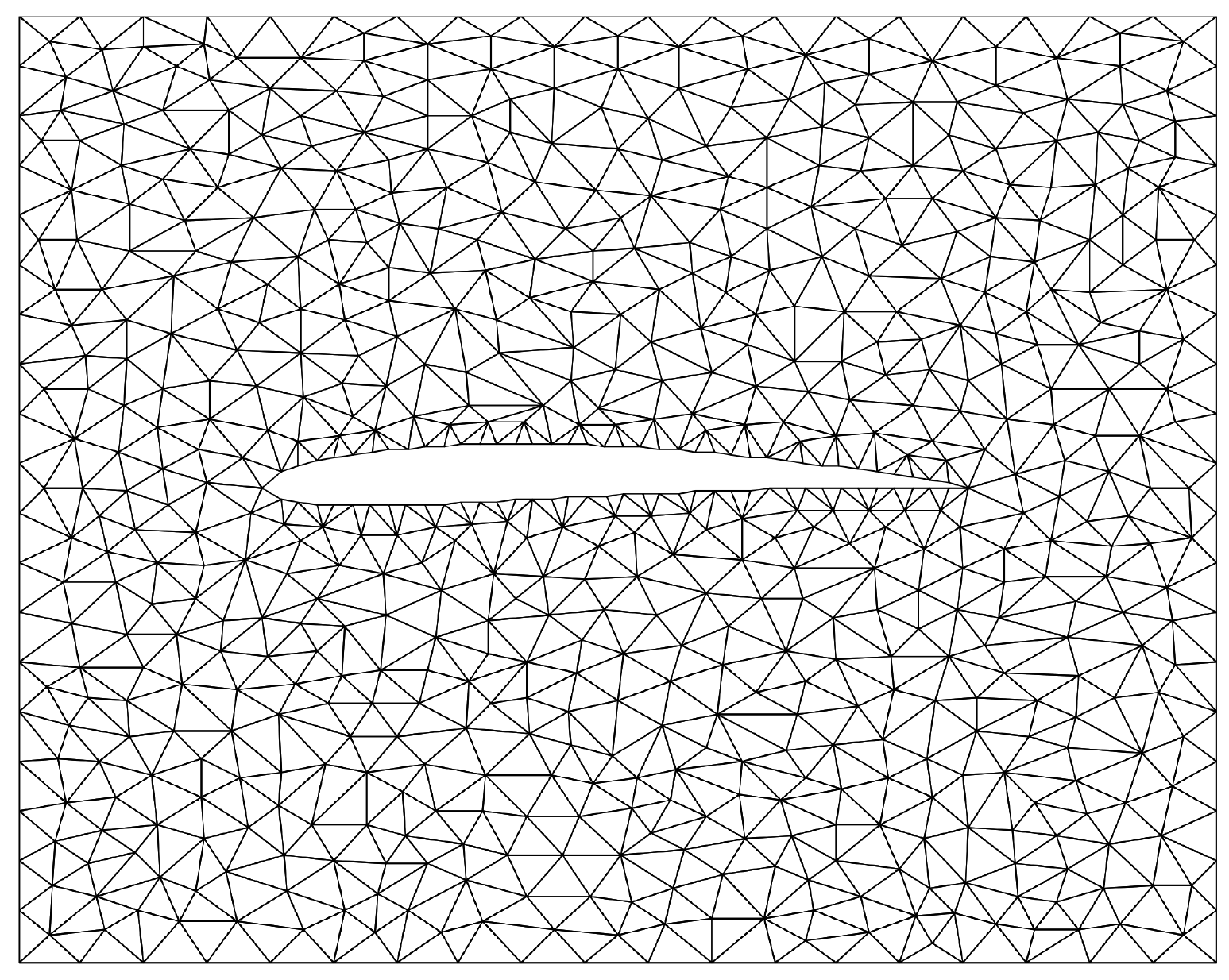}\\
\caption{Meshes of Example \ref{ex5}. Meshsizes $h=0.143$ and $0.073$.}

\label{meshes_ex5}
\end{center}
\end{figure}

\begin{figure}[H]
\begin{center}
\includegraphics[scale=0.15]{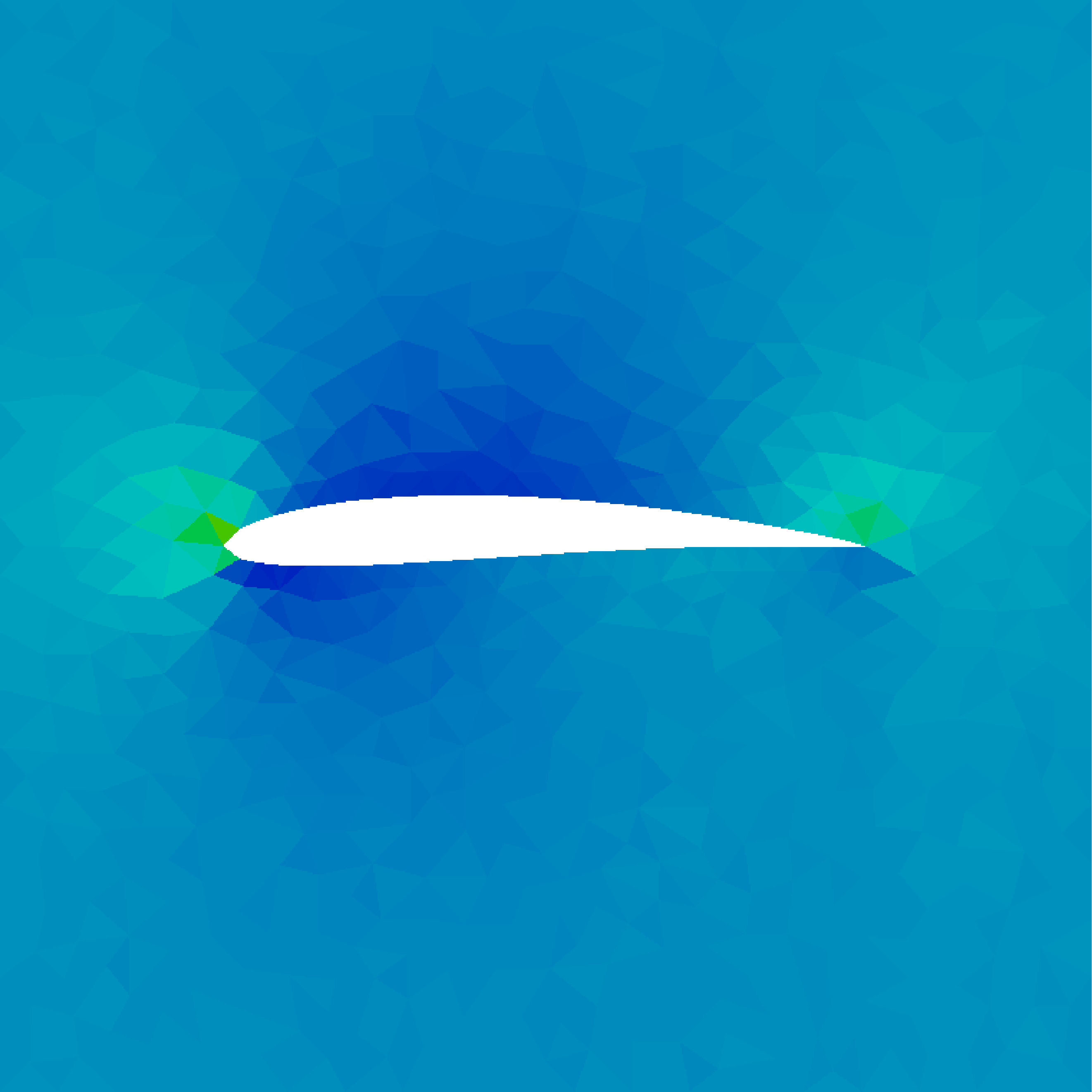}
\includegraphics[scale=0.15]{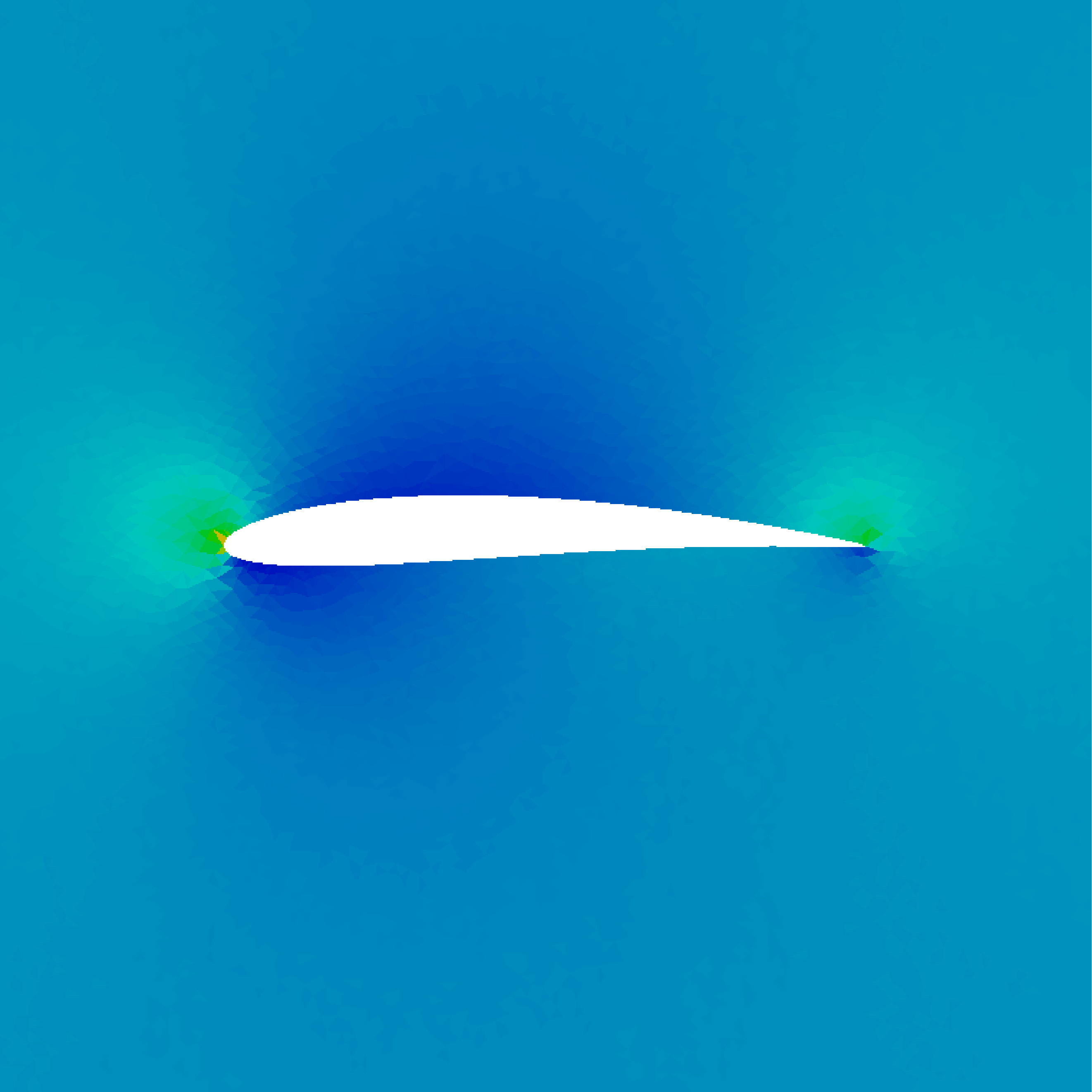}\\
\vspace{0.3cm}
\includegraphics[scale=0.15]{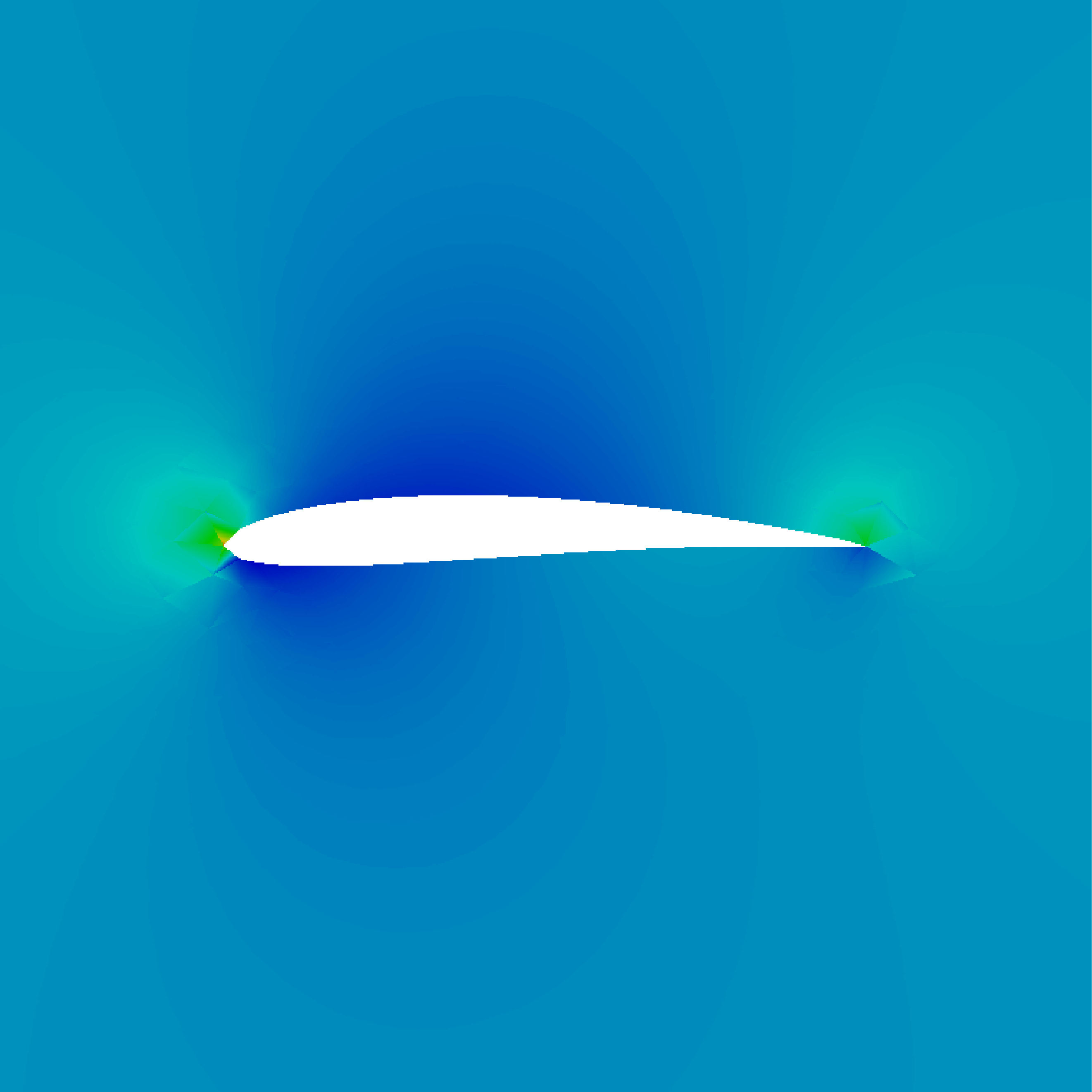}
\includegraphics[scale=0.15]{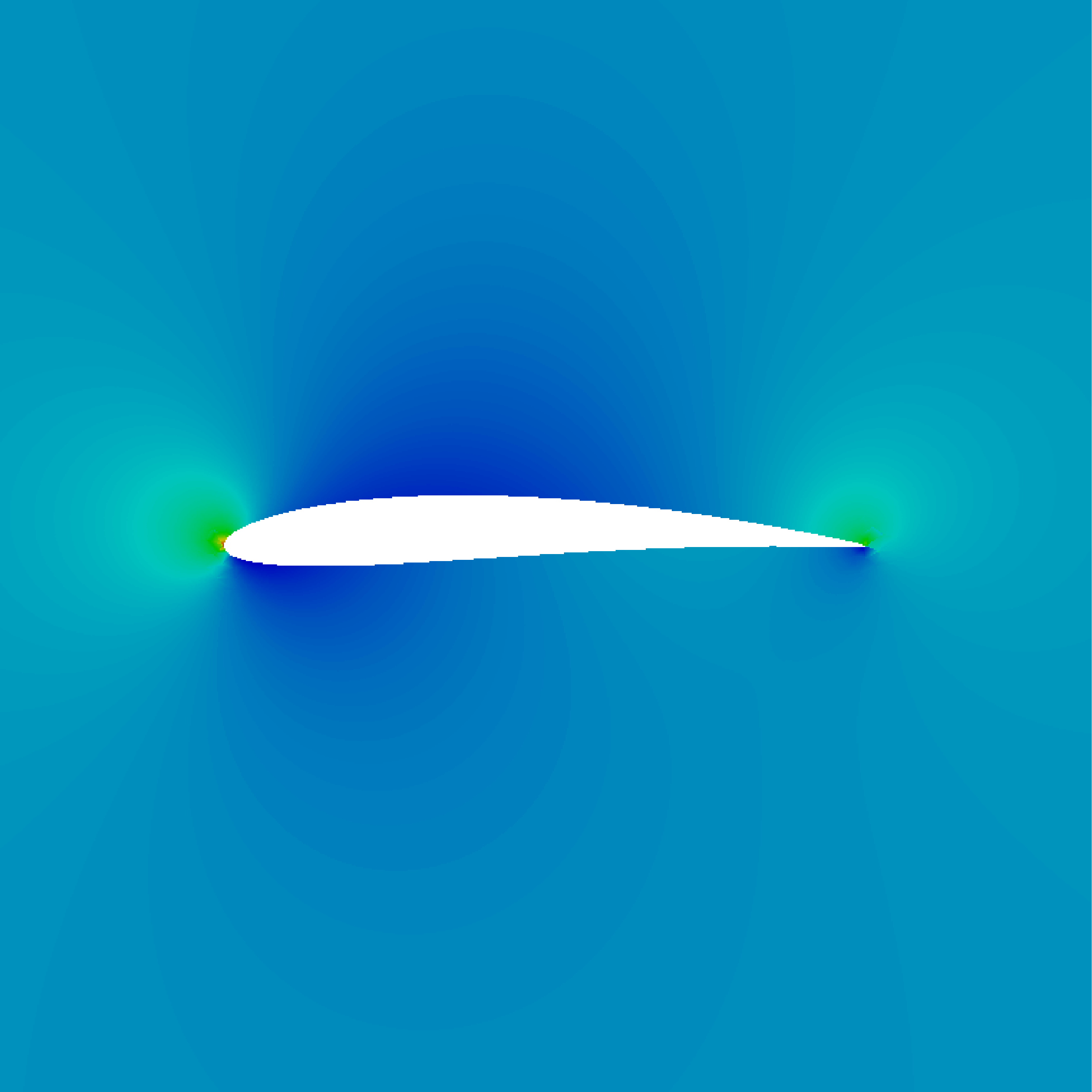}\\
\vspace{0.3cm}
\includegraphics[scale=0.15]{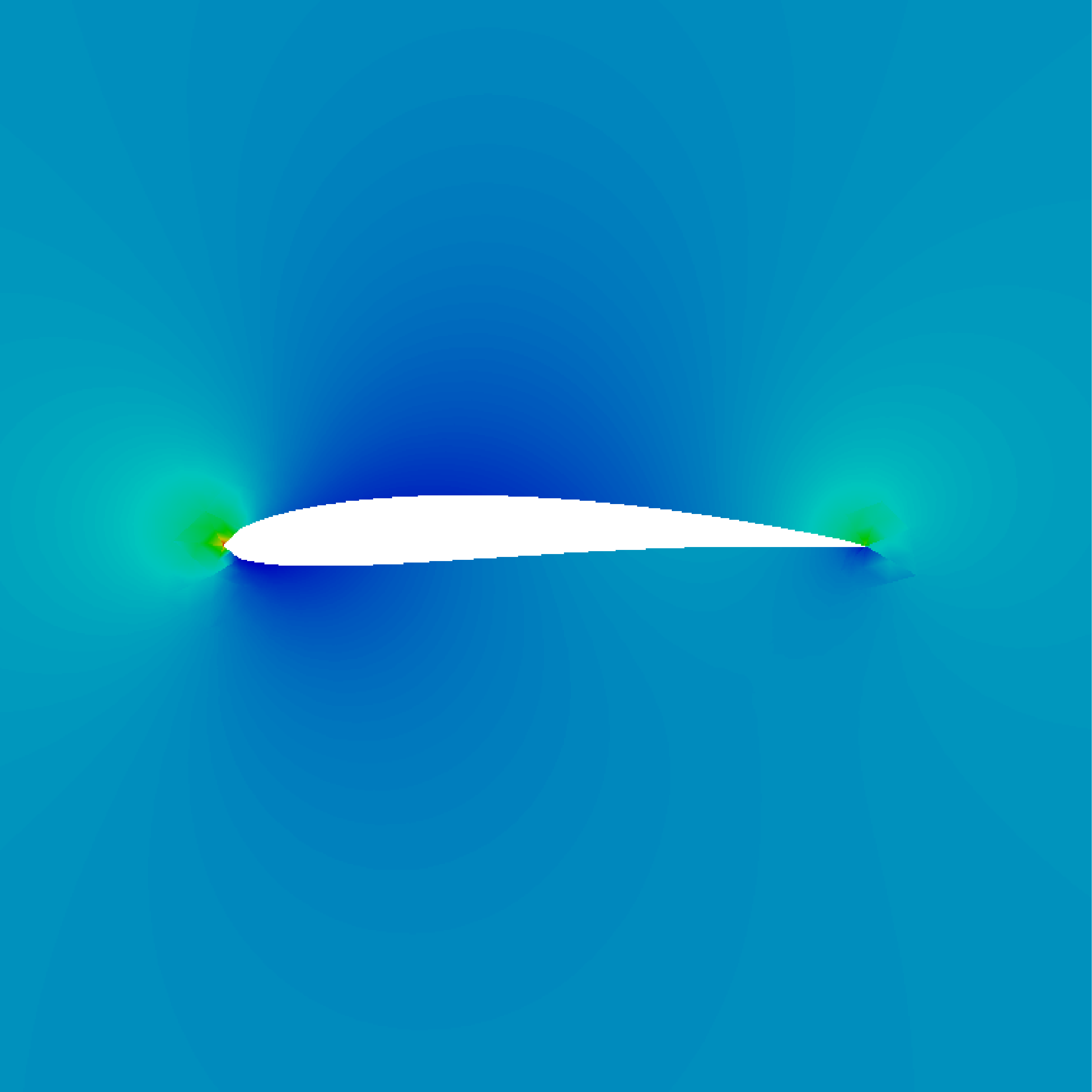}
\includegraphics[scale=0.15]{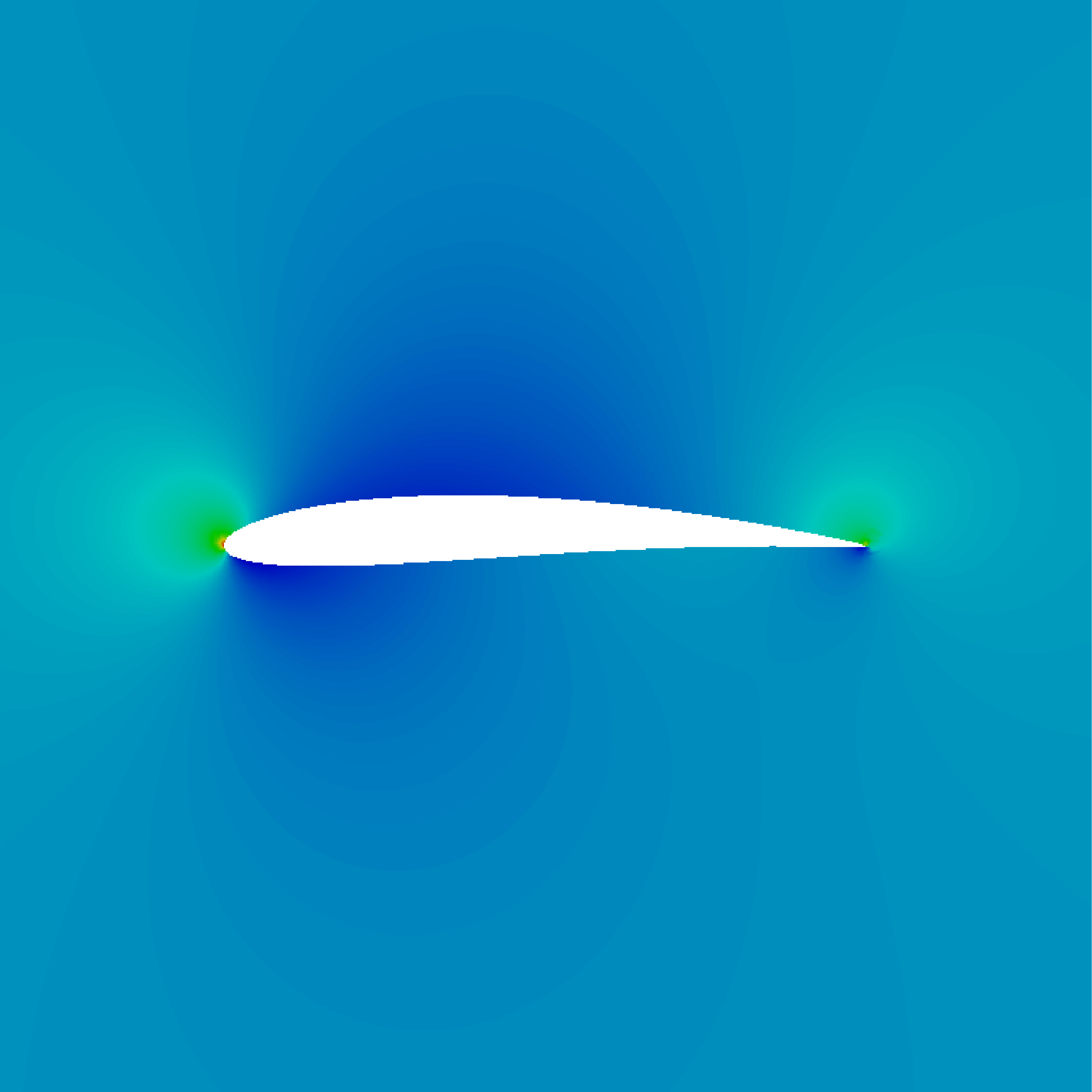}\\
\caption{Approximation of the $x$-component of $\qn$ Example \ref{ex5} (non-smooth solution). Columns: meshsize $h=0.143$ and $0.024$. Rows: Polynomial of degree $k=0,1$ and $2$.}
\label{fig_ex5_nonsmooth}
\end{center}
\end{figure}

\begin{example}\label{ex5}
{\rm Now we test the performance of the method where $\Omega$ is a bounded domain exterior to an airfoil. This is the most difficult case in our examples since the domain has a boundary with a curved, re-entrant corner. The airfoil is obtained by using the Joukowsky
transformation:
\begin{equation*}
J(z)=z + \frac{\lambda^2}{z},
\end{equation*}
where $z\in\mathbb{C}$ and $\lambda \in \mathbb{R}$. It is well known that this transformation maps the disc
centered at $(s_1,s_2)$ of radius $R$ to an airfoil when we set $\lambda=R-\sqrt{s_1^2+s_2^2}$.
Here, we take
$R=0.1605$ and $s_1=s_2=0.01$. In Fig. \ref{meshes_ex5} we show two triangulations of the domain with meshsizes  $h=0.143$ and $0.073$. Neumann boundary conditions are imposed around the airfoil and Dirichlet data in the remaining part of the boundary. \\

We consider the following two examples: \\

\begin{enumerate}
\item[a)]{\bf Smooth solution}. We set $f$ and $g$ such that
$u(x,y)=\sin(x)\sin(y)$ is the exact solution as in previous example.  In Table~\ref{table:ex_5_smooth} we observe that similar conclusions to those in previous examples can be drawn,  even though in the case the domain is more complicated. 
\item[b)]{\bf Non-smooth solution}. We now consider a potential flow around the airfoil where the exact solution in polar coordinates is
$u(r,\theta)=\displaystyle r\cos(\theta) \bigg(1+\frac{R^2}{r^2}\bigg)$. Here $g_N= 0$ around the airfoil. In this case $\nabla u$ has singularities at the leading and trailing edges, hence we do not expect high order convergence rates. In fact, this can be seen on Table~\ref{table:ex_5_nonsmooth} where in all the cases $u$ converges with order one and $\qn$ converges with order less than one. However, for a fixed mesh, the errors decrease when the polynomial degree increases. In Fig. \ref{fig_ex5_nonsmooth} we show the approximation of the $x$-component of $\qn$ considering $h=0.143$ and $0.024$ and $k=0,1$ and $2$.
\end{enumerate}
}
\end{example}

\begin{table}[H]\renewcommand{\arraystretch}{1.3}\addtolength{\tabcolsep}{-5pt}
\begin{center}
{\scriptsize\begin{tabular}{c|c|cc|cc|cc|cc}
  \hline \hline
  \multicolumn{2}{c}{}     &
  \multicolumn{2}{c}{$\|e_u\|_{\ltwoint}$}     &
  \multicolumn{2}{c}{$\|e_{\qn}\|_{\ltwoint}$} &
  \multicolumn{2}{c}{$\|e_{\widehat{u}}\|_{\mathcal{E}_h}$}&
  \multicolumn{2}{c}{$ \|e_{u^*}\|_{\ltwoint}$}\\
  $k$& $h$& error & order &error&order &error &order &error &order \\
  \hline \hline
                		    &0.143&5.69E-03  & -     &2.25E-02 &  -    &1.35E-03& -      &5.76E-03  &- \\
							&0.113&4.78E-03  &0.75 &1.71E-02 & 1.18&7.52E-04 &2.50 &4.81E-03 & 0.77   \\ 
0						    &0.073&3.12E-03  &0.98 &1.05E-02 & 1.11&4.30E-04 &1.29 &3.14E-03 & 0.98   \\ 
                  			&0.038&1.59E-03  &1.03 &5.36E-03 & 1.04&1.97E-04 &1.19 &1.59E-03 & 1.04   \\                  			 
                  			&0.024&9.93E-04  &1.02 &3.25E-03 & 1.08&1.21E-04 &1.06 &9.94E-04 & 1.02   \\ 
  \hline
                		    &0.143&1.41E-04  & -     &2.91E-04 &  -   &1.46E-05 & -      &1.48E-05  &- \\
							&0.113&8.04E-05  &2.38 &1.68E-04 &2.33&8.36E-06 &2.39 &8.46E-06 &2.37   \\ 
1						    &0.073&3.36E-05  &2.01 &6.72E-05 & 2.11&1.95E-06 &3.35 &1.96E-06 & 3.36   \\ 
                  			&0.038&8.51E-06  &2.11 &1.74E-05 & 2.07&5.30E-07 &2.00 &5.14E-07 & 2.05   \\                  			 
                  			&0.024&3.21E-06  &2.11 &6.50E-06 & 2.12&1.32E-07 &3.00 &1.28E-07 & 3.00   \\ 
  \hline
                		    &0.143&1.89E-06  & -     &3.58E-06 &  -      &1.92E-07 & -     &1.85E-07  &- \\
							&0.113&8.56E-07  &3.37 &1.55E-06 & 3.56  &6.58E-08 &4.56 &6.34E-08 & 4.56   \\ 
2					        &0.073&2.27E-07  &3.06 &4.06E-07 & 3.09  &5.65E-09 &5.65 &5.67E-09 & 5.56   \\ 
                  			&0.038&2.96E-08  &3.12 &5.30E-08 & 3.12  &6.17E-10 &3.39 &5.97E-10 & 3.45   \\                  			 
                  			&0.024&6.87E-09  &3.15 &1.24E-08 & 3.14  &7.78E-11 &4.47 &7.57E-11 & 4.45   \\ 
  \hline
                		    &0.143&2.13E-08  & -      &3.00E-08 &  -     &1.04E-08 & -     &9.98E-10  &- \\
							&0.113&7.16E-09  &4.64  &1.06E-08 &4.44  &3.33E-09 &4.86 &3.20E-10 & 4.85   \\ 
3						    &0.073&1.32E-09  &3.89  &1.80E-09 &4.08  &1.89E-10 &6.61 &1.83E-11 & 6.58   \\ 
                  			&0.038&8.65E-11  &4.18  &1.20E-10 &4.14  &1.47E-11 &3.91 &1.40E-12 & 3.95   \\                  			 
                  			&0.024&1.25E-11  &4.17  &1.75E-11 &4.16  &3.52E-12 &3.09 &3.32E-13 & 3.10   \\ 
\hline \hline                  			
 \end{tabular}  }
\end{center}
\caption{History of convergence of the approximation in Example \ref{ex5}a) (smooth solution).}\label{table:ex_5_smooth}
\end{table}

 \begin{table}[H]\renewcommand{\arraystretch}{1.3}\addtolength{\tabcolsep}{-5pt}
\begin{center}
{\scriptsize\begin{tabular}{c|c|cc|cc|cc|cc}
  \hline \hline
  \multicolumn{2}{c}{}     &
  \multicolumn{2}{c}{$\|e_u\|_{\ltwoint}$}     &
  \multicolumn{2}{c}{$\|e_{\qn}\|_{\ltwoint}$} &
  \multicolumn{2}{c}{$\|e_{\widehat{u}}\|_{\mathcal{E}_h}$}&
  \multicolumn{2}{c}{$ \|e_{u^*}\|_{\ltwoint}$}\\
  $k$& $h$& error & order &error&order &error &order &error &order \\
  \hline \hline
                		    &0.143&2.49E-03  & -      &2.20E-02 &  -      &1.40E-03 & -     &2.53E-0  &- \\
							&0.113&1.81E-03  &1.35  &1.62E-02 & 1.29  &7.08E-04 &2.92 &1.84E-03 & 1.36   \\ 
0						    &0.073&1.11E-03  &1.11  &1.10E-02 & 0.90  &2.94E-04 &2.02 &1.12E-03 & 1.14   \\ 
                  			&0.038&5.75E-04  &1.01  &7.23E-03 & 0.64  &1.63E-04 &0.91 &5.77E-04 & 1.02   \\                  			 
                  			&0.024&3.49E-04  &1.08  &5.73E-03 & 0.50  & 9.09E-05&1.26 &3.50E-04 & 1.08   \\ 
  \hline
                		    &0.143&4.04E-04  & -    &8.38E-03 &  -      &4.29E-04  & -    &3.97E-04  &- \\
							&0.113&1.80E-04  &3.45&5.60E-03 & 1.72  &2.08E-04 &3.09 &1.89E-04  & 3.15   \\ 
1						    &0.073&7.93E-05  &1.88&3.38E-03 & 1.16  &8.83E-05 &1.97 &8.07E-05 & 1.96   \\ 
                  			&0.038&4.52E-05  &0.86&2.00E-03 & 0.80  &4.82E-05 &0.93 &4.53E-05 & 0.88   \\                  			 
                  			&0.024&3.03E-05  &0.86&1.63E-03 & 0.45  &3.23E-05 &0.87 &3.03E-05 & 0.87   \\ 
  \hline
                		    &0.143&1.55E-04  & -        &4.37E-03& -     & 1.77E-04 &  -     &1.57E-04 & - \\
							&0.113&8.10E-05  &2.78    &3.02E-03&1.58 &9.16E-05 & 2.81  &8.12E-05 &2.82    \\ 
2					        &0.073&4.91E-05  &1.15    &1.72E-03&1.30 &5.35E-05 & 1.24  &4.91E-05 &1.16    \\ 
                  			&0.038&2.70E-05  &0.92    &9.71E-04&0.87 &2.87E-05 & 0.95  &2.70E-05 &0.92    \\                  			 
                  			&0.024&1.70E-05  &0.99    &8.32E-04&0.33 &1.81E-05 & 0.99  &1.70E-05 &0.99    \\ 
  \hline
                		    &0.143&7.94E-05  & -      &2.73E-03 &  -     &9.13E-05 & -     &8.02E-05  &- \\
							&0.113&4.89E-05  &2.06  &1.84E-03 & 1.68 &5.45E-05 &2.19 &4.92E-05 & 2.08   \\ 
3						    &0.073&3.59E-05  &0.71  &1.09E-03 & 1.21 &3.90E-05 &0.77 &3.60E-05 & 0.72   \\ 
                  			&0.038&1.79E-05  &1.07  &6.34E-04 & 0.83 &1.90E-05 &1.10 &1.79E-05 & 1.07   \\                  			 
                  			&0.024&1.07E-05  &1.10  &5.15E-04 &0.45  &1.14E-05 &1.10 &1.08E-05 & 1.10   \\ 
\hline \hline                  			
 \end{tabular}  }
\end{center}
\caption{History of convergence of the approximation in Example \ref{ex5}b) (Non smooth solution).}\label{table:ex_5_nonsmooth}
\end{table}

\section{Elliptic interface problem}
\begin{figure}[ht!]
\begin{center}
\includegraphics[scale=0.37]{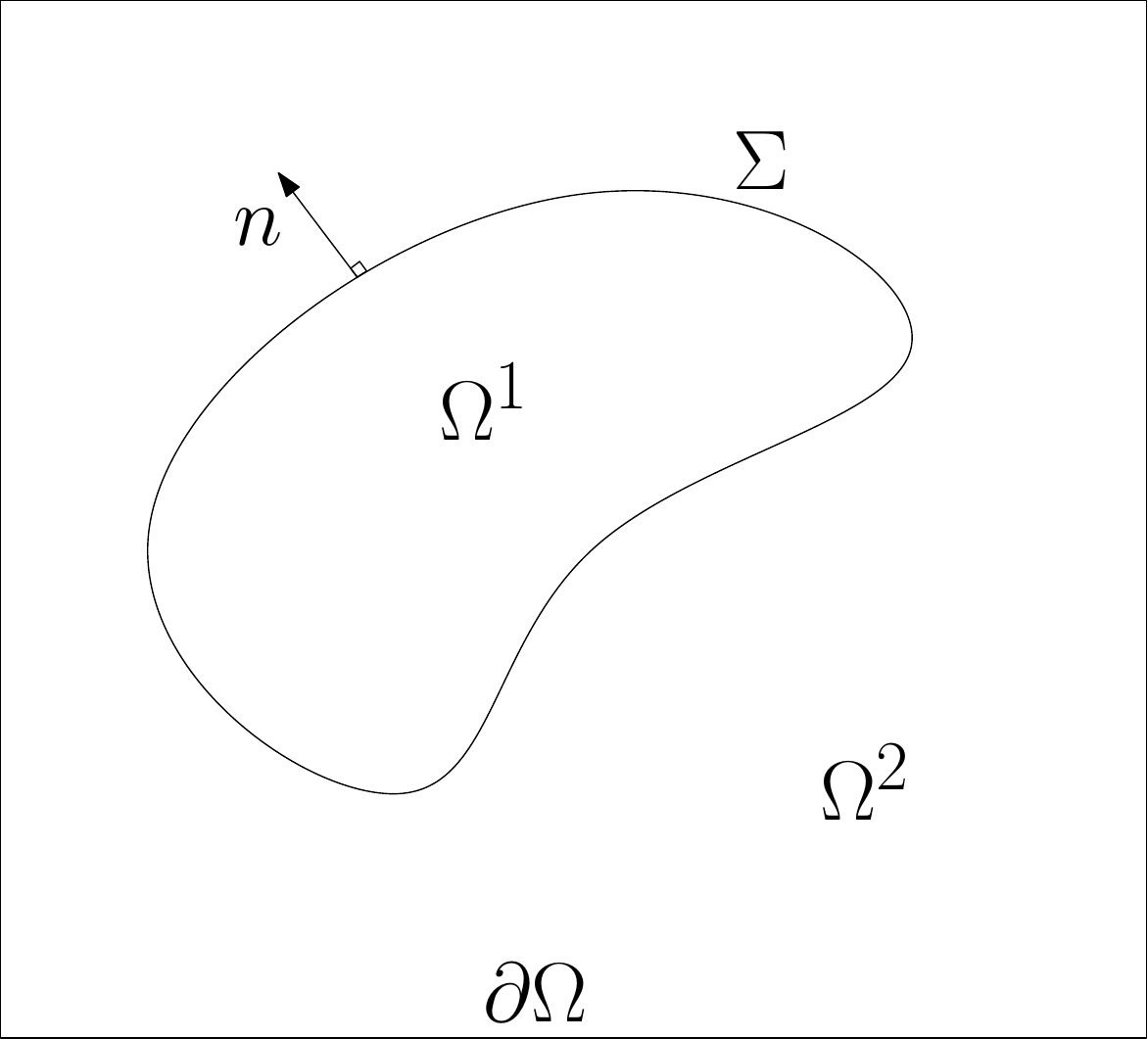}
\caption{Example domain $\Omega$ divided in two regions $\Omega^1$ and $\Omega^2$ by an interface $\Sigma$}
\label{fig_interface}
\end{center}
\end{figure}

Let us now consider and interface $\Sigma$ that divides the domain $\Omega$ in two disjoint subdomains $\Omega^{1}$ and $\Omega^{2}$ as Figure \ref{fig_interface} show. Then, problem (\ref{model_problem}) becomes
\begin{subequations}
\begin{align}
-\divv\bld{q}&=f\quad\text{in }\Omega,\\
\bld{q}+\Kn\nabla u&=0\quad\text{in }\Omega,\\
u&=g_{D}\ \text{on }\Gamma_{D},\\
\bld{q}\cdot\bld{n}&=g_{N}\ \text{on }\Gamma_{N},\\
u\vert_{\Sigma^{1}}-u\vert_{\Sigma^{2}}&=s_{D}\ \text{on }\Sigma,\label{Djump}\\
\bld{q}\vert_{\Sigma^{1}}\cdot\bld{n^{1}}+\bld{q}\vert_{\Sigma^{2}}\cdot\bld{n^{2}}&=s_{N}\ \text{on }\Sigma.
\end{align}\label{model_interface_problem}
\end{subequations}
Here $\Sigma^{1}$ and $\Sigma^{2}$ are defined by
\begin{align*}
\Sigma^{1}&:=\{\xx-\epsilon\bld{n}^{1}:\xx\in\Sigma\text{ and }\epsilon\to 0\},\\
\Sigma^{2}&:=\{\xx-\epsilon\bld{n}^{2}:\xx\in\Sigma\text{ and }\epsilon\to 0\},
\end{align*}
where $\bld{n}^{j}$ ($j \in \{1,2\}$) is the unit outward normal unit vector of the subdomain  $\Omega^j$,  $s_{D}\in H^{1/2}(\Sigma)$ and $s_{N}\in H^{-1/2}(\Sigma)$ are prescribed jumps at the  interface.  At $\Sigma$ we adopt the convention $\bld{n}:=\bld{n}^{1}$.

For the sake of simplicity we assume $\partial \Omega$ to be polygonal (if not, we apply the technique explained in previous section).  However, the interface $\Sigma$ is not necessarily piecewise flat. The numerical results provided in section (\ref{sec:numerical_results}) for a boundary value problem, suggested that the distance between the computational domain and the boundary  should be of order $O(h^2)$  with a family of paths normal to the computational boundary. That is why we interpolate the interface $\Sigma$ by  piecewise linear segments. The computational interface, denoted by $\Sigma_h$, divides the computational domain $\textsf{D}_{h}$ in two disjoint unions of elements $\textsf{D}_{h}^{1}$ and $\textsf{D}_{h}^{2}$. $\Sigma_{h}^{j}$  ($j \in \{1,2\}$) is defined as $\Sigma^{j}_h:=\{\xx-\epsilon\bld{n}_h^{j}:\xx\in\Sigma_h\text{ and }\epsilon\to 0\}$,
where $\bld{n}^{j}_h$ is the unit outward normal  vector of the computational domain   $\textsf{D}_{h}^{j}$.

The main idea is to impose the jump of the scalar variable, denoted by  $\widetilde{s}_{D}^{h}$, on the computational interface $\Sigma_h$. On the other hand, the jump $s_N$ will be imposed at $\Sigma$ by using the idea explained in Section \ref{sec:gNh}.

Following the approach by \cite{MIT1},  the method HDG applied to the interface problem seeks an approximation $(\bld{q}_h,u_h,\lambda_h)  \in \bld{V}_{h}\times W_{h}\times M_{h}$ such that 
\begin{subequations}
\begin{align}
(\Kn^{-1}\bld{q}_{h},\bld{v})_{\textsf{D}_{h}}-(u_{h},\divh\bld{v})_{\textsf{D}_{h}}+\langle \lambda_{h},\bld{v}\cdot\bld{n}\rangle_{\partial \textsf{D}_{h}}&=0,\\
(w,\divh \bld{q}_{h})_{\textsf{D}_{h}}+\langle(\widehat{\bld{q}}_{h}-\bld{q}_{h})\cdot\bld{n},w\rangle_{\partial \textsf{D}_{h}}&=(f,w)_{\textsf{D}_{h}},\\
\langle\widehat{\bld{q}}_{h}\cdot\bld{n},\mu\rangle_{\partial \textsf{D}_{h}\setminus(\Gamma \cup\Sigma_{h})}&=0,\\
\langle \lambda_{h},\mu\rangle_{\Gamma_{D}}&=\langle g_{D},\mu\rangle_{\Gamma_{D}},\label{discDc}\\
\langle \widehat{\bld{q}}_{h}\cdot\bld{n},\mu\rangle_{\Gamma_{N}}&=\langle g_{N},\mu\rangle_{\Gamma_{N}},\label{discNc}
\end{align}
for all $(\bld{v},w,\mu)\in\bld{V}_{h}\times W_{h}\times M_{h}$. We still need to specify the jump of the normal component of $\qn$ at $\Sigma$.\\

Here $\lambda_h$ is a single-valued function, however the approximation of $u$ must be  double-valued on $\Sigma_h$. Then, similarly to \cite{MIT1}, we let  $\lambda_h$ be the approximation of $u_h|_{\Sigma_h^2}$ and consider $\lambda_h+\widetilde{s}_{D}^h$ as an approximation of $u_h|_{\Sigma_h^1}$. Thus, we define 
\begin{equation}
\widehat{u}_{h}:=\lambda_{h}+\delta_{\Sigma_{h}}\widetilde{s}^{h}_{D},\label{numfluxu}
\end{equation} 
\label{discreteproblem}
\end{subequations}
where  $\delta_{\Sigma_{h}}$,  defined on $\partial \textsf{D}_{h}$, satisfies
\begin{equation}
\delta_{\Sigma_{h}}=\begin{cases}1\quad\text{on }\partial K\cap\Sigma_{h},\text{ if }\partial K\cap\Sigma_{h}\neq\emptyset\text{ and }K\in \textsf{D}_{h}^{1},\\ 0\quad\text{otherwise}.\end{cases}\label{sqitchfunc}
\end{equation}

To complete the method we define the numerical flux as usual
$$
\widehat{\bld{q}}_{h} : = \widetilde{\bld{q}}_h -\tau (u_h - \widehat{u}_h)\bld{n} \quad\textrm{on}\quad \partial \textsf{D}_h.
$$

\subsection{Approximation $\widetilde{s}_{D}^h$}
In order to define an approximation of $\widetilde{s}_{D}^h$, we use the same {\it transferring} technique used  for the Dirichlet data on a curved boundary (\ref{new_bc}). Let $e\subset\Sigma_{h}$ such that $e=\partial K_1 \cap \partial K_2$ and, without loss of generality, assume that $e$ lies completely inside of $\Omega^2$. We denote by $(\qn_h^j,u_h^j)$ the approximation $(\qn_h,u_h)$ restricted to the domain $\textsf{D}_h^j$. Now, for each $\xx \in e$, we observe that $\sigma(\xx)\subset K_1\cap\Omega^2$ and then, according to the approximation given in (\ref{gh_D}),

\begin{eqnarray}
u_h^2(\xx)\approx u_h^2(\bar{\xx})+\int_{\sigma(\xx)}\Kn^{-1}E^{K_2}(\bld{q}_h^2)\cdot\bld{m},
\end{eqnarray}
where $E^{K_2}(\bld{q}_h^2)$ is the standard extrapolation of $\bld{q}_h^2$ to the whole $\mathbb{R}^2$ space defined in (\ref{DBC}). Similarly,
\begin{eqnarray}
u_h^1(\xx)\approx u_h^1(\bar{\xx})+\int_{\sigma(\xx)}\Kn^{-1}E^{K_1}(\bld{q}_h^1)\cdot\bld{m},
\end{eqnarray}
In this case $E^{K_1}(\bld{q}_h^1) =\bld{q}_h^1$.

Combining both equations,
\begin{align*}
u_h^1(\xx) - u_h^2(\xx)\approx u_h^1(\bar{\xx})- u_h^2(\bar{\xx})
+\int_{\sigma(\xx)}\Kn^{-1}E^{K_1}(\bld{q}_h^1)\cdot\bld{m}
-\int_{\sigma(\xx)}\Kn^{-1}E^{K_2}(\bld{q}_h^2)\cdot\bld{m}.
\end{align*}

This expression suggest the  following approximation
\begin{eqnarray}
 s^{h}_{D}(\bx):=s_{D}(\bar{\xx})
+\int_{\sigma(\xx)}\Kn^{-1}E^{K_1}(\bld{q}_h^1)\cdot\bld{m}
-\int_{\sigma(\xx)}\Kn^{-1}E^{K_2}(\bld{q}_h^2)\cdot\bld{m}.
\end{eqnarray}

\subsection{Imposition of  $s_{N}$}
For approximating $s_{N}$ we use the same idea that we applied for a Neumann boundary edge.  For each interface edge $e\in\Sigma_{h}$, we consider $\Sigma_{e}\subset\Sigma$, the part of $\Sigma$ associated to $e$. We denote by $K_{e}^{1}$ and $K_{e}^{2}$ the element of $\textsf{D}_{h}^{1}$ and $\textsf{D}_{h}^{2}$  where $e$ belongs. Then, we impose the following condition at the interface $\Sigma$:

\begin{eqnarray}
\langle E^{K_e^1}(\bld{q}_h)\cdot\bld{n}^{1}+E^{K_e^2}(\bld{q}_h)\cdot\bld{n}^{2},\mu\rangle_{\Sigma_{e}}=\langle s_{N},\mu\rangle_{\Sigma_{e}}\quad\forall\mu\in M_{\boldsymbol{\phi}}(\Sigma_{e}),\label{approxdiscNjump}
\end{eqnarray}
where $M_{\boldsymbol{\phi}}(\Sigma_{e})$ is defined similarly as in (\ref{spaceM}).

\subsection{Numerical results: Interface problem}
\label{sec:numerical_results_interface}
Finally, in this section we consider three numerical examples showing the performance of our technique in elliptic interface problems.  
Since the computational domains $\textsf{D}_{h}^{1}$ and $\textsf{D}_{h}^{2}$ do not exactly fit $\Omega^1$ and $\Omega^2$, we exclude from the computation of the errors the triangles intersecting the interface. Let $\widetilde{\textsf{D}}_h$ the set of triangles whose faces are not interface edges. We measure the errors using the following norms $
\|\cdot\|_{\ltwotilde}$ and
\begin{eqnarray*}
\|e_{\hat{u}}\|_{\tildeE}:&=&\left(\sum_{K\in
\widetilde{\textsf{D}}_h: K\cap \Sigma_h = \emptyset} h_K \|\textsf{P}_{\partial} u - \hat{u}_h
\|^2_{L^2(\partial K)}\right)^{1/2}.
\end{eqnarray*}

\begin{figure}[H]
\begin{center}
\includegraphics[width=4.5cm]{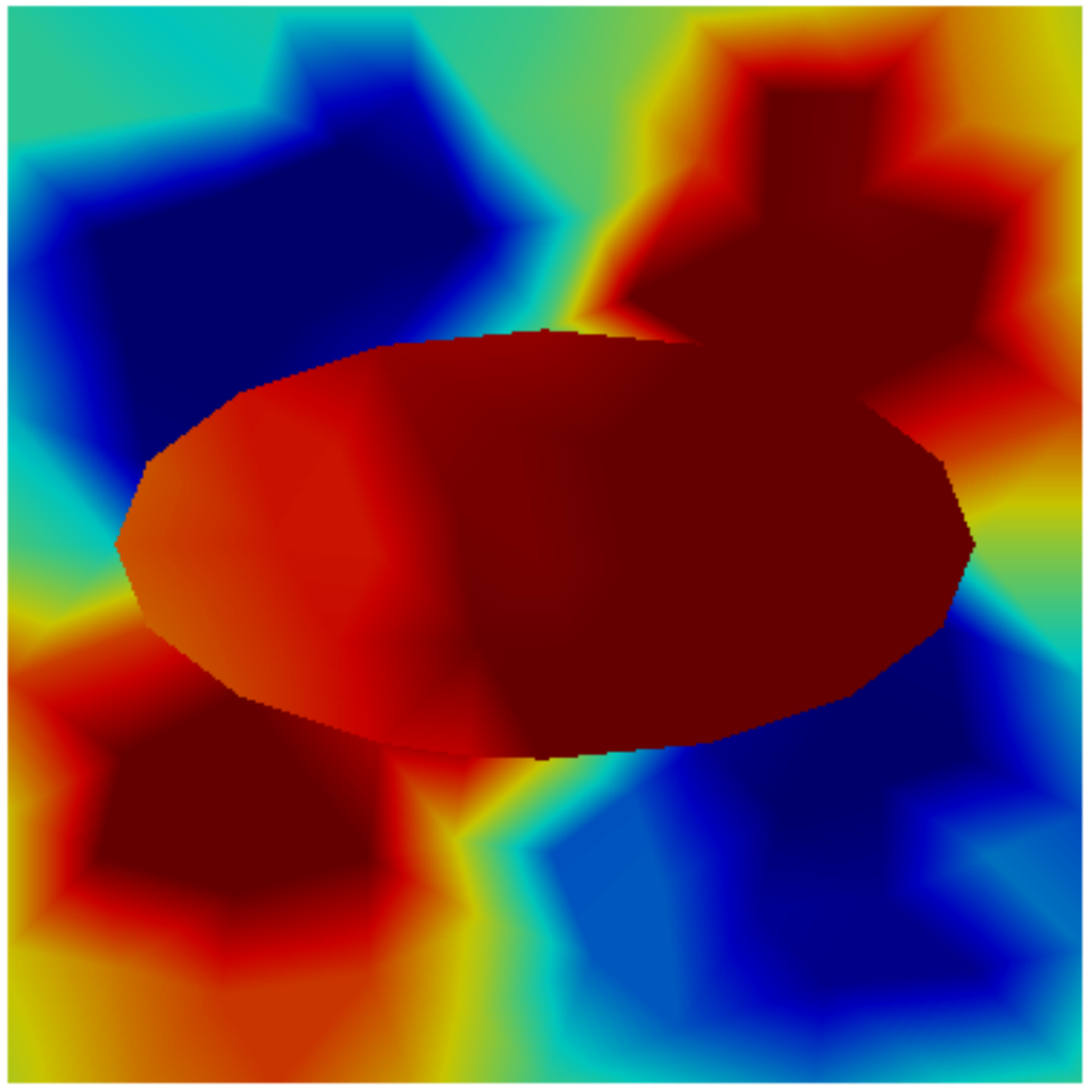}
\includegraphics[width=4.5cm]{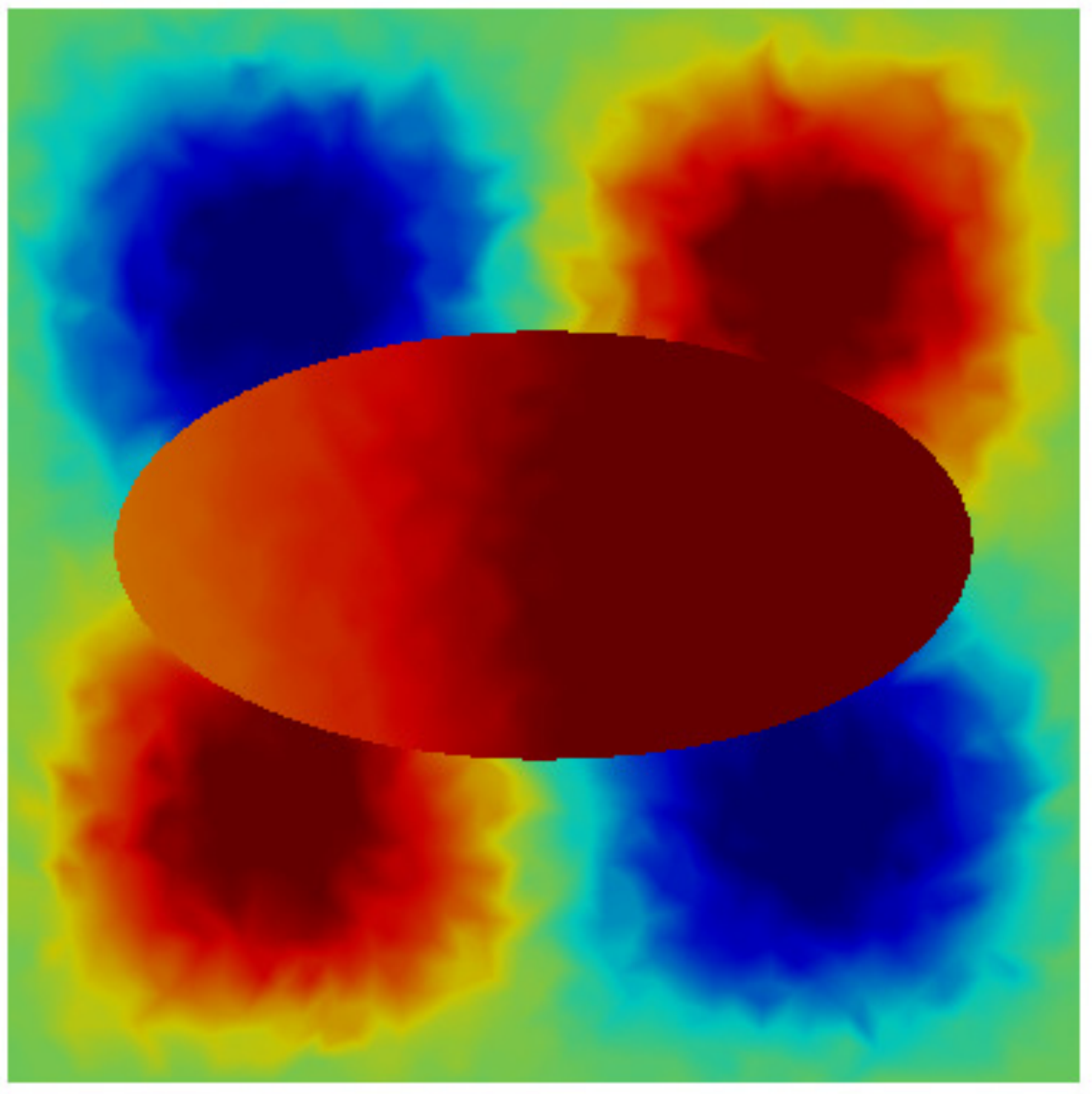}\\
\includegraphics[width=4.5cm]{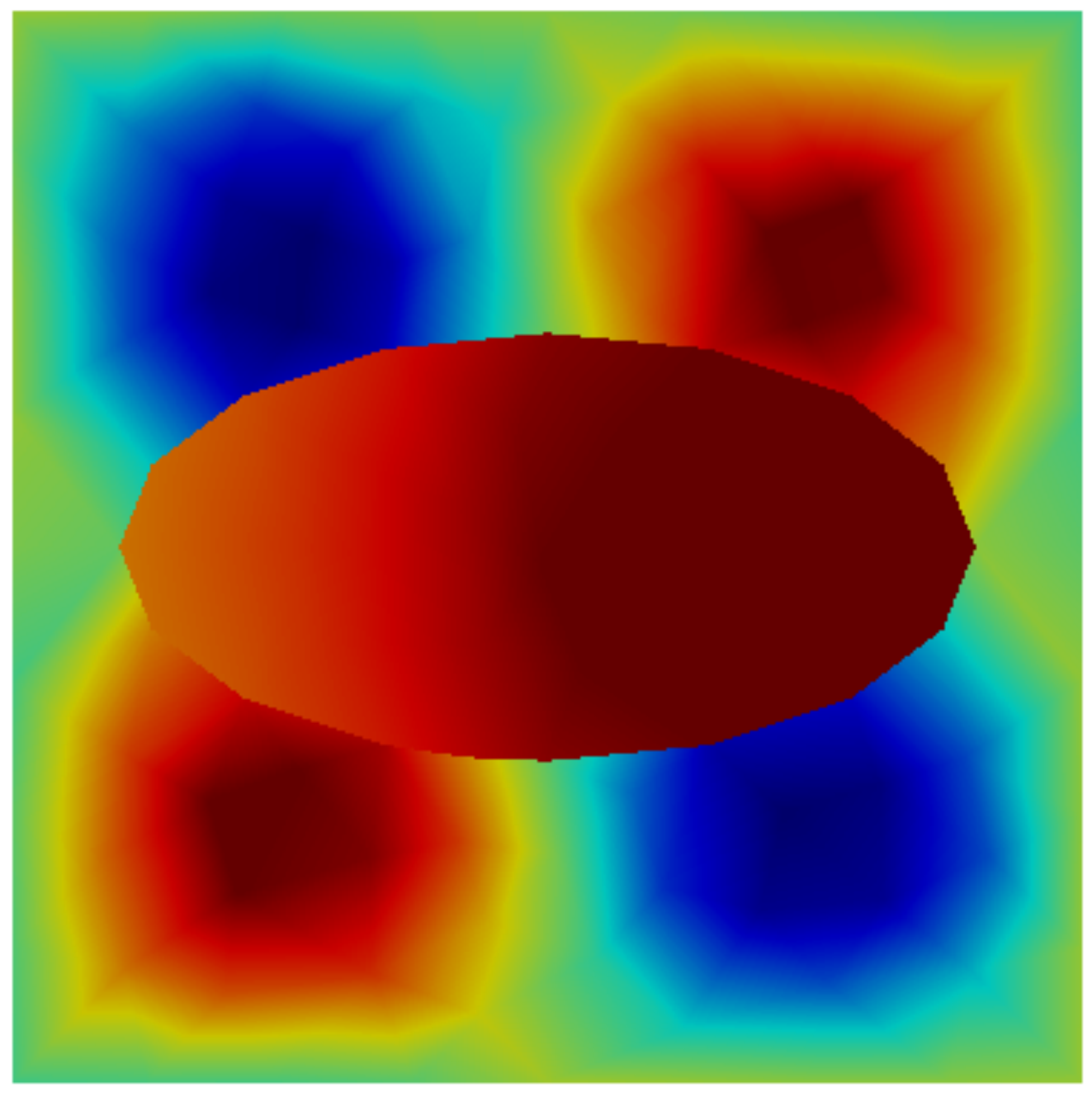}
\includegraphics[width=4.5cm]{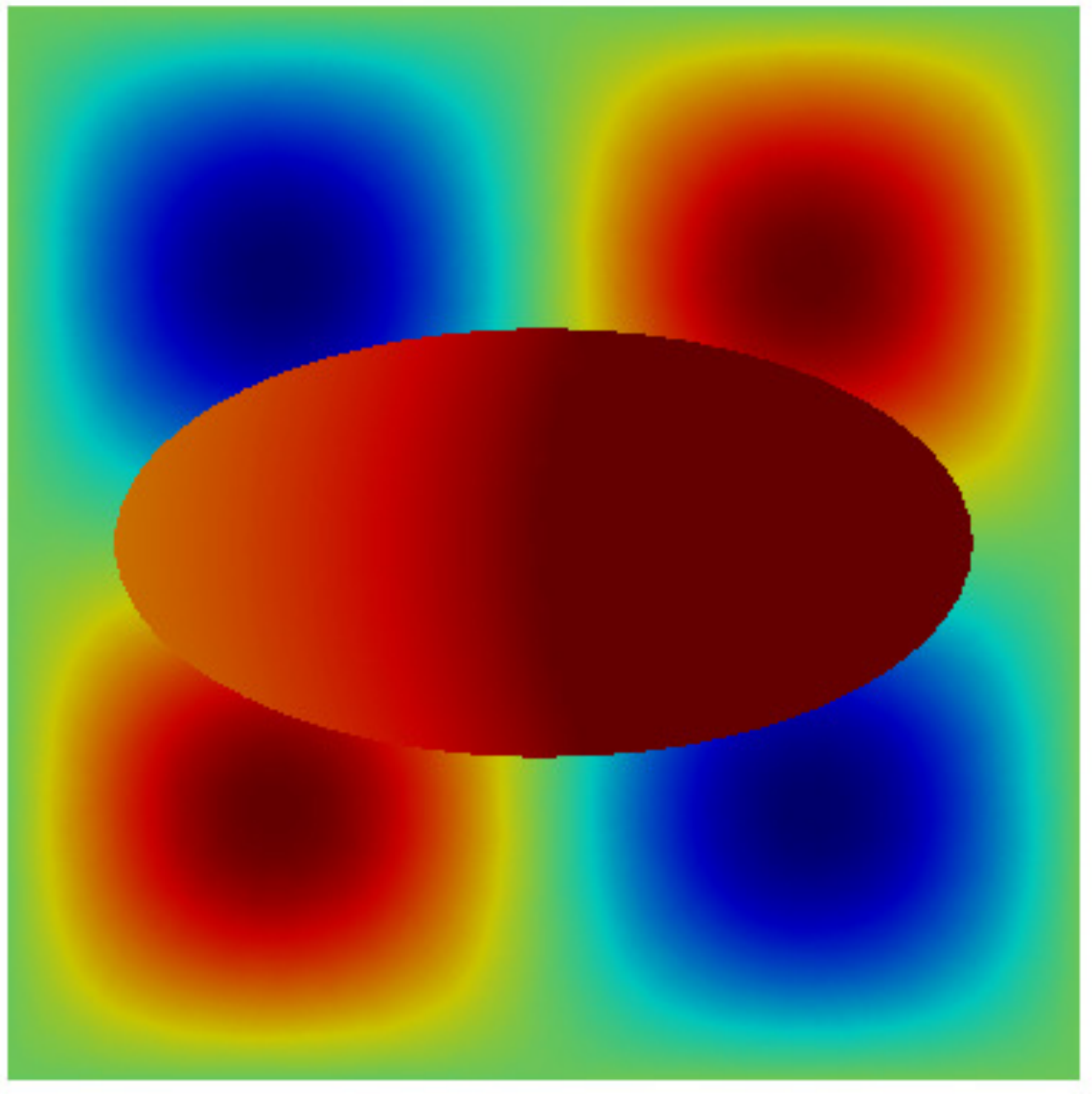}\\
\includegraphics[width=4.5cm]{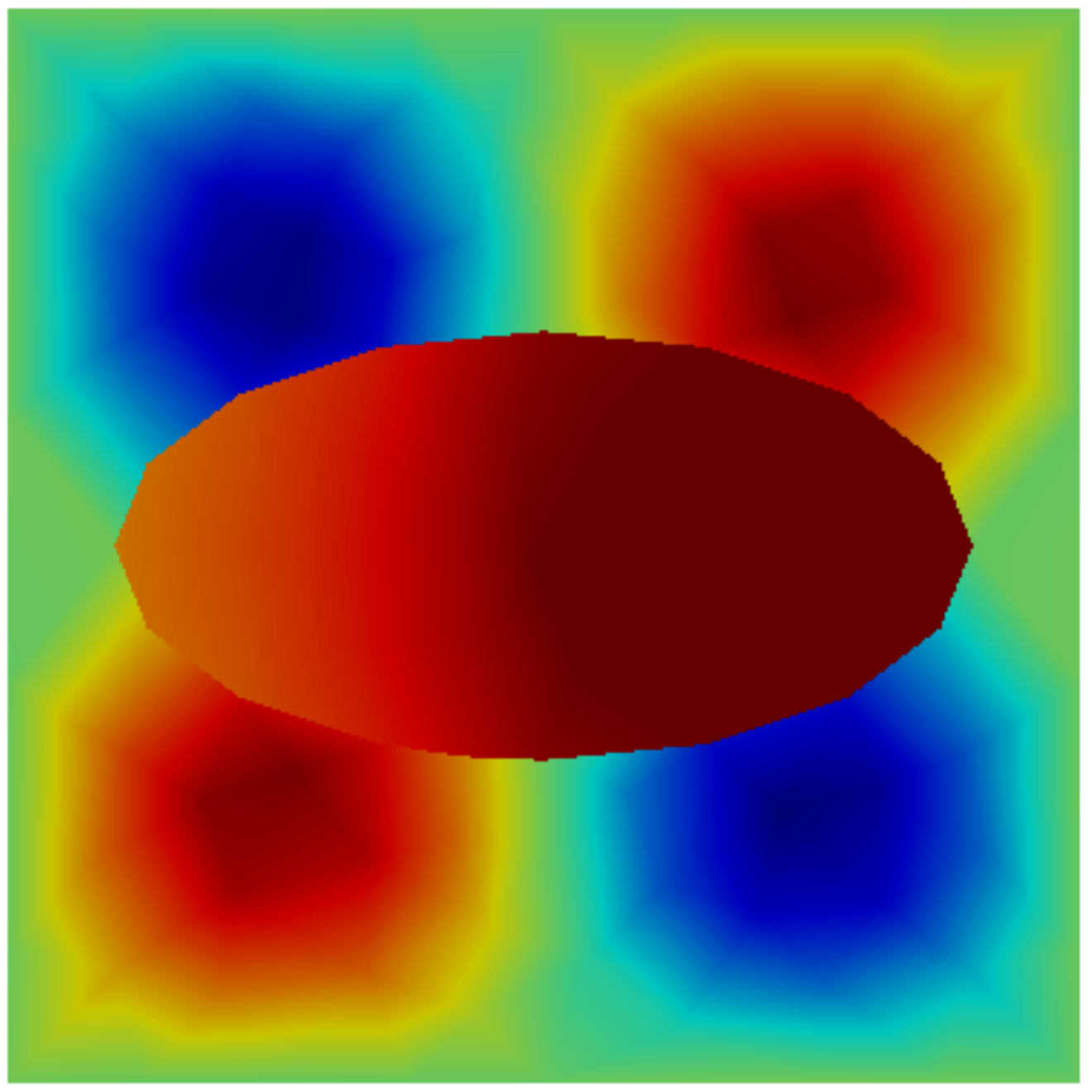}
\includegraphics[width=4.5cm]{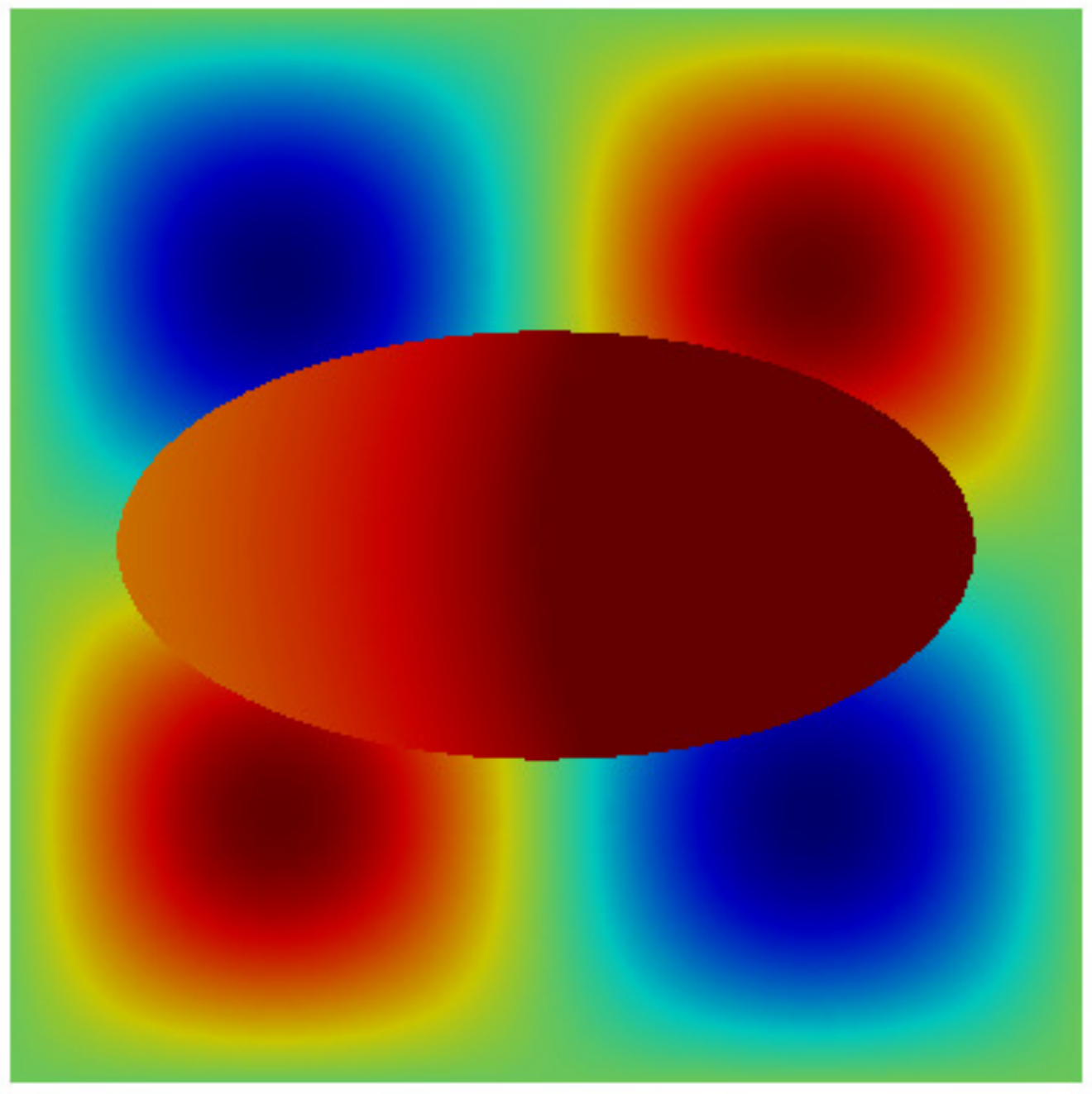}\\
\caption{Approximation of the scalar variable in Example \ref{ex6}. Columns: meshsize of $h=0.072$ and $0.018$. Rows: Polynomial of degree $k=0$, $1$ and $2$.}\label{fig:ex6}
\end{center}
\end{figure}

\begin{example}[Elliptical-shaped domain]\label{ex6}
{\rm We first solve a Poisson equation in a  the domain  $\Omega=(-1,1)^2$ divided by the elliptical interface $\Sigma$ described by $(x/0.8)^2+(y/0.4)^2 =1$. We take $\bld{\rm{K}}=\bld{\rm{I}}$ and \begin{eqnarray*}
u=\begin{cases}e^{x}\cos(y)\quad\hfill \textrm{in}\:\Omega^{1}\\ \sin(\pi x)\sin(\pi y)\quad\hfill\textrm{in}\:\Omega^{2}\end{cases}.
\end{eqnarray*}
as exact solution. The source term, transmission and Dirichlet boundary conditions are obtained from this exact solution.\\

 In Table~\ref{table:ex6} the history of convergence for this example is displayed. Similarly to the examples involving Neumann boundary data, the order of convergence for $u$ and $\qn$ are optimal whereas the convergence of the numerical trace is suboptimal, i.e., $O(h^{k+1})$. Moreover, even though superconvergence of the postprocessed solution $u^*_h$ is lost, it provides a more accurate approximation of $u$. Figure \ref{fig:ex6} shows the approximation $u_h$ obtained with  meshsizes of $h=0.072$ and $0.018$; and polynomial degree $k=0$, $1$ and $2$.}
\end{example}

\begin{figure}[H] 
\begin{center}
\includegraphics[width=4.5cm]{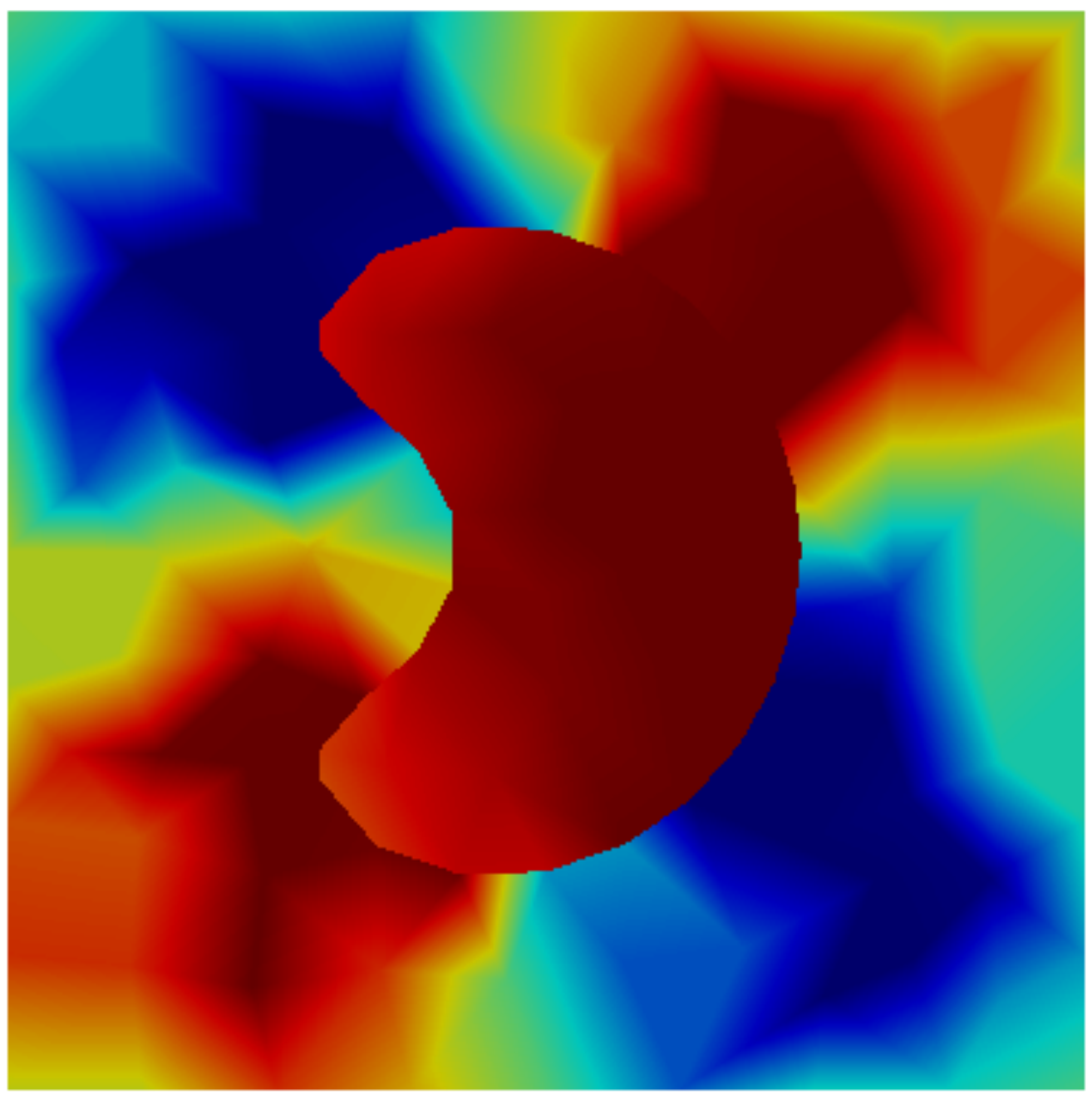}
\includegraphics[width=4.5cm]{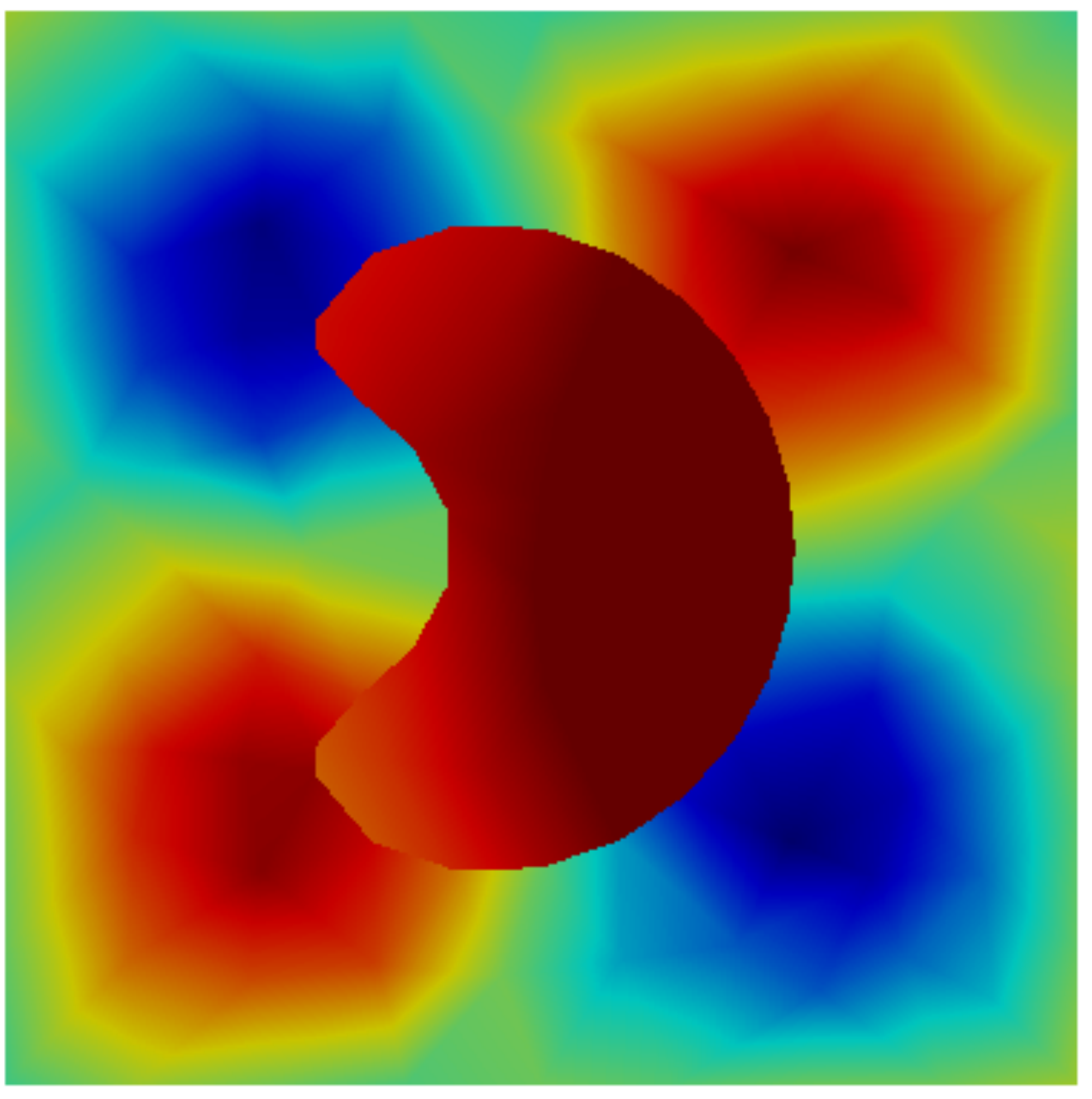}
\includegraphics[width=4.5cm]{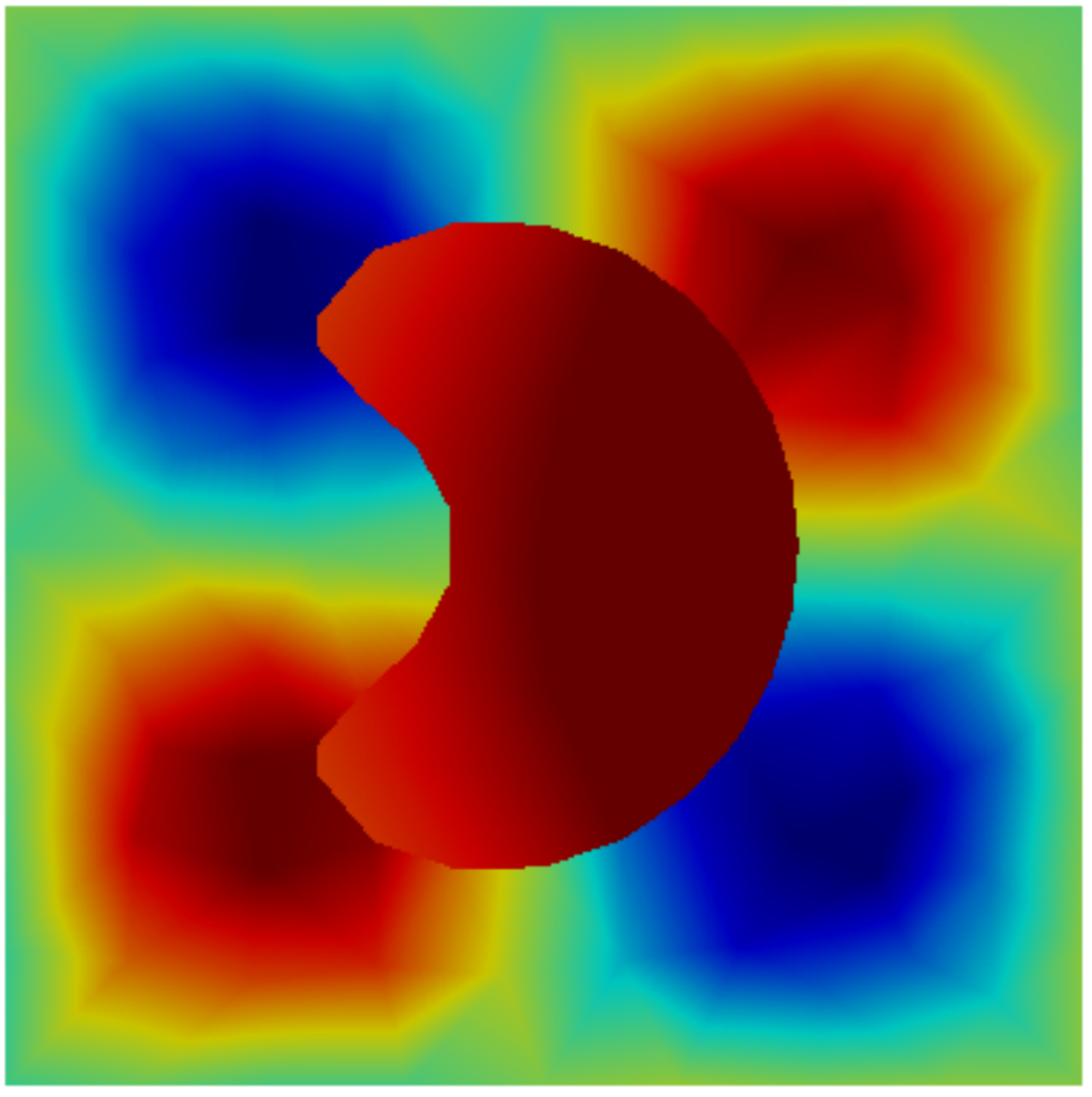}
\includegraphics[width=4.5cm]{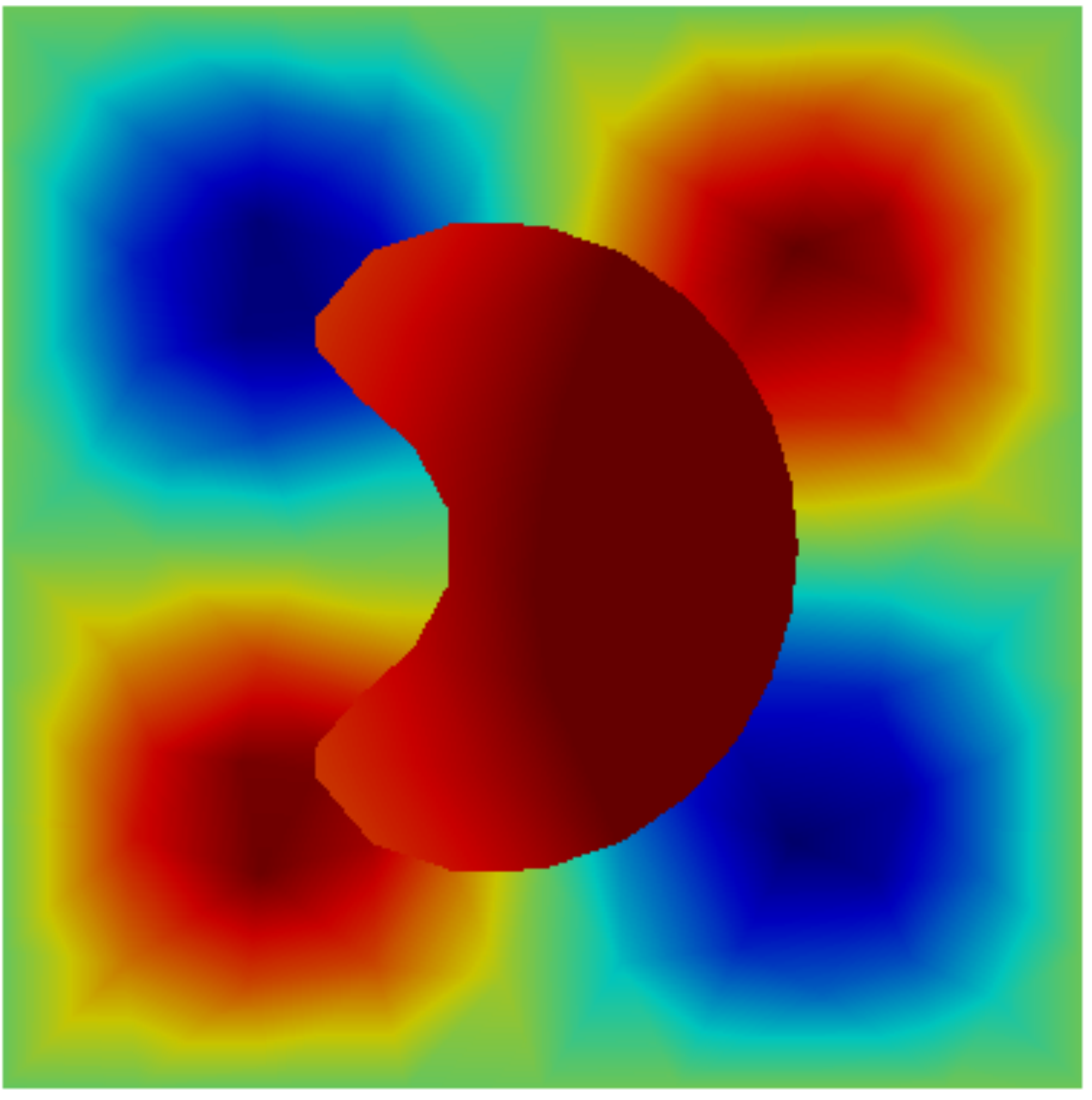}\\ 
\includegraphics[width=4.5cm]{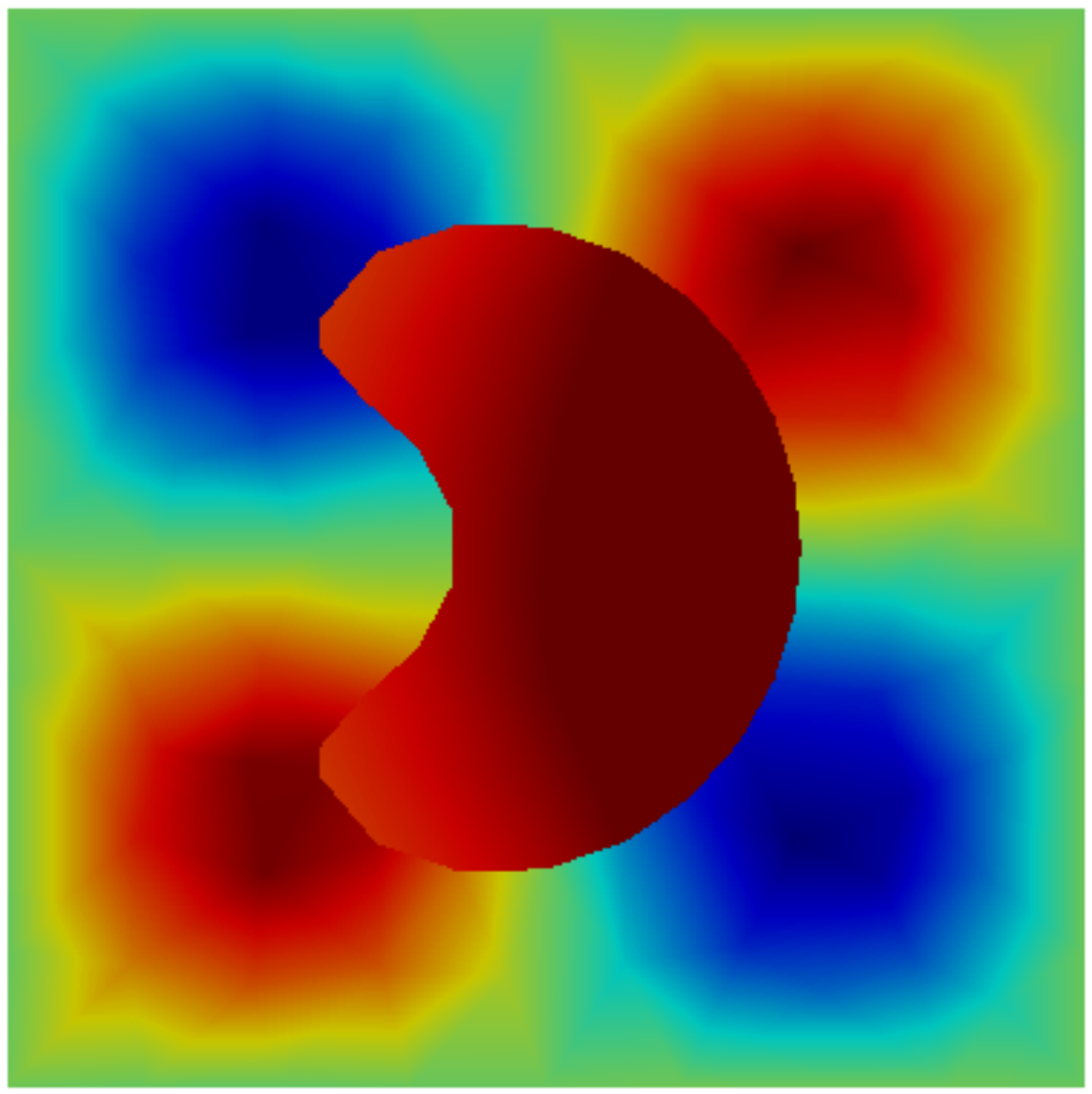}
\includegraphics[width=4.5cm]{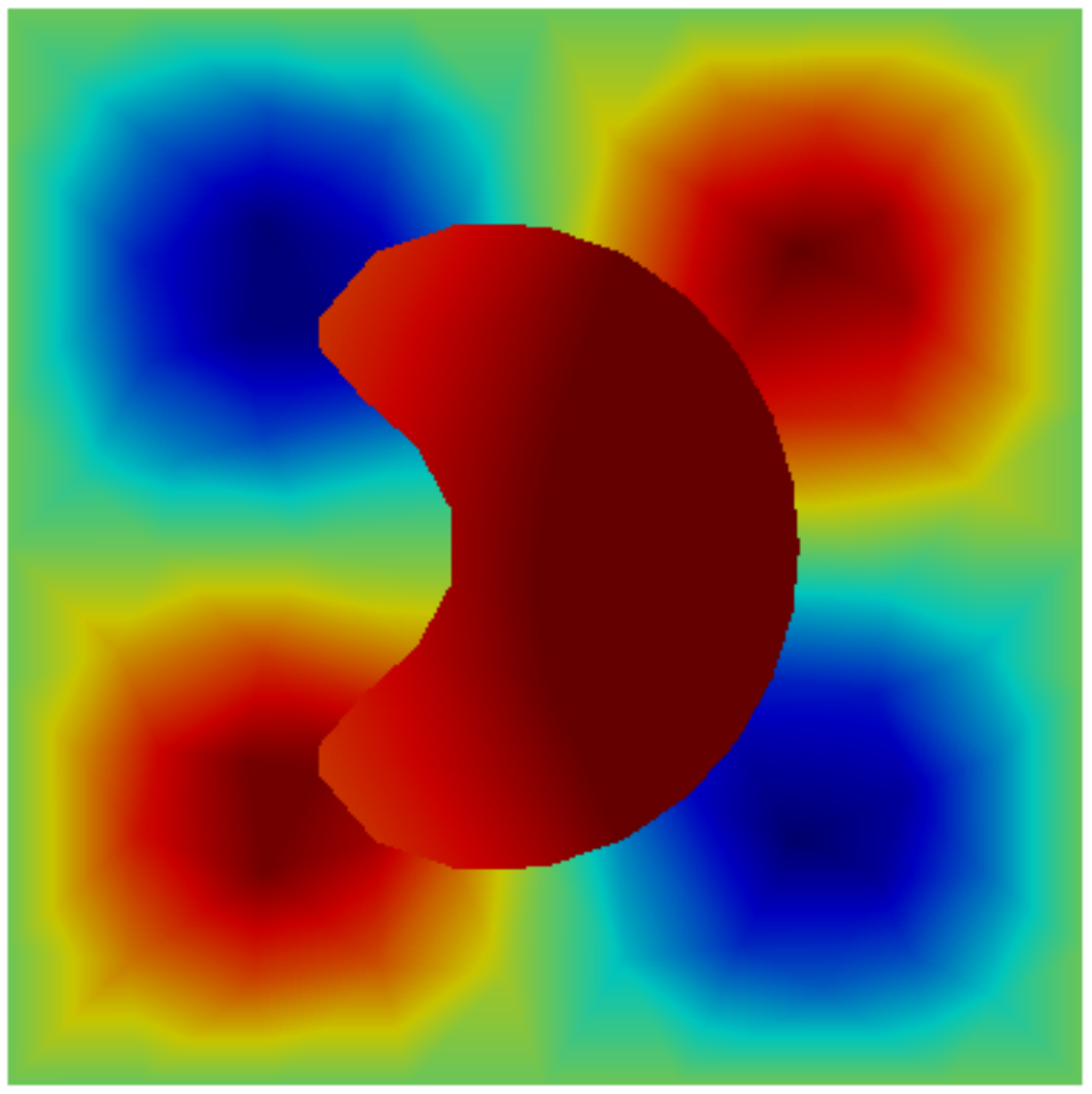}\\ 
\caption{Approximations $u_{h}$ (left) and  $u_{h}^{*}$ of the scalar variable $u$ of Example \ref{ex7}. Columns: meshsize $h=0.069$. Rows: Polynomial of degree $k=0$, $1$ and $2$.}\label{fig:ex7}
\end{center}
\end{figure}

\begin{example} [Kidney-shaped domain]\label{ex7}
{\rm We now consider the same exact solution as in previous example, but considering a kidney-shaped described by $(2[(x+0.5)^2+y^2]-x-0.5)^2-[(x+0.5)^2+y^2]+0.1=0$. In despite of the changes of convexity of this geometry, Table~\ref{table:ex7} shows similar accuracy on the approximations as the ones obtained in Example \ref{ex6}.  
Figure \ref{fig:ex7} shows the quality of the approximations of the scalar variable $u_h$ and its postprocessing $u^*_h$ obtained with  a meshsize  of $h=0.069$ and polynomial degree $k=0$, $1$ and $2$. As expected, $u^*_h$ provides a more accurate approximation of $u_h$ without significantly increase the computational cost.
}
\end{example}

\begin{table}[H]\renewcommand{\arraystretch}{1.3}\addtolength{\tabcolsep}{-4pt}
\begin{center}
{\scriptsize\begin{tabular}{c|c|cc|cc|cc|cc}
  \hline \hline
  \multicolumn{2}{c}{}     &
  \multicolumn{2}{c}{$\|e_u\|_{\ltwotilde}$}     &
  \multicolumn{2}{c}{$\|e_{\qn}\|_{\ltwotilde}$} &
  \multicolumn{2}{c}{$\|e_{\widehat{u}}\|_{\tildeE}$}&
  \multicolumn{2}{c}{$ \|e_{u^*}\|_{\ltwotilde}$}\\
  $k$& $h$& error & order &error&order &error &order &error &order \\
  \hline \hline
& $0.072$ & $2.37E-01$ & $-$ & $3.53E-01$ & $-$ & $3.66E-02$ & $-$ & $4.16E-02$ & $-$\\ 
& $0.035$ & $1.22E-01$ & $0.94$ & $1.92E-01$ & $0.87$ & $2.01E-02$ & $0.85$ & $2.16E-02$ & $0.93$\\ 
$0$ & $0.018$ & $5.97E-02$ & $1.03$ & $9.22E-02$ & $1.04$ & $1.05E-02$ & $0.93$ & $1.09E-02$ & $0.98$\\ 
& $0.009$ & $2.98E-02$ & $1.01$ & $4.66E-02$ & $0.99$ & $5.60E-03$ & $0.92$ & $5.67E-03$ & $0.95$\\ 
& $0.004$ & $1.50E-02$ & $1.00$ & $2.34E-02$ & $1.00$ & $2.84E-03$ & $0.98$ & $2.86E-03$ & $0.99$\\  \hline
& $0.072$ & $1.98E-02$ & $-$ & $4.20E-02$ & $-$ & $1.75E-03$ & $-$ & $2.24E-03$ & $-$\\ 
& $0.035$ & $4.95E-03$ & $1.97$ & $1.01E-02$ & $2.03$ & $1.63E-04$ & $3.37$ & $2.72E-04$ & $3.00$\\ 
$1$ & $0.018$ & $1.24E-03$ & $1.97$ & $2.36E-03$ & $2.07$ & $2.05E-05$ & $2.96$ & $3.06E-05$ & $3.12$\\ 
& $0.009$ & $3.12E-04$ & $2.01$ & $5.78E-04$ & $2.05$ & $7.79E-06$ & $1.41$ & $8.02E-06$ & $1.95$\\ 
& $0.004$ & $7.85E-05$ & $2.00$ & $1.43E-04$ & $2.02$ & $1.24E-06$ & $2.67$ & $1.22E-06$ & $2.73$\\  \hline
& $0.072$ & $1.44E-03$ & $-$ & $3.93E-03$ & $-$ & $1.25E-04$ & $-$ & $1.58E-04$ & $-$\\ 
& $0.035$ & $2.00E-04$ & $2.80$ & $5.24E-04$ & $2.86$ & $2.69E-05$ & $2.18$ & $2.78E-05$ & $2.47$\\ 
$2$ & $0.018$ & $2.43E-05$ & $3.01$ & $5.99E-05$ & $3.10$ & $1.89E-06$ & $3.80$ & $1.96E-06$ & $3.79$\\ 
& $0.009$ & $3.07E-06$ & $3.01$ & $7.54E-06$ & $3.01$ & $3.00E-07$ & $2.67$ & $3.03E-07$ & $2.72$\\ 
& $0.004$ & $3.92E-07$ & $2.98$ & $9.54E-07$ & $2.99$ & $3.99E-08$ & $2.92$ & $4.01E-08$ & $2.93$\\  \hline
& $0.072$ & $1.10E-04$ & $-$ & $3.02E-04$ & $-$ & $7.05E-06$ & $-$ & $8.09E-06$ & $-$\\ 
& $0.035$ & $7.78E-06$ & $3.76$ & $2.16E-05$ & $3.75$ & $1.79E-07$ & $5.21$ & $3.06E-07$ & $4.65$\\ 
$3$ & $0.018$ & $4.49E-07$ & $4.08$ & $1.17E-06$ & $4.16$ & $6.33E-09$ & $4.78$ & $8.78E-09$ & $5.08$\\ 
& $0.009$ & $2.82E-08$ & $4.02$ & $7.17E-08$ & $4.06$ & $4.60E-10$ & $3.81$ & $4.96E-10$ & $4.18$\\ 
& $0.004$ & $1.80E-09$ & $3.98$ & $4.49E-09$ & $4.01$ & $1.93E-11$ & $4.59$ & $2.02E-11$ & $4.63$\\  \hline\hline
\end{tabular}}
\end{center}
\caption{History of convergence of the approximation in Example \ref{ex6} (elliptical-shaped)}\label{table:ex6}
\end{table}

\begin{table}[h]\renewcommand{\arraystretch}{1.3} \addtolength{\tabcolsep}{-4pt}
\begin{center}
{\scriptsize\begin{tabular}{c|c|cc|cc|cc|cc}
  \hline \hline
  \multicolumn{2}{c}{}     &
  \multicolumn{2}{c}{$\|e_u\|_{\ltwotilde}$}     &
  \multicolumn{2}{c}{$\|e_{\qn}\|_{\ltwotilde}$} &
  \multicolumn{2}{c}{$\|e_{\widehat{u}}\|_{\tildeE}$}&
  \multicolumn{2}{c}{$ \|e_{u^*}\|_{\ltwotilde}$}\\
  $k$& $h$& error & order &error&order &error &order &error &order \\
  \hline \hline
& $0.069$ & $2.37E-01$ & $-$ & $3.76E-01$ & $-$ & $3.80E-02$ & $-$ & $4.39E-02$ & $-$\\ 
& $0.035$ & $1.23E-01$ & $0.97$ & $2.05E-01$ & $0.90$ & $1.97E-02$ & $0.98$ & $2.12E-02$ & $1.08$\\ 
$0$ & $0.018$ & $6.07E-02$ & $1.03$ & $9.73E-02$ & $1.09$ & $1.09E-02$ & $0.86$ & $1.13E-02$ & $0.92$\\ 
& $0.009$ & $3.01E-02$ & $1.01$ & $4.79E-02$ & $1.02$ & $5.69E-03$ & $0.94$ & $5.77E-03$ & $0.97$\\ 
& $0.004$ & $1.51E-02$ & $1.00$ & $2.41E-02$ & $1.00$ & $2.89E-03$ & $0.99$ & $2.91E-03$ & $1.00$\\   \hline
& $0.069$ & $2.13E-02$ & $-$ & $4.35E-02$ & $-$ & $1.89E-03$ & $-$ & $2.59E-03$ & $-$\\ 
& $0.035$ & $5.30E-03$ & $2.06$ & $1.08E-02$ & $2.05$ & $4.31E-04$ & $2.19$ & $4.94E-04$ & $2.45$\\ 
$1$ & $0.018$ & $1.33E-03$ & $2.03$ & $2.63E-03$ & $2.07$ & $1.03E-04$ & $2.10$ & $1.08E-04$ & $2.23$\\ 
& $0.009$ & $3.28E-04$ & $2.01$ & $6.24E-04$ & $2.07$ & $2.22E-05$ & $2.20$ & $2.27E-05$ & $2.24$\\ 
& $0.004$ & $8.28E-05$ & $2.00$ & $1.56E-04$ & $2.01$ & $5.74E-06$ & $1.97$ & $5.78E-06$ & $1.99$\\  \hline
& $0.069$ & $1.56E-03$ & $-$ & $3.72E-03$ & $-$ & $1.39E-04$ & $-$ & $1.74E-04$ & $-$\\ 
& $0.035$ & $2.02E-04$ & $3.02$ & $5.62E-04$ & $2.79$ & $1.78E-05$ & $3.04$ & $1.93E-05$ & $3.25$\\ 
$2$ & $0.018$ & $2.58E-05$ & $3.02$ & $6.53E-05$ & $3.15$ & $2.53E-06$ & $2.86$ & $2.60E-06$ & $2.93$\\ 
& $0.009$ & $3.19E-06$ & $3.01$ & $7.76E-06$ & $3.07$ & $4.17E-07$ & $2.60$ & $4.19E-07$ & $2.63$\\ 
& $0.004$ & $4.04E-07$ & $3.00$ & $9.80E-07$ & $3.01$ & $5.03E-08$ & $3.07$ & $5.04E-08$ & $3.08$\\   \hline
& $0.069$ & $1.31E-04$ & $-$ & $3.50E-04$ & $-$ & $1.27E-05$ & $-$ & $1.40E-05$ & $-$\\ 
& $0.035$ & $7.96E-06$ & $4.13$ & $2.11E-05$ & $4.15$ & $9.42E-07$ & $3.84$ & $9.68E-07$ & $3.94$\\ 
$3$ & $0.018$ & $4.92E-07$ & $4.08$ & $1.27E-06$ & $4.11$ & $4.03E-08$ & $4.61$ & $4.10E-08$ & $4.63$\\ 
& $0.009$ & $2.92E-08$ & $4.06$ & $7.37E-08$ & $4.10$ & $2.22E-09$ & $4.18$ & $2.23E-09$ & $4.19$\\ 
& $0.004$ & $1.87E-09$ & $3.99$ & $4.71E-09$ & $4.00$ & $1.44E-10$ & $3.98$ & $1.44E-10$ & $3.98$\\ \hline\hline
\end{tabular}}
\end{center}
\caption{History of convergence of the approximation in Example \ref{ex7} (kidney-shaped)}\label{table:ex7}
\end{table}

%

\begin{figure}[H]
\begin{center}
\includegraphics[width=4.5cm]{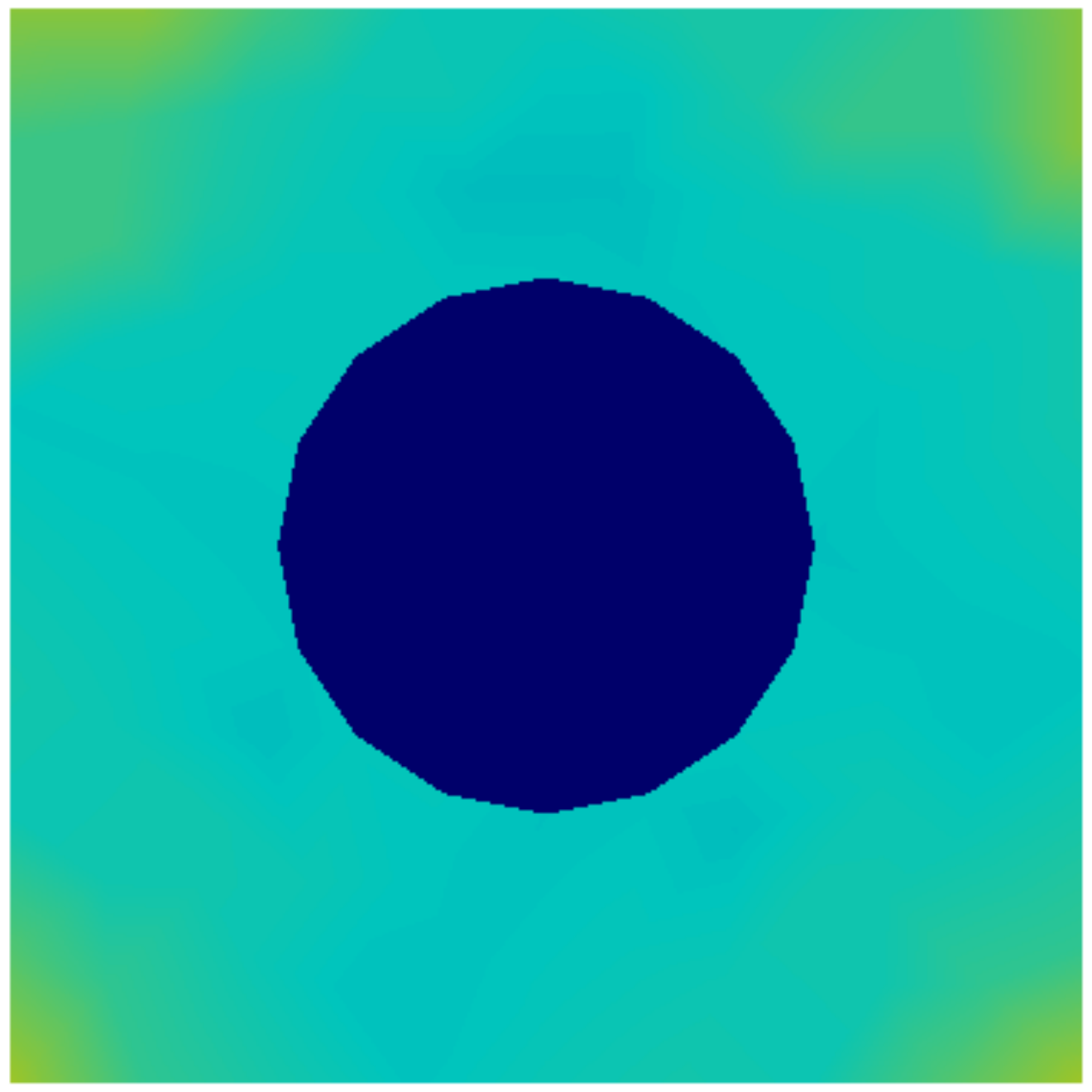}
\includegraphics[width=4.5cm]{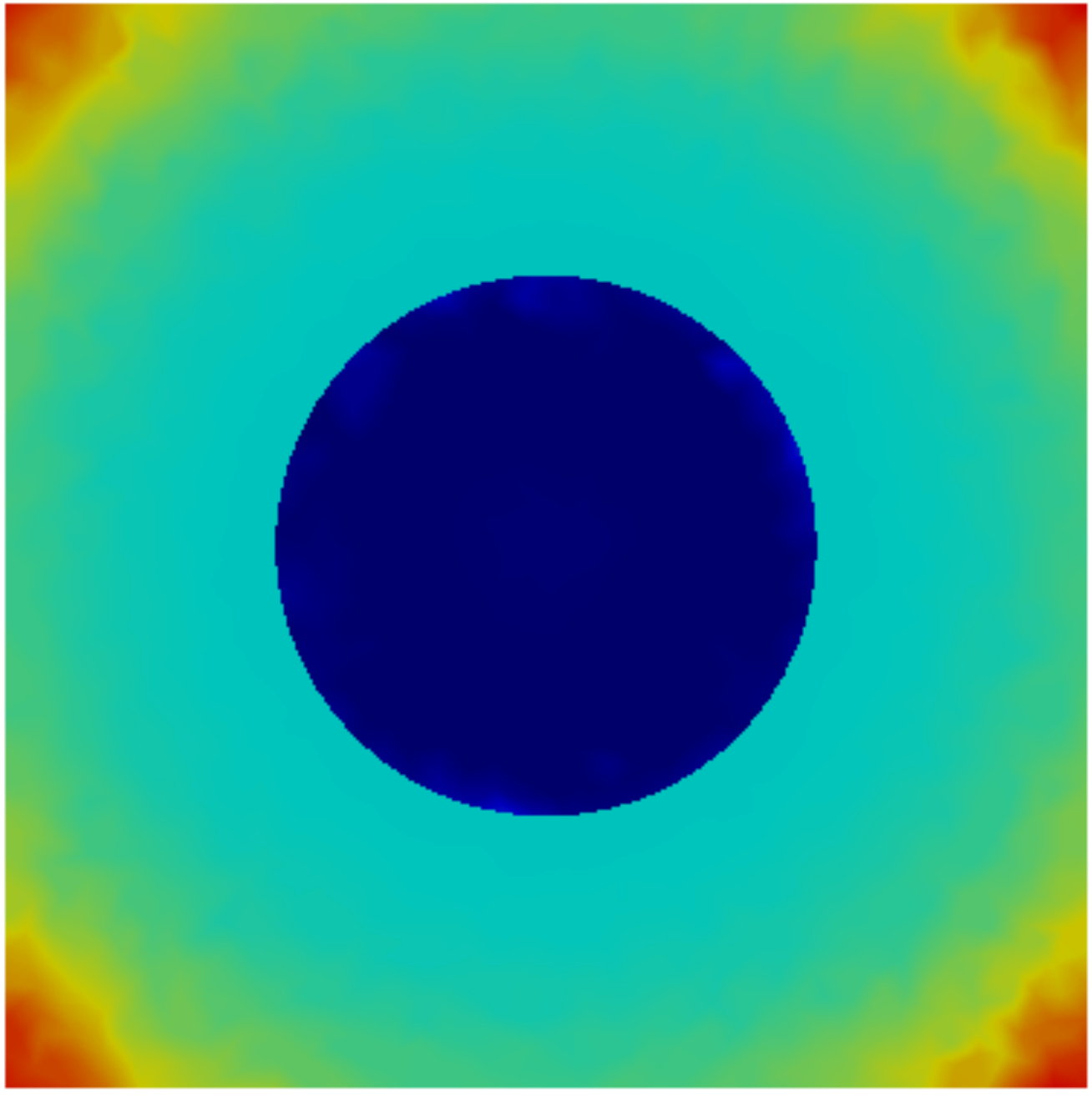}\\ 
\includegraphics[width=4.5cm]{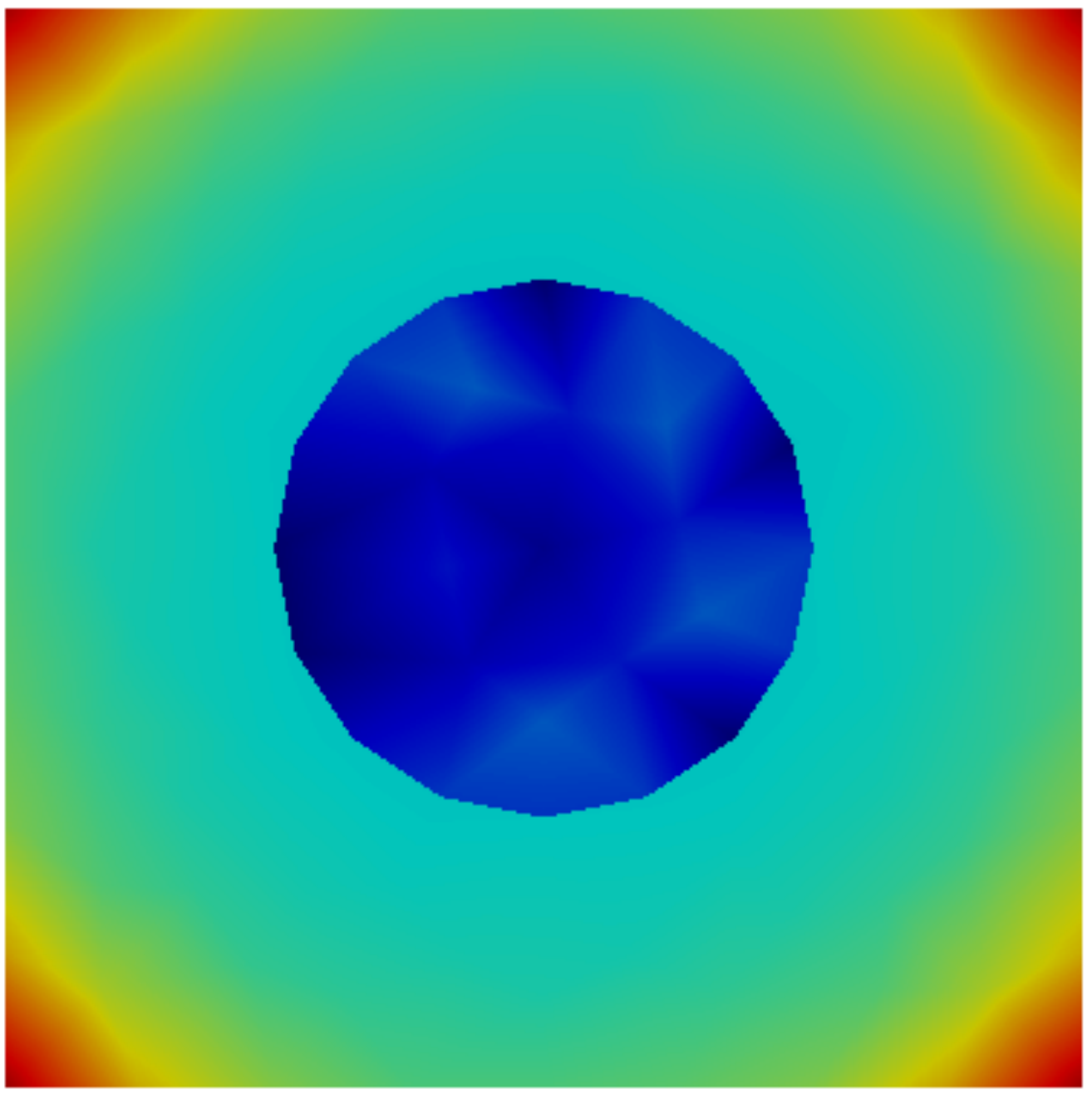}
\includegraphics[width=4.5cm]{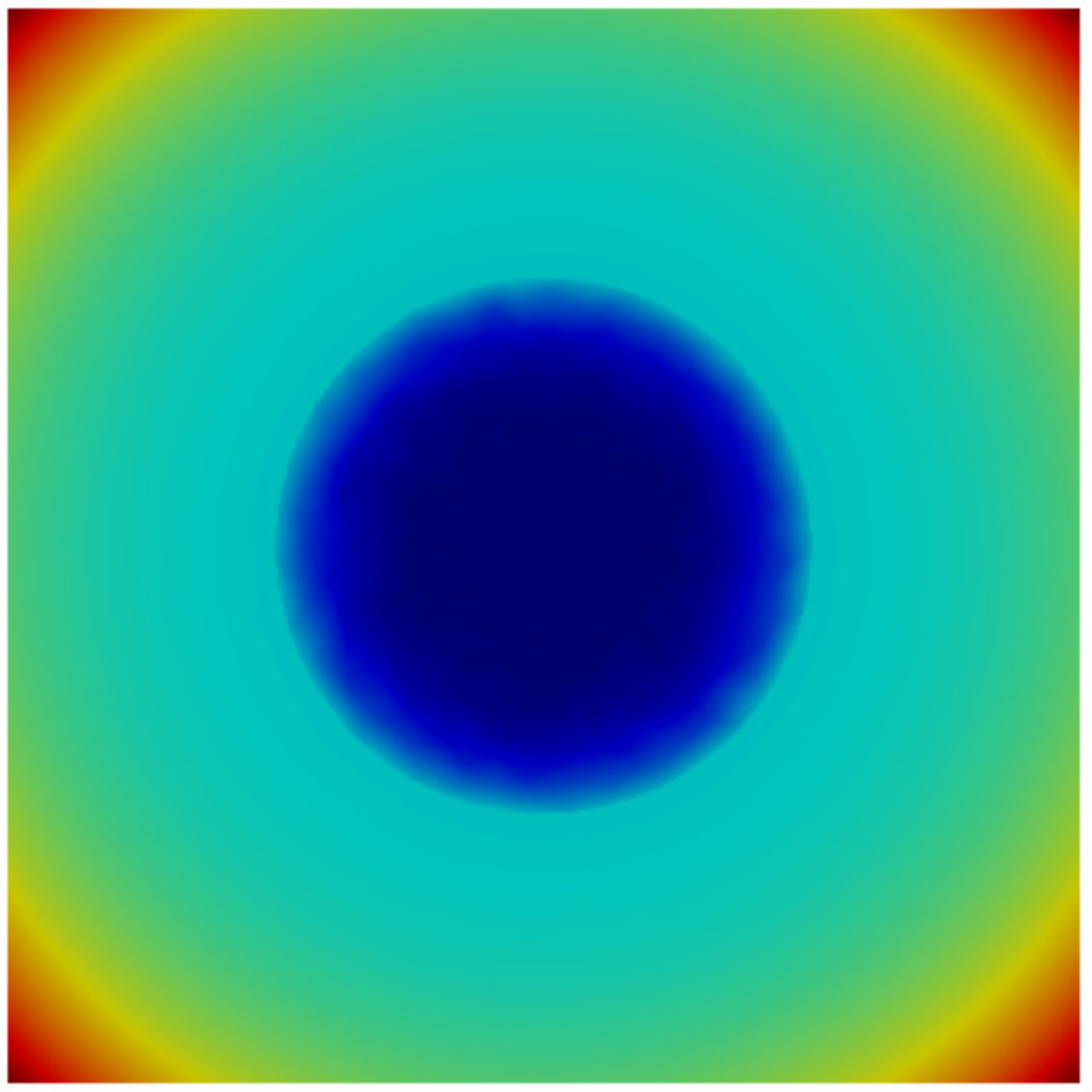}\\ 
\includegraphics[width=4.5cm]{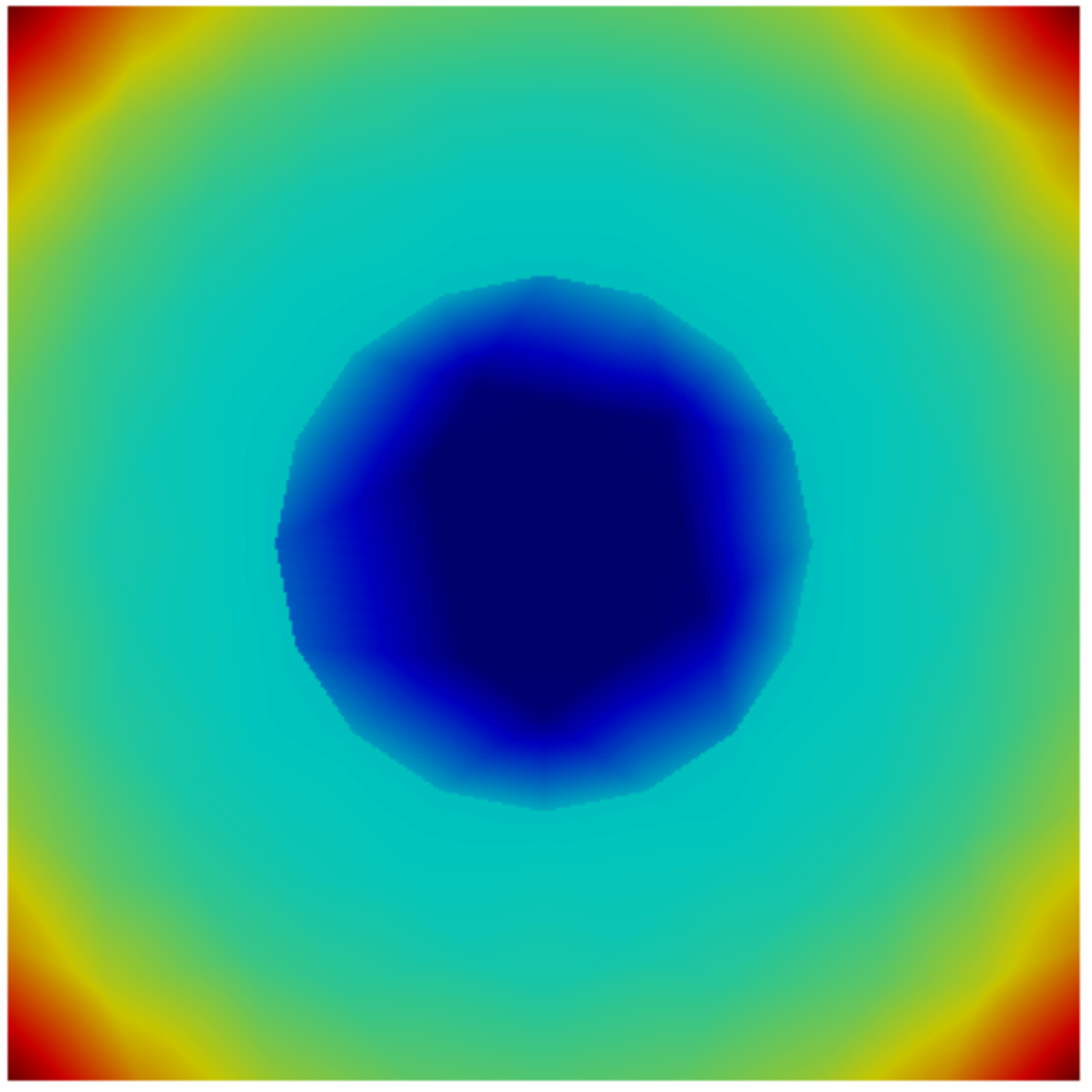}
\includegraphics[width=4.5cm]{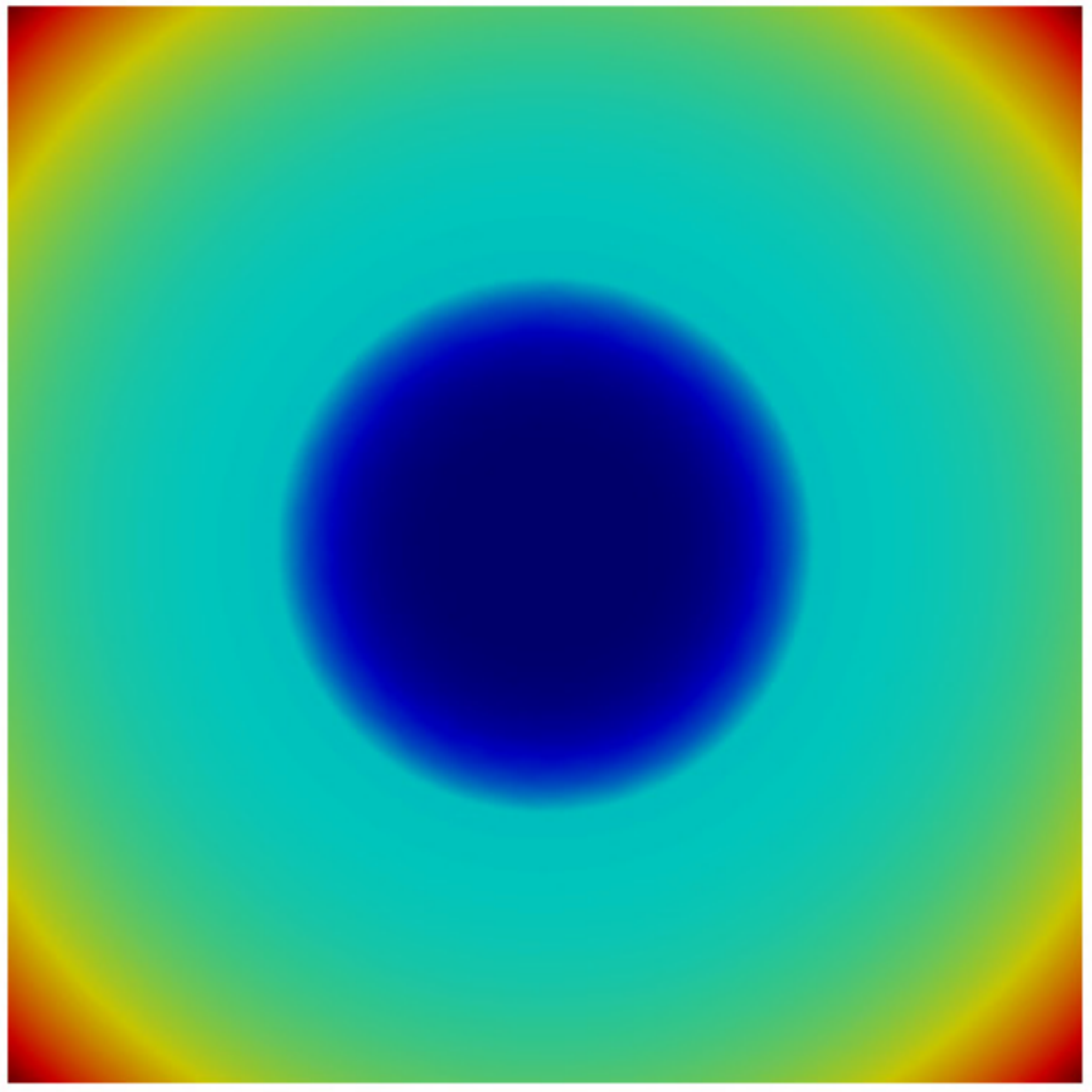}\\ 
\caption{Approximation of the scalar variable in Example \ref{ex8} (thermal conductivity). Columns: meshsize $h=0.072$ and $0.018$. Rows: Polynomial of degree $k=0$, $1$ and $2$.}\label{fig:ex8}
\end{center}
\end{figure}

\begin{example}[Thermal conductivity]\label{ex8}
{\rm Finally, considering the example provided by \cite{MIT1}, we simulate the heat distribution $u$ at steady state, due to the heat source $f$, over the domain $\Omega=(-1,1)^2$ divided by a circular interface of radius $R=0.5$ centered at the origin. The source term $f$ and the thermal conductivity tensor are given by
\begin{eqnarray*}
f(x,y)=-10(x^{2}+y^{2})^{3/2}-15x^{2}(x^{2}+y^{2})^{1/2}-15y^{2}(x^{2}+y^{2})^{1/2}
\end{eqnarray*}
and $\bld{\rm{K}}=\kappa_{j} \In$ in $\Omega^{j}$ ($j=1,2)$. The exact solution of this problem is
\begin{eqnarray*}
u=\begin{cases}\frac{1}{\kappa_{1}}(x^{2}+y^{2})^{5/2}\quad\hfill\text{ in }\in\Omega^{1}\\ \frac{1}{\kappa_{2}}(x^{2}+y^{2})^{5/2}+\left(\frac{1}{\kappa_{1}}-\frac{1}{\kappa_{2}}\right)R^{5}\quad\hfill\text{ in }\Omega^{2}\end{cases},
\end{eqnarray*}
and we consider $\kappa_{1}=1$, $\kappa_{2}=100$. Dirichlet boundary condition on $\Gamma$ is derived from the previous equation. In this case the jumps  $s_D$ and $s_N$ are both equal to zero. Table~\ref{table:ex8} validates the optimal convergence rates of order $h^{k+1}$ for the heat distribution $u$ and the flux $\qn$. Figure \ref{fig:ex8} shows the approximated heat distribution considering meshes of size $h=0.072$ and $0.018$, and polynomials of degree $k=0$, $1$ and $2$.
}
\end{example}


\begin{remark}
If the mesh is fine enough, the errors $e_u$, $e_{\qn}$ and $e_{u^*}$ can be computed in the entire computational domain $\textsf{D}_{h}$ since the quadrature points of a triangle $K\in \textsf{D}_{h}^j$ will eventually lie in $\Omega^j$. This happens in all previous examples. In fact, we computed the errors $\|e_u\|_{L^2(\textsf{D}_{h})}$, $\|e_{\qn}\|_{L^2(\textsf{D}_{h})}$ and  $\|e_{u*}\|_{L^2(\textsf{D}_{h})}$. Their behavior and magnitude are similar to ones displayed in the convergence tables.
\end{remark}

\section{Conclusions}\label{sec:conclusions}

We have proposed a technique for high order approximation of boundary value problems in curved domains with mixed boundary conditions.  We have provided numerical evidence suggesting that the technique performs properly if the family of paths is normal to the computational boundary. A practical way to always satisfy this restriction is to define  $\Gamma^h$ by  interpolating $\Gamma$ using only piecewise linear segments. Moreover, we have extend this technique to elliptic interface problems where the interface is not necessarily polygonal. We have presented numerical results indicating that the order of convergence of  are optimal for the error of $u$ and $\qn$ if the interface is interpolated by  piecewise linear segments.


\section*{Acknowledgments}
W. Qiu is partially by the GRF of Hong Kong (Grant No. 9041980 and 9042081) and a grant from the Research Grants Council of the Hong Kong Special Administrative Region, China (Project No. CityU 11302014). 
M. Solano Partially supported by CONICYT-Chile
through grant FONDECYT-11130350, BASAL project CMM, Universidad de Chile and Centro de Investigaci\'on en Ingenier\'ia  Matem\'atica (CI$^2$MA). P. Vega acknowledges the Scholarship Program of CONICYT-Chile.

\begin{table}[ht!]\renewcommand{\arraystretch}{1.3}\addtolength{\tabcolsep}{-4pt}
\begin{center}
{\scriptsize\begin{tabular}{c|c|cc|cc|cc|cc}
  \hline \hline
  \multicolumn{2}{c}{}     &
  \multicolumn{2}{c}{$\|e_u\|_{\ltwotilde}$}     &
  \multicolumn{2}{c}{$\|e_{\qn}\|_{\ltwotilde}$} &
  \multicolumn{2}{c}{$\|e_{\widehat{u}}\|_{\tildeE}$}&
  \multicolumn{2}{c}{$ \|e_{u^*}\|_{\ltwotilde}$}\\
  $k$& $h$& error & order &error&order &error &order &error &order \\
  \hline \hline
& $0.072$ & $9.10E-03$ & $-$ & $5.66E-01$ & $-$ & $1.42E-03$ & $-$ & $1.63E-03$ & $-$\\ 
& $0.035$ & $7.30E-03$ & $0.31$ & $3.09E-01$ & $0.84$ & $7.58E-04$ & $0.88$ & $8.22E-04$ & $0.96$\\ 
$0$ & $0.018$ & $4.56E-03$ & $0.68$ & $1.46E-01$ & $1.08$ & $4.04E-04$ & $0.91$ & $4.21E-04$ & $0.97$\\ 
& $0.009$ & $2.50E-03$ & $0.87$ & $7.32E-02$ & $1.01$ & $2.09E-04$ & $0.96$ & $2.13E-04$ & $0.99$\\ 
& $0.004$ & $1.28E-03$ & $0.96$ & $3.60E-02$ & $1.02$ & $1.06E-04$ & $0.97$ & $1.07E-04$ & $0.99$\\   \hline
& $0.072$ & $1.39E-03$ & $-$ & $5.99E-02$ & $-$ & $5.95E-05$ & $-$ & $1.55E-04$ & $-$\\ 
& $0.035$ & $4.51E-04$ & $1.57$ & $1.45E-02$ & $1.98$ & $1.23E-05$ & $2.20$ & $1.87E-05$ & $2.95$\\ 
$1$ & $0.018$ & $1.36E-04$ & $1.74$ & $3.36E-03$ & $2.12$ & $2.34E-06$ & $2.39$ & $3.26E-06$ & $2.52$\\ 
& $0.009$ & $3.73E-05$ & $1.88$ & $8.24E-04$ & $2.04$ & $5.18E-07$ & $2.19$ & $6.20E-07$ & $2.41$\\ 
& $0.004$ & $9.42E-06$ & $1.98$ & $2.06E-04$ & $1.99$ & $6.66E-08$ & $2.95$ & $8.03E-08$ & $2.94$\\    \hline
& $0.072$ & $1.69E-04$ & $-$ & $3.94E-03$ & $-$ & $1.28E-05$ & $-$ & $2.07E-05$ & $-$\\ 
& $0.035$ & $2.33E-05$ & $2.77$ & $4.62E-04$ & $3.00$ & $9.34E-07$ & $3.66$ & $1.19E-06$ & $3.99$\\ 
$2$ & $0.018$ & $3.40E-06$ & $2.78$ & $5.36E-05$ & $3.11$ & $1.05E-07$ & $3.16$ & $1.19E-07$ & $3.32$\\ 
& $0.009$ & $4.73E-07$ & $2.87$ & $6.64E-06$ & $3.04$ & $1.38E-08$ & $2.95$ & $1.45E-08$ & $3.07$\\ 
& $0.004$ & $5.94E-08$ & $2.98$ & $7.84E-07$ & $3.07$ & $1.39E-09$ & $3.30$ & $1.43E-09$ & $3.33$\\    \hline
& $0.072$ & $1.35E-05$ & $-$ & $1.58E-04$ & $-$ & $8.50E-07$ & $-$ & $1.24E-06$ & $-$\\ 
& $0.035$ & $8.36E-07$ & $3.89$ & $7.36E-06$ & $4.28$ & $1.80E-08$ & $5.39$ & $2.51E-08$ & $5.45$\\ 
$3$ & $0.018$ & $5.64E-08$ & $3.90$ & $4.05E-07$ & $4.19$ & $1.03E-09$ & $4.13$ & $1.23E-09$ & $4.37$\\ 
& $0.009$ & $3.94E-09$ & $3.87$ & $2.43E-08$ & $4.09$ & $6.98E-11$ & $3.92$ & $7.39E-11$ & $4.08$\\ 
& $0.004$ & $2.49E-10$ & $3.97$ & $1.41E-09$ & $4.10$ & $2.02E-12$ & $5.10$ & $2.17E-12$ & $5.08$\\    \hline\hline
\end{tabular}}
\end{center}
\caption{History of convergence of the approximation in Example \ref{ex8} (thermal conductivity)}
\label{table:ex8}
\end{table}


\end{document}